\documentclass[a4paper,11pt]{article}
\usepackage[UKenglish]{babel}
\usepackage{authblk}
\usepackage[colorlinks=true]{hyperref}
\hypersetup{citecolor=blue}
\usepackage{color}
\usepackage[utf8]{inputenc}
\usepackage{amsfonts,amsthm,amssymb}
\usepackage{mathtools}
\usepackage{dsfont}
\usepackage{float}
\usepackage{epsfig}
\usepackage{graphicx}
\usepackage{caption}
\usepackage{subcaption}
\usepackage{kotex}
\usepackage{tikz}
\usepackage{pgfplots}
\pgfplotsset{
    tick align=outside,
    x grid style={white},
    xmajorgrids,
    y grid style={white},
    ymajorgrids,
    axis line style={white},
    axis background/.style={fill=white!92!black},
    legend style={draw=white, fill=white},
    legend cell align={left}
}
\newtheorem{remark}{Remark}
\newtheorem{thm}{Theorem}[section]
\newtheorem{lem}[thm]{Lemma}

\numberwithin{equation}{section}
\numberwithin{thm}{section}
\numberwithin{problem}{section}
\newcommand{\beq}{\begin{equation}}
\newcommand{\eeq}{\end{equation}}
\newcommand {\f}{\frac}
\newcommand {\pa}{\partial}
\newcommand {\e}{\varepsilon}

\newcommand {\init}{\text{in}}
\newcommand {\MR}{\text{MR}}
\newcommand {\ADN}{\text{ADN}}
\newcommand {\GN}{\text{GN}}
\newcommand{\R}{\mathbb{R}}

\newcommand\N{{\mathbb N}}

\newcommand\caFu{{\cal F}_U}
\newcommand\caFv{{\cal F}_V}
\newcommand\caFw{{\cal F}_W}
\newcommand\caK{{\cal K}}
\newcommand\caP{{\cal P}}
\newcommand\caL{{\cal L}}
\newcommand\caQ{{\cal Q}}

\DeclareMathOperator{\In}{\,in\,}

\newcommand\Ep{{\mathcal E}_p}
\newcommand\afast{a_{\text{fast}}}
\newcommand\bfast{b_{\text{fast}}}
\newcommand\cfast{c_{\text{fast}}}
\newcommand\dfast{d_{\text{fast}}}
\newcommand\MScN[1]{\href{http://www.ams.org/mathscinet-getitem?mr=#1}{\nolinkurl{(#1)}}}
\newcommand\DOI[1]{\href{http://dx.doi.org/#1}{(doi: \nolinkurl{#1})}}
\newcommand\LINK[1]{\href{#1}{(link: \nolinkurl{#1})}}
\newenvironment{acknowledgment}{\noindent{\bf Acknowledgment}}{}
\title
{On a class of triangular cross-diffusion systems and its fast reaction approximation}

\author{E. Brocchieri$^1$, L. Corrias$^2$}
\date{\today}
\providecommand{\keywords}[1]{\small\textit{{Keywords.}} #1\\}
\providecommand{\subjclass}[1]{\small\textit{{2010 Mathematics Subject Classification.}} #1}
\begin{document}
\maketitle
\begin{abstract}
The purpose of this article is to investigate the emergence of cross-diffusion in the time evolution of two slow-fast species in competition.  A class of  triangular cross-diffusion system is obtained as the singular limit of a fast reaction-diffusion system.  We first prove the convergence of the unique strict solution of the fast reaction-diffusion system towards a (weak, strong) solution of the cross-diffusion system, as the reaction rate $\e^{-1}$ goes to $+\infty$. Furthermore, under the assumption of small cross-diffusion, we obtain a convergence rate  as well as the influence of the initial layer, due to initial data, on the convergence rate itself. Both results are obtained through energy functionals that handle the fast reaction terms uniformly in~$\e$. 
\end{abstract}
\keywords{Cross-diffusion, singular limits, dynamical systems, slow-fast manifold.}
\subjclass{Primary : 35B25, 35B40, 35K57, 35Q92, 92D25. Secondary 35B45, 35K45}. 
\section{Introduction}
This article deals with the emergence of cross-diffusion in the singular limit of a reaction-diffusion system with multiple time scales. The system models two species, say ${\bf u}$ and ${\bf v}$, in competition in a bounded region of $\R^N$, $N\!\ge\!1$, with reflecting boundary. Due to the inter et intra competition, the individuals of the species ${\bf u}$ may switch between two different states, ${\bf a}$ and ${\bf b}$, with switching rate of order $\e^{-1}$, $\e\!>\!0$, and change the diffusivity together with their state. Hence, the density $u^\e$ of the population ${\bf u}$ writes as $u^\e\!=\!u^\e_a+u^\e_b$, where $u^\e_a\!\ge\!0$, $u^\e_b\!\ge\!0$ are the densities of the two subpopulations (one for each state). The slow population~${\bf v}$ has density $v^\e$ and the global dynamic is modelled by the fast reaction system below
\begin{equation}\label{meso system}
\begin{cases}
\partial_tu_a^\e-d_a\Delta u_a^\e=f_a(u_a^\e,u_b^\e,v^\e)+\e^{-1} Q(u_a^\e,\,u_b^\e,\, v^\e\,),&\In\,\,(0,\infty)\times\Omega,\\
\partial_tu_b^\e-d_b\Delta u_b^\e=f_b(u_a^\e,u_b^\e,\,v^\e)-\e^{-1} Q(u_a^\e,\, u_b^\e,\,v^\e\,),&\In\,\,(0,\infty)\times\Omega,\\
\partial_tv^\e-d_v\Delta v^\e=f_v(u_a^\e,u_b^\e,\,v^\e),&\In\,\,(0,\infty)\times\Omega,
\end{cases}
\end{equation}
where $d_a,d_b,d_v\!>\!0$, $d_a\neq d_b$ and $\Omega$ is a smooth bounded domain of $\mathbb{R}^N$, $N\!\ge\!1$. 

System  \eqref{meso system} is supplemented with homogeneous Neumann boundary conditions 
\beq\label{HNBC}
\nabla u_a^\e\cdot \vec{n}=\nabla u_b^\e\cdot \vec{n}=\nabla v^\e\cdot \vec{n}=0\,,\qquad\In\,(0,\infty)\times\pa\Omega,
\eeq
where $\vec{n}$ is the unit outward normal vector, and nonnegative initial data
\begin{equation}\label{ID ua ub v}
u_a^\e(0)=u_a^\init\ge0,\quad
u_b^\e(0)=u_b^\init\ge0,\quad 
v^\e(0)=v^\init\ge0,\quad\In\,\Omega.
\end{equation}
The competitive dynamics are given by the reaction functions 
\beq\label{def react functions}
\begin{split}
f_a(u_a,u_b,v)&:= \eta_au_a(1-a u_a-cv)-\gamma_au_au_b,\\
f_b(u_a,u_b,v)&:= \eta_bu_b(1-b u_b-dv)-\gamma_bu_au_b,\\
f_v(u_a,u_b,v)&:= \eta'_vv(1-au_a-cv)+\eta''_vv(1-bu_b-dv)\,,
\end{split}
\eeq
with $a,b>0$, $c,d\in\R_+$, $\eta_a,\eta_b>0$, $\eta'_v,\eta''_v,\in\R_+$, $(c\eta'_v,d\eta''_v)\neq(0,0)$, $\gamma_a,\gamma_b\in\R_+$, while the switching dynamic between the states ${\bf a}$ and ${\bf b}$ is modelled  by
\begin{equation}\label{def Q}
\begin{split}
&Q(u_a,u_b,v):=q(u_a,u_b,v)/\Lambda(u_a,u_b,v)\,,\\
&q(u_a,u_b,v):=\phi(\bfast u_b+\dfast v)\, u_b -\psi(\afast u_a+\cfast v)\,u_a\,,\\
&\Lambda(u_a,u_b,v):=\phi(\bfast u_b+\dfast v)+\psi(\afast u_a+\cfast v)\,,
\end{split}
\eeq
where $\psi,\phi$ are chosen so that $\Lambda>0$ and $\afast,\bfast,\cfast,\dfast\in\R_+$, $(\afast,\cfast)\neq(0,0)$, $(\bfast,\dfast)\neq(0,0)$. 

We are interested in the case where the increase in density $u_b^\e$ (respectively $u_a^\e$) pushes individuals of species ${\bf u}$ to migrate towards the state ${\bf a}$ (respectively ${\bf b}$). Therefore, it is natural to consider increasing transition functions $\psi,\phi\,:\,\R_+\mapsto\R_+$. In addition, this choice gives us relative satisfaction measures 
\beq\label{def rel sat measures}
\begin{split}
&\Lambda_a(u_a,u_b,v):=\psi(\afast u_a+\cfast v)/\Lambda(u_a,u_b,v)\,,\\
&\Lambda_b(u_a,u_b,v):=\phi(\bfast u_b+\dfast v)/\Lambda(u_a,u_b,v)\,,
\end{split}
\eeq
increasing with respect to $u_a$ and $u_b$, respectively, as well as a family of energy functionals well fitted to handle the fast reaction term $\e^{-1}Q$.  Furthermore, as power laws are meaningful from the biological point of view, and general enough to allow to extend our results to any pair $(\psi,\phi)$ behaving like power functions as $x\searrow0$ and $x\nearrow\infty$, we choose 
\beq\label{def phi psi}
\psi(x)=(A+ x)^{\alpha},\quad\phi(x)=(B+x)^{\beta},\quad x\ge0\,,\quad \alpha,\beta>0\,.
\eeq

Due to the symmetry of \eqref{meso system}--\eqref{def phi psi}, i.e. to the interchangeable role of the subpopulations with densities $u_a^\e$ and $u_b^\e$, we can assume without loss of generality that $\alpha\le\beta$. As the transition rates $\alpha, \beta$ can be different from each other, we have an additional slow-fast intra-dynamic that deeply affects the mathematical analysis of the system. Hence, in order to have well defined relative satisfaction measures \!\eqref{def rel sat measures} and well defined and manageable energy functionals, we need to assume the following
\begin{align}
\tag{H1}\label{(H1)}
0<\alpha\le\beta\,,\quad& A>0\,,\quad B\ge0\quad\text{and}\quad B>0\quad\text{if}\quad\beta<1\,,\\
\tag{H2}\label{(H2)}
&0\le\beta-\alpha< 2(\alpha+3)\,,\\
\tag{H3}\label{(H3)}
\afast\le a\,,\qquad &\bfast\le b\,,\qquad \cfast\le c\,,\qquad \dfast\le d\,.
\end{align}
Note that assumption \eqref{(H2)} is written in a way that highlights the difference between the transition rates, since this gap will be crucial in the further analysis. Obviously, if $\beta\le\alpha$, one has to switch the role of $A$ and $B$ and the role of $\alpha$ and $\beta$ in \eqref{(H1)}, \eqref{(H2)}.  Furthermore, assumption \eqref{(H3)} beyond being relevant from the mathematical point of view, it is also biologically meaningful since it implies that the switching between the two states of the individuals of the species ${\bf u}$ cannot occur too fast.

The  fast reaction-diffusion system system \eqref{meso system}--\eqref{def phi psi} is the natural generalisation of the system introduced in \cite{Brocchieri2021} to investigate the impact of dietary diversity of populations in competition. However, in  \cite{Brocchieri2021} only the subpopulation with density $u_b^\e$ is in direct competition with the population ${\bf v}$ and no intra-specific competition in the population ${\bf u}$ is taken into account, i.e. $c=\gamma_a=\gamma_b=\eta_v'=0$,  $b=d$, $\eta_v''=\eta_b$, so that $u_b^{-1}f_b=v^{-1}f_v$. Moreover, ${\bf u}$ is in direct competition with the species ${\bf v}$, i.e. $a=b$, $c=d$, $\gamma_a=\gamma_b$, $\eta_a=\eta_b$ and $\eta_v'=\eta_v''$, so that $u_a^{-1}f_a=u_b^{-1}f_b$. Finally, in \cite{Brocchieri2021}, $\afast$ and $\bfast$ are zero, meaning that the fast switching between individuals of population ${\bf u}$ is uniquely determined by ${\bf v}$. 

The aim of the paper is to investigate the singular limit of \eqref{meso system}--\eqref{def phi psi}. More specifically, we are concerned with the convergence analysis of the problem
\beq\label{meso system ue ve}
\begin{cases}
\partial_tu_b^\e-d_b\Delta u_b^\e=f_b(u_a^\e,u_b^\e,\,v^\e)-\e^{-1} Q(u_a^\e,\, u_b^\e,\,v^\e\,),&\In\,\,(0,\infty)\times\Omega,\\
\partial_tu^\e-\Delta(d_a u_a^\e+d_b u_b^\e)=f_u(u_a^\e,u_b^\e,v^\e),\qquad&\In\,\,(0,\infty)\times\Omega,\\
\partial_tv^\e-\Delta(d_v v^\e)=f_v(u_a^\e,u_b^\e,\,v^\e),&\In\,\,(0,\infty)\times\Omega,\\
\nabla u_b^\e\cdot \vec{n}=\nabla u^\e\cdot \vec{n}=\nabla v^\e\cdot \vec{n}=0\,,\qquad&\In\,(0,\infty)\times\pa\Omega\\
u_b^\e=u_b^\init,\quad u^\e(0)=u^\init:=u_a^\init+u_b^\init,\quad v^\e(0)=v^\init,&\In\,\Omega
\end{cases}
\eeq
as $\e\to0$ and the relative rate of convergence, where we denote
\beq\label{def fu}
f_u(u_a,u_b,v):=f_a(u_a,u_b,v)+f_b(u_a,u_b,v)\,.
\eeq

The class of system obtained in the limit and for which we investigate the existence and uniqueness issue of global solutions, reads as the cross-diffusion system
\beq\label{macro system}
\begin{cases}
\partial_tu-\Delta(A(u,v))=f_u(u_a^*(u,v),u_b^*(u,v),v),\qquad&\In\,\,(0,\infty)\times\Omega,\\
\partial_tv-\Delta(d_v\, v)=f_v(u_a^*(u,v),u_b^*(u,v),v),&\In\,\,(0,\infty)\times\Omega,\\
\nabla A(u,v)\cdot \vec{n}=\nabla v\cdot \vec{n}=0\,,\qquad&\In\,(0,\infty)\times\pa\Omega,\\
u(0)=u^\init,\quad v(0)=v^\init,&\In\,\Omega
\end{cases}
\eeq
where 
\beq\label{structure}
A(u,v):=d_a u_a^*(u,v)+d_b u_b^*(u,v)
\eeq
and $(u_a^*,u_b^*)$ is the $C^1$ maps from $\R^2_+$ to $\R^2_+$ such that, for all $(\tilde u,\tilde v)\in\R^2_+$, the pair $(u_a^*(\tilde u,\tilde v),u_b^*(\tilde u,\tilde v))$ is the unique nonnegative  solution of the nonlinear system
\beq\label{nonlinearsys}
\left\{
\begin{split}
&u_a^*+u_b^*=\tilde u\\
&Q(u_a^*,u_b^*,\tilde v)=0\,,
\end{split}
\right. 
\eeq
with $Q$ defined in \eqref{def Q}. 

Note that it is not necessary to consider initial data $u_a^\init,u_b^\init,v^\init$ belonging to the critical manifold given by \eqref{nonlinearsys}, i.e. initial data satisfying $Q(u_a^\init,u_b^\init,v^\init)=0$. Moreover, 
\[
u_a^*(u,v)=\Lambda_b^*(u,v)u\,,\qquad u_b^*(u,v)=\Lambda_a^*(u,v)u\,,
\] 
where $\Lambda_a^*,\Lambda_b^*$ are the satisfaction measures \eqref{def rel sat measures} evaluated at $(u_a^*(u,v),u_b^*(u,v),v)$. Therefore, the reaction terms in \eqref{macro system} still write as Lotka-Volterra competitive reactions, with coefficients depending on $(u,v)$. 

The emergence of cross-diffusion as a fast reaction singular limit has been observed in several mathematical models for ecology, biology or chemistry, and more specifically  in the context of competitive interactions \cite{Desvillettes2015a}, predator-prey interactions \cite{Conforto2018,Iida2023,Soresina2023}, dietary diversity and starvation \cite{Brocchieri2021} and enzyme reaction \cite{BaoBao2024}. On the other hand, cross-diffusion systems have attracted significant interest, at least since the seminal paper of  Shigesada-Kawasaki-Teramoto \cite{Shigesada1979}, because cross-diffusion can induce instability and thus explain pattern formation whereas linear diffusion cannot (see also \cite{BreKueSore2021,GAMBINO2012,Iida2006,Moussa2019} and the references therein). The mathematical analysis of cross-diffusion systems is delicate and has received a lot of attention in the last two decades. A fundamental contribution to the existence of solutions issue was given by Chen and J$\ddot{\text{u}}$ngel \cite{ChenJuengel2004, ChenJuengel2006}, for cross-diffusion systems with entropy structure, including the SKT model in \cite{Shigesada1979}. These entropy methods were shown to be robust enough to treat generalisations of the SKT system \cite{Desvillettes2014, Desvillettes2015}. Afterwards, the relation between the structure of systems involving cross-diffusion and the existence of an entropy functional has been deeply investigated (see e.g.~\cite{Daus2019,Lep2017}).  

In this article, the authors aim to give a contribution to the aforementioned tusk, showing that there are biologically relevant cross-diffusion systems, that, even though they do not have an entropy structure, they have a family of underlying energies helpful to obtain a positive answer to the existence issue. Moreover, we also address the issue of the convergence rate, as $\e$ goes to 0, of $(u^\e,v^\e)$ towards $(u,v)$ and of $(u_a^\e,u_b^\e)$ towards $(u_a^*(u,v),u_b^*(u,v))$. This problem is strictly linked with the crucial topic of invariant slow manifolds of slow-fast dynamical systems in finite dimension. The underlying theory has been extended to the case of slow-fast PDE systems, e.g. fast reaction-diffusion system, (see \cite{Desvillettes2025} and the references therein). Here, employing purely analytical tools, different from those used in \cite{Desvillettes2025}, we obtain the convergence rate and quantify its dependence on the initial layer due to initial data not lying on the critical manifold given by \eqref{nonlinearsys}.

Triangular cross-diffusion systems, similar to \eqref{macro system}, has been considered also in \cite{Brocchieri2024}, where the authors investigate the existence and weak-strong stability issues, using regularisation techniques and fixed-point arguments. However, system \eqref{macro system}--\eqref{nonlinearsys} don't fit in the class analysed in \cite{Brocchieri2024}. 

The article is organised as follows. In {\it Section} \ref{Main results} we state the main results and a set of notations. In {\it Section} \ref{energy} we introduce the family of energy functionals for \eqref{meso system}--\eqref{def phi psi}. The time evolution of the energies along the trajectories of the unique strict solution $(u_a^\e,u_b^\e,v^\e)$ of \eqref{meso system}--\eqref{def phi psi} is analysed in {\it Sections} \ref{Iprea}--\ref{sect proof energy lemma}. This analysis will provide us a priori estimates on $(u_a^\e,u_b^\e,v^\e)$, uniform in $\e$, stated in Lemma~\ref{energy estimate}. The proof of the existence result for the cross-diffusion system \eqref{macro system}--\eqref{nonlinearsys}, stated in Theorem~\ref{thm existence}, is given in {\it Section} \ref{sect proof existence macro}, while {\it Section} \ref{sect proof uniqueness macro} is devoted to the proof of the uniqueness of bounded solutions.  We conclude the article with {\it Section} \ref{sect proof rate of conv} where we obtain a convergence rate for the singular limit $\e\to0$, stated in Theorem \ref{thm rate conv}. For the sake of completeness, in Appendix \ref{appendix A} we give the existence result of the unique strict solution $(u_a^\e,u_b^\e,v^\e)$ of the fast reaction-diffusion system, stated in Theorem~\ref{th exist meso}. Finally, Appendix \ref{appendix B} is devoted to the proof of the solvability of the nonlinear system \eqref{nonlinearsys}. 
\medskip

\begin{acknowledgment}
The authors warmly thank Laurent Desvillettes and Helge Die\-tert for the fruitful discussions they had all along the preparation of the article. 
\end{acknowledgment}
\section{Statements and main results}\label{Main results}
The starting results are the existence of a unique nonnegative strict solution of system \eqref{meso system}--\eqref{def phi psi} together with basic estimates (independent of $\e$) on the solution. All of this is stated in Theorem~\ref{th exist meso} and obtained applying the classical theory of analytic semigroups, the maximal (parabolic) regularity and taking advantage of  the competition dynamics \eqref{def react functions}. We refer mainly to \cite{Lunardi1995} and we sketch the proof in Appendix \ref{appendix A}, for the reader convenience. 

Let us define
\beq\label{def Dp}
D_p:=\{w\in W^{2,p}(\Omega) : \nabla w\cdot \vec{n}=0\ \text{on}\ \partial\Omega\}\,,\quad p\in(1,+\infty)\,,
\eeq
\beq\label{notations 1}
\overline\eta:=\eta_a \vee\eta_b\,,\quad 
\eta:=a\eta_a\wedge b\eta_b\,,\quad\eta_v:=\eta_v'+\eta_v''\,,\quad 
r_v:=c\,\eta_v'+d\,\eta_v''\,.
\eeq
\begin{thm}[Well-posedness of the fast reaction-diffusion system]\label{th exist meso} Let $\Omega\!\subset\! \R^N, N \ge 1$, be a bounded open set with $C^2$ boundary $\partial\Omega$ and let $\e>0$. Assume \eqref{def react functions}--\eqref{def phi psi}, \eqref{(H1)} and let 
\beq\label{hyp ID}
u_a^\init,u_b^\init,v^\init\in \bigcap_{1<p<\infty} D_p
\eeq 
be non-negative initial data, with $u_a^\init,u_b^\init$ not identically zero. It follows that there exists a triplet
\beq\label{strict sol}
(u_a^\e,u_b^\e,v^\e)\in C^1([0,\infty);(L^p(\Omega))^3)\cap C^0([0,\infty);D_p^3)\,,\quad\forall\ p\in(1,+\infty)
\eeq  
with $u_a^\e,u_b^\e\!>\!0$ and $v^\e\!\ge\!0$ on $(0,\infty)\!\times\Omega$, which is the unique strict solution of~\eqref{meso system} with boundary conditions \eqref{HNBC} and initial conditions \eqref{ID ua ub v}. In addition, the solution satisfies the following estimates (independent of $\e>0$)
\beq\label{LinftyL1}
\| u_a^\e+u_b^\e\|_{L^\infty(0,\infty;L^1(\Omega))}\le \max\{ \| u^{\init}_a+u^{\init}_b\|_{L^1(\Omega)},2|\Omega|\overline\eta{\eta}^{-1}\}=:K_1\,,
\eeq
\beq\label{bound v}
\| v^\e\|_{L^{\infty}((0,\infty)\times\Omega)}\le \max\Big\{\|v^\init\|_{L^\infty(\Omega)},\f{\eta_v}{r_v}\Big\}=: K_\infty\,,
\eeq 
and, for all $T>0$,
\beq\label{result p=1}
\| u_a^\e\|_{L^2((0,T)\times\Omega)}^2+\|u_b^\e\|_{L^2((0,T)\times\Omega))}^2
\le {\eta}^{-1}\| u^{\init}_a+u^{\init}_b\|_{L^1(\Omega)}+\overline\eta{\eta}^{-1}\, K_1\,T=: K_2\,,
\eeq
\beq\label{max reg p=2}
\|\pa_tv^\e\|_{L^2((0,T)\times\Omega)}+\sum_{i,j}\|\pa_{x_i,x_j}v^\e\|_{L^2((0,T)\times\Omega)}\le C_1( K_2, K_\infty,T,|\Omega|)\,,
\eeq
\beq\label{est grad v L4}
\|\nabla v^\e\|_{L^4((0,T)\times\Omega)}\le C_2(K_2, K_\infty,N,T,|\Omega|)\,.
\eeq
\end{thm}

As mentioned in the introduction, the main results of the article are the existence and uniqueness of a global (weak, strong) solution $(u,v)$ of the cross-diffusion system \eqref{macro system}--\eqref{nonlinearsys}, obtained when $\e\to0$ in \eqref{meso system ue ve}, and the rate of convergence. 

The existence result for the cross-diffusion system,  stated in Theorems \ref{thm existence} and proved in Section \ref{sect proof existence macro}, needs the proof of further (uniform in $\e$) estimates. For this purpose, we construct in Sections \!\ref{energy} a well fitted family of energy functionals  \eqref{def  energy}--\eqref{def hap hbp}. The analysis of the time evolution of the energies is developed in Sections~\ref{Iprea}--\ref{sect proof energy lemma} and leads to Lemma~\ref{energy estimate} below. For the proof of this Lemma, we need to define the exponents increasing in $\alpha,\beta$
\beq\label{def q r}
\begin{split}
&q(p):=p+\alpha(p-1)=(\alpha+1)(p-1)+1\,,\\
&r(p):=p+\beta(p-1)=(\beta+1)(p-1)+1\,,
\end{split}
\eeq
and the critical values 
\beq\label{def palpha pbeta}
p_\alpha\coloneqq 1+\f1 {1+\alpha}\in (1,2),\qquad p_\beta \coloneqq 1+\f1 {1+\beta}\in (1,2)\,.
\eeq
Note that $p_\beta\le p_\alpha$ and $q(p)\le r(p)$, since $0<\alpha\le\beta$, and that
\beq\label{q=r=2}
q(p_\alpha)=r(p_\beta)=2\,.
\eeq
Furthermore, writing $r(p)=q(p)+(\beta-\alpha)(p-1)$, one see that the gap between $r(p)$ and $q(p)$ is controlled by the gap $\beta-\alpha$ between the transition rates. Therefore, the size of $\beta-\alpha$ will be crucial in the bootstrap procedure performed in the proof of Lemma~\ref{energy estimate}. To carry out the bootstrap, we define the decreasing family of intervals
\beq\label{def In}
I_n:=\big(2(\alpha+1),2(\alpha+1)+\f4{(\alpha+1)^n})\,,\quad n\in\N\cup\{0\}\,,
\eeq
so that the admissible set $[0,2(\alpha+3))$ for $\beta-\alpha$ in $\eqref{(H2)}$ reads as $[0,2(\alpha+1)]\cup(\cup_nI_n)$. Moreover, if $\beta-\alpha>2(\alpha+1)$, we also denote 
\beq\label{def nalphabeta}
n_{\alpha,\beta} \text{ the largest integer such that $\beta\!-\!\alpha\in I_{n_{\alpha,\beta}}$}\,.
\eeq
\begin{lem}[Energy estimates]\label{energy estimate}
Under hypothesis of Theorem \ref{th exist meso} and assuming in addition \eqref{(H2)}, \eqref{(H3)}, for all $T>0$, it holds
\begin{itemize}
\item[(i)] there exists $C(T)>0$ such that, for all $p\in[p_\beta,p_\alpha]$ and $\e>0$,
\beq\label{basic energy estimate lemma 2.2}
\begin{split}
&\|u_a^\e\|_{L^\infty(0,T;L^{ q(p)}(\Omega))}+\|u_b^\e\|_{L^\infty(0,T;L^{ r(p)}(\Omega))}\\ 
&\qquad+\|\nabla(u_a^\e)^{ q(p)/2}\|_{L^2((0,T)\times\Omega)}+ \|\nabla(u_b^\e)^{ r(p)/2}\|_{L^2((0,T)\times\Omega)}\\
&\qquad\qquad+\|u_a^\e\|_{L^{ q(p)+1}((0,T)\times\Omega)}+\|u_b^\e\|_{L^{ r(p)+1}((0,T)\times\Omega)}\le C(T)\,;
\end{split}
\eeq
\item[(ii)] if $\beta-\alpha\in[0, 2(\alpha+1)]$, for all $p\in[2,+\infty)$, there exists $C(p,T)>0$ such that, for all $\e>0$,
\beq\label{energy estimate 1 lemma 2.2}
\begin{split}
&\|u_a^\e\|_{L^\infty(0,T;L^{q(p)}(\Omega))}+\|u_b^\e\|_{L^\infty(0,T;L^{r(p)}(\Omega))}\\
&\quad+\|\nabla(u_a^\e)^{q(p)/2}\|_{L^2((0,T)\times\Omega)}+ \|\nabla(u_b^\e)^{r(p)/2}\|_{L^2((0,T)\times\Omega)}\le C(p,T)\,;
\end{split}
\eeq
\item[(iii)]  if $\beta\!-\!\alpha\in(2(\alpha+1),2(\alpha+3))=\cup_nI_n$, there exists $C(n_{\alpha,\beta},T)\!>\!0$ such that, for all $p\!\in\![2,1+(\alpha\!+\!1)^{n_{\alpha,\beta}}]$ and~$\e>0$, 
\begin{align}\label{energy estimate 2 lemma 2.2}
&\|u_a^\e\|^{ q(p)}_{L^\infty(0,T;L^{ q(p)}(\Omega))}+\|u_b^\e\|^{ r(p)}_{L^\infty(0,T;L^{ r(p)}(\Omega))}\notag\\ 
&\qquad+\|\nabla(u_a^\e)^{ q(p)/2}\|_{L^2((0,T)\times\Omega)}^2+ \|\nabla(u_b^\e)^{ r(p)/2}\|_{L^2((0,T)\times\Omega)}^2\\
&\qquad\qquad+\|u_a^\e\|_{L^{ q(p)+1}((0,T)\times\Omega)}^{ q(p)+1}+\|u_b^\e\|_{L^{ r(p)+1}((0,T)\times\Omega)}^{ r(p)+1}\le C(n_{\alpha,\beta},T)\notag\,.
\end{align}
\item[(iv)] there exists $C(T)>0$ such that, for all $\e>0$,
\beq\label{evol Ep p=2 final}
\e^{-\f12}\|\Lambda^{1/2}Q\|_{L^2((0,T)\times\Omega)}\le C(T)(1+\|u_a^\e\|_{L^2((0,T)\times\Omega)} +\|u_b^\e\|_{L^2((0,T)\times\Omega)})\,. 
\eeq
\end{itemize}

Finally, for all $T>0$ and $p\in(1,+\infty)$ if $\beta-\alpha\in[0,2(\alpha+1)]$, or $p\in(1,2+(\alpha\!+\!1)^{n_{\alpha,\beta}+1}])$ if $\beta-\alpha\in(2(\alpha+1),2(\alpha+3))$, there exist $C_1(p,T)\!>\!0$ and $C_2(p,T)\!>\!0$ such that, for all $\e>0$, it holds
\begin{align}
&\|\partial_{t}v^\e\|_{L^p((0,T)\times\Omega)}+\sum_{i,j}\|\partial_{x_ix_j} v^\e\|_{L^{{p}}((0,T)\times\Omega)}\le C_1\label{lap v}\\
&\|\nabla v^\e\|_{L^{2p}((0,T)\times\Omega)}\leq C_2\label{gradient v}\,.
\end{align}
\end{lem}
\begin{remark}\label{remark grad u in L2}
Thanks to identity \eqref{q=r=2}, we obtain the $L^2((0,T)\times\Omega)$ uniform estimate of $\nabla u_a^\e$ and $\nabla u_b^\e$ from \eqref{basic energy estimate lemma 2.2}, taking $p=p_\alpha$ and $p=p_\beta$, respectively.  

\end{remark}

Next, the cross-diffusion term in \eqref{macro system} is due to the convergence of the pair $(u_a^\e,u_b^\e)$ towards the unique solution of the nonlinear system~\eqref{nonlinearsys}. Estimate \eqref{conv slow manifold} in Lemma \ref{converg (ua,ub)} (proved in Section \ref{sect proof existence macro}) is the key tool to obtain this convergence. The solvability of \eqref{nonlinearsys} and the regularity of the map $(u_a^*,u_b^*)\!:\R^2_+\mapsto~\R^2_+$, in turn, are straightforward consequences of the regularity of $Q$ and of the implicit function theorem. These properties are resumed in Lemma~\ref{property cross-diff} below and the proof is given for completeness in Appendix~\ref{appendix B}. 

\begin{lem}[Existence and regularity of the map $(u_a^*,u_b^*)$]\label{property cross-diff}
Assume \eqref{def Q},\eqref{def phi psi}, \eqref{(H1)}. For all $(\tilde u,\tilde v)\!\in\!\R_+^2$, there exists a unique nonnegative solution $(u_a^*(\tilde u,\tilde v),u_b^*(\tilde u,\tilde v))$ of \eqref{nonlinearsys}. Moreover, $(u_a^*,u_b^*)\in(C^1(\R_+^2))^2$, with 
\beq\label{partial u ua* ub*}
\pa_{\tilde u}u_a^*(\tilde u,\tilde v),\,\pa_{\tilde u}u_b^*(\tilde u,\tilde v)\in(0, 1)
\eeq
and, assuming $\afast$ and $\bfast$ both non-zero,
\beq\label{partial v ua* ub*}
-\f{\cfast}{\afast}\le \pa_{\tilde v}u_a^*(\tilde u,\tilde v)=-\pa_{\tilde v}u_b^*(\tilde u,\tilde v)\le\f{\dfast}{\bfast}\,.
\eeq
\end{lem}

\begin{lem}[Convergence of $(u_a^\e,u_b^\e)$]\label{converg (ua,ub)} 
Under assumptions of Lemma \ref{energy estimate}, let $(u_a^\e,u_b^\e,v^\e)$ be the unique strict solution of \eqref{meso system}--\eqref{def phi psi}, given by Theorem \ref{th exist meso}, and let  $u^\e=u_a^\e+u_b^\e$. For all $T\!>\!0$, there exists $C(T)\!>\!0$ such that, for all $\e\!>\!0$, it holds
\beq\label{conv slow manifold}
\|u_b^\e-u_b^*(u^\e,v^\e)\|_{L^2((0,T)\times\Omega)}\le \|Q(u_a^\e,u_b^\e,v^\e)\|_{L^2((0,T)\times\Omega)}\le A^{-\f\alpha2}\,C(T)\sqrt\e\,. 
\eeq
As $|u_a^\e-u_a^*(u^\e,v^\e)|=|u_b^\e-u_b^*(u^\e,v^\e)|$, the same inequality holds for $u_a^\e-u_a^*(u^\e,v^\e)$. 
\end{lem}

Having these results at hand, we can state the existence of a global (weak, strong) solution $(u,v)$ of the cross-diffusion system. When $\beta-\alpha\!\in\!(2(\alpha+1),\,2(\alpha+3))$, we will: (i) use the best integrability for $u_a^\e$ given in \eqref{energy estimate 2 lemma 2.2}, namely the $L^{q(p)+1}$ integrability with $p\!=\!1+(\alpha+1)^{n_{\alpha,\beta}}$, that gives $q(p)+1\!=\!2+(\alpha+1)^{n_{\alpha,\beta}+1}$; (ii) add a dimension depending condition on $\beta-\alpha$. 

%

\begin{thm}[Existence for the cross-diffusion system]\label{thm existence}
Let $\Omega\subset \R^N, N \ge 1$, be a bounded open set with $C^2$ boundary $\partial\Omega$ and assume \eqref{(H1)}, \eqref{(H2)}, \eqref{(H3)}. Let $u_a^\init,u_b^\init, v^\init$ be non-negative initial data satisfying \eqref{hyp ID}, with $u_a^\init,u_b^\init$ not identically zero. Furthermore, if $N>6$ and $\beta-\alpha\!\in\!(2(\alpha+1),\,2(\alpha+3))$, assume that $(\alpha+1)^{n_{\alpha,\beta}+1}\ge2$. There exists a pair of nonnegative measurable functions $(u,v):(0,\infty)\times\Omega\to\R^2_+$ such that, for all $p$ satisfying
\begin{equation}\label{def p thm existence macro}
     \begin{cases}
         p\in(2,+\infty), \,\,\,\,\,\qquad\qquad\qquad\text{ if } \beta-\alpha\in [0,\,2(\alpha+1)],\\
         p=2+(\alpha+1)^{n_{\alpha,\beta}+1},\qquad\quad\, \text{if } \beta-\alpha\in(2(\alpha+1),\,2(\alpha+3)),
     \end{cases}
 \end{equation}
(see \eqref{def In}, \eqref{def nalphabeta}), $s=\f p2\wedge 2$ and all $T>0$, it holds
\begin{itemize}
\item [(i)] $u\in L^2(0,T;H^1(\Omega))\cap L^p((0,T)\times\Omega)\cap L^\infty(0,T;L^{p-1}(\Omega))$ and\\ $\pa_tu\in L^s(0,T;H^{-1}(\Omega))$,
\item [(ii)] $v\in W^{1,p}((0,T)\times\Omega)\cap L^{p}\big(0,T; W^{2,p}(\Omega)\big)\cap L^{\infty}((0,\infty)\times\Omega)$,
\item [(iii)] up to the extraction of a subsequence from the strict solutions $(u_a^{\e},u_b^{\e},v^\e)_\e$
\[
(u_a^{\e},u_b^{\e})\to(u_a^*(u,v),u_b^*(u,v))\,,\ \text{ a.e. in } (0,\infty)\times\Omega\,,\text{ as }\e\to0\,,
\]
where $(u^*_a(u,v),u^*_b(u,v))$ is the unique solution of \eqref{nonlinearsys}\,,
\item [(iv)] $u^*_a(u,v),u^*_b(u,v)\in L^2(0,T;H^1(\Omega))\cap L^p((0,T)\times\Omega)\cap L^\infty(0,T;L^{p-1}(\Omega))$\,,
\item [(v)] $(u,v)$ is a global (weak, strong) solution of \eqref{macro system}, \eqref{structure}, \eqref{nonlinearsys}, i.e. for all $T>0$ and all $w\in C^1([0,T];H^1(\Omega))$ such that $w(T)=0$, it holds 
\beq\label{weak eq u}
\begin{cases}
&-\int_0^T\!\!\int_\Omega u\,\partial_{t}w\,dxdt
+\int_0^T\!\!\int_\Omega\nabla(d_au^*_a(u,v)+d_bu^*_b(u,v))\cdot\nabla w\,dxdt\\[1.5ex]
&\qquad\quad=\int_{\Omega}u^{\init}\,w(0)\,dx+ \int_0^T\!\!\int_\Omega f_u(u_a^*(u,v), u_b^*(u,v),v)\,w\,dxdt\\[1.5ex]
&\pa_tv=d_v\Delta v+f_v(u_a^*(u,v),u_b^*(u,v),v)\quad\text{ in } L^p((0,T)\times\Omega)\\[1.5ex]
& \sum_{i}\gamma_0(\pa_{x_i}v)\,n_i=0\,,\quad \text{a.e. in } (0,T)\times\pa\Omega\\[1.5ex]
&u(0)= u^\init\,,\quad v(0)=v^\init\quad \text{a.e. in } \Omega
\end{cases}
\eeq
\end{itemize}
\end{thm}

One can shows that the global (weak, strong) solution $(u,v)$ obtained in Theorem~\ref {thm existence} enjoys additional regularity properties.  However, exploring the regularity of $(u,v)$ lies beyond the scope of the paper. Instead, we conclude the analysis by proving the uniqueness of solutions of \eqref{macro system}--\eqref{nonlinearsys} that are bounded in both components. Note that the boundedness of the component $u$ is required only to handle the reaction. Theorem~\ref{Thm uniqueness} is proved in Section~\ref{sect proof uniqueness macro}. 
\begin{thm}[Uniqueness and stability]\label{Thm uniqueness}
Under the assumptions of Theorem \ref{thm existence}, let $(u_i,v_i)$, $i=1,2$, be two solutions with initial data $(u_i^\init,v_i^\init)$. Assume in addition that, for all $T>0$, $u_i\in L^\infty((0,T)\times\Omega)$, $i=1,2$, and $\afast, \bfast$ are both non-zero. It follows that, for all $T>0$, there exists $C(T, \|u_i\|_{L^\infty((0,T)\times\Omega)},\|v_i\|_{L^\infty((0,\infty)\times\Omega)})>0$ such that  
\beq\label{est stability}
\begin{split}
\|u_1-u_2\|^2_{L^2((0,T)\times\Omega)}&+\|v_1-v_2\|^2_{L^2((0,T)\times\Omega)}\\
&\le\! C\!\left(\|u_1^\init-u_2^\init\|^2_{L^2(\Omega)}+\|v_1^\init-v_2^\init\|^2_{L^2(\Omega)}\right).
\end{split}
\eeq
\end{thm}

Concerning the rate of convergence issue, we analyse the time evolution of the $L^2(\Omega)$ norms of $u^\e\!-u$, $v^\e\!-\!v$, $u_b^\e-u_b^*(u,v)$ via the functional \eqref{def L} which include the ad hoc sub-functional \eqref{def E}--\eqref{def caP} designed to handle the fast reaction. The result is established under additional regularity assumptions of the solutions and small cross-diffusion, i.e. small $|d_b-d_a|$. A careful reading of the proof makes it clear that the diffusivity coefficient $d_v$ can only help to optimise the smallness condition \eqref{hyp rate of convergence}, but cannot remove it. The obtained estimates also illustrate how the initial layer
\beq\label{initial layer}
\e_\init:=\|u_b^\init-u_b^*(u^\init,v^\init)\|_{H^1(\Omega)}
\eeq
slows down the convergence rate. 
\begin{thm}[Rate of convergence]\label{thm rate conv}
Under the hypothesis of Theorem \ref{thm existence}, let $(u_a^\e,u_b^\e,v^\e)$ be the unique nonnegative global strict solution of \eqref{meso system}--\eqref{def phi psi}. 
Assume in addition $d_b>d_a$, $A,B>0$ for all $\alpha,\beta$ and that $u^\e=u_a^\e+u_b^\e$ is uniformly bounded locally in time, i.e.  for all $T>0$ there exists $M_T>0$ such that, for all $\e\in(0,1)$, it holds
\beq\label{bound assumpt1}
\|u^\e\|_{L^\infty((0,T)\times\Omega)}\le M_T\,.
\eeq
Let $(u,v)$ be a nonnegative global classical solution of \eqref{macro system}--\eqref{nonlinearsys} such that, for all $T>0$, 
\beq\label{bound assumpt2}
u,v\in C^0([0,T];C^3(\overline\Omega))\cap C^1([0,T];C^1(\overline\Omega))\,.
\eeq
 Then, for all $T>0$, there exists a constant $C_1(\alpha,\beta,A,B,T)>0$ such that, if
\beq\label{hyp rate of convergence}
1<\f{d_b}{d_a}<1+C_1(\alpha,\beta,A,B,T)\,,
\eeq
there exists $C_2(T)>0$ such that, for all $\e\in(0,1)$ and for $\e_\init$ defined in \eqref{initial layer}, it holds 
\begin{align}
&\|u^\e-u\|_X+\| v^\e-v\|_X\le C_2(T)(\e+\e^{\f12}\,\e_\init)\,,\label{u-ue v-ve}\\
&\|u_a^\e-u_a^*(u,v)\|_Y+\|u_b^\e-u_b^*(u,v)\|_Y\le C_2(T)\,(\e+\e^{\f12}\,\e_\init)\,,\label{uae-ua* ube-ub*}
\end{align}
where $\|\cdot\|_X:=\|\cdot\|_{L^2(0,T;H^1(\Omega))}+\|\cdot\|_{ L^\infty(0,T;L^2(\Omega))}$ and $\|\cdot\|_Y:=\|\cdot\|_{L^2(0,T;H^1(\Omega))}$. 
\end{thm}

A by-product of the above convergence result is again the uniqueness of smooth solutions, without the requirements $\afast\neq0$ and $\bfast\neq0$.
\medskip

\noindent{\bf Notations.} Hereafter, $\nabla$ and $\Delta$ will denote the gradient with respect to the spatial variable $x$ and the Laplacian, while $D$ and $D^2$ will denote the gradient and the Hessian with respect to non spatial variables. Moreover, $\pa_i$ and $\pa_{ij}$ will denote respectively  the partial derivative with respect to the $i$-$th$ variable and the partial derivative with respect to the $i$-$th$ and $j$-$th$ variables, whatever they are. For the sake of simplicity, in the computations of Sections \ref{energy}--\ref{sect proof energy lemma} we will keep explicit only the constants depending on $p\ge2$, and we will write $x\lesssim y$ meaning that there exists a universal constant $C>0$, not depending on $p$, such that $x\le C\, y$. Finally, when there is no risk of confusion, we will denote $\Omega_T$, $T>0$, the usual open cylinder $(0,T)\times\Omega$ and we will omit the $\e$ superscript. 
\section{A family of energy functionals $\Ep(u_a,u_b,v)$}\label{energy}
System \!\eqref{meso system}--\eqref{def phi psi} is naturally endowed with the following family\! of \! energy \!functionals
\beq\label{def  energy}
\Ep(u_a,u_b,v):=\int_{\Omega}h_{p}(u_a,u_b,v)\,dx,\qquad u_a,u_b,v\in\R_+\,,\qquad p\ge1\,,
\eeq
with the total energy density $h_p$ given by
\beq\label{def hp}
h_p(u_a,u_b,v):=h_{a,p}(u_a,v)+h_{b,p}(u_b,v)
\eeq
and the partial energy densities $h_{a,p},h_{b,p}$ defined as
\beq\label{def hap hbp}
\begin{split}
h_{a,p}(u_a,v)&\coloneqq \int_{0}^{u_a}\psi^{p-1}(\afast z+\cfast v)z^{p-1}dz\,,\\
h_{b,p}(u_b,v)&\coloneqq\int_{0}^{u_b}\phi^{p-1}(\bfast z+\dfast v)z^{p-1}dz\,.
\end{split}
\eeq
Moreover, setting 
\beq\label{def theta omega}
\theta(z,v)\coloneqq A+ \afast z+\cfast v,\qquad
\omega(z,v)\coloneqq B+ \bfast z+\dfast v,
\eeq
by \eqref{def phi psi}, the partial energy densities in \eqref{def hap hbp} rewrite  as
\beq\label{def hap hbp bis}
\begin{split}
h_{a,p}(u_a,v)&= \int_{0}^{u_a}\theta^{\alpha(p-1)}(z,v)z^{p-1}dz\,,\\
h_{b,p}(u_b,v)&=\int_{0}^{u_b}\omega^{\beta(p-1)}(z,v)z^{p-1}dz\,.
\end{split}
\eeq

\begin{remark}
It is worth noticing that in definition \eqref{def hap hbp} the transition functions $\psi,\phi$ are not renormalised by $\Lambda(u_a,u_b,v)$ unlike in $Q(u_a,u_b,v)$ (see \eqref{def Q}). This choice will be clear later.  
\end{remark}

The interest in the family of energies \eqref{def  energy}--\eqref{def hap hbp bis} is threefold. It allows us to obtain further a priori estimates on the densities $u_a^\e$, $u_b^\e$ and their gradients in Lebesgue spaces, to handle easily the contribution due to the fast reaction $\e^{-1}Q$ in the aforementioned estimates, and to obtain the convergence of  $\Lambda^{\f12}Q$  towards $0$, as $\e\to0$. Indeed, on the one hand, using for all $z,v\ge0$, 
\beq\label{ineq theta omega}
\theta(z,v)\ge\afast\,z\qquad\text{ and }\qquad \omega(z,v)\ge\bfast\,z\,,
\eeq
it is easily seen from \eqref{def hp}, \eqref{def hap hbp bis}, that, for $p\ge1$,
\[
\afast^{\alpha(p-1)}\f{u_a^{q(p)}}{q(p)}+ \bfast^{\beta(p-1)}\f{u_b^{r(p)}}{r(p)} \le h_p(u_a,u_b,v)\le \afast^{-p}\f{\theta(u_a,v)^{q(p)}}{q(p)}+ \bfast^{-p}\f{\omega(u_b,v)^{r(p)}}{r(p)}
\]
implying that $\mathcal{E}_p$ is well defined, for all $p\ge1$, along the trajectories of the solution $(u_a^\e,u_b^\e,v^\e)$ obtained in Theorem \ref{th exist meso}, and also that, for all $T>0$, 
\beq\label{est from E}
\|u_a^\e\|^{q(p)}_{L^\infty(0,T;L^{q(p)}(\Omega))}+\|u_b^\e\|^{r(p)}_{L^\infty(0,T;L^{r(p)}(\Omega))}\le \f{r(p)}{\afast^{\alpha(p-1)}\wedge\bfast^{\beta(p-1)}}\,\mathcal{E}_p(T).
\eeq

On the other hand, denoting ${\mathcal F}:=(f_a,f_b,f_v)^T$, the evolution of $\Ep$ along the solution is described by the differential equation
\beq\label{evol energy} 
\begin{split}
\Ep'(t)&=\,\f d{dt}\int_{\Omega}h_p(u_a^\e,u_b^\e,v^\e)\,dx=\int_{\Omega}(\pa_1h_p\pa_tu_a^\e+\pa_2h_p\pa_tu_b^\e+\pa_3h_p\pa_tv^\e)dx\\
&=\int_{\Omega}(d_a\pa_1h_p\Delta u_a^\e+d_b\pa_2h_p\Delta u_b^\e+d_v\pa_3h_p\Delta v^\e)dx\\
&\quad+\int_\Omega D h_p\cdot {\mathcal F}\,dx
+\f1\e\int_\Omega(\pa_1h_p-\pa_2h_p)Q\,dx\\
& \eqqcolon I^p_{\text{diff}}+I^p_{\text{rea}}+I^p_{\text{fast}}\,.
\end{split}
\eeq
Then, we see that, for all $p\ge1$, it holds
\[
\begin{split}
 I^p_{\text{fast}}&:=\f1\e\int_\Omega(\pa_1h_p-\pa_2h_p)Q\,dx\\
 &\ =-\f1\e\int_\Omega\big[\big(\phi(\bfast u_b^\e+\dfast v^\e)u_b^\e\big)^{p-1}-\big(\psi(\afast u_a^\e+\cfast v^\e)u_a^\e\big)^{p-1}\big]\,Q\,dx\,.
 \end{split}
 \]
As $x\mapsto x^{p-1}$ is an increasing function on $\R_+$, the latter and \eqref{def Q} gives
\beq\label{Ipfats neg}
I^p_{\text{fast}}\le0\,,\qquad \forall\ p\ge1\,,
\eeq
so that $I^p_{\text{fast}}$ can be neglected in the evolution equation \eqref{evol energy}, whenever it is useless, i.e. $p\neq2$.  When $p=2$, $I^2_{\text{fast}}$ reads as
\beq\label{Ifreap2}
I^2_{\text{fast}}=-\f1\e\int_\Omega q(u_a^\e,u_b^\e,v^\e)\,Q(u_a^\e,u_b^\e,v^\e)\,dx=-\f1\e\int_\Omega \Lambda\,Q^2(u_a^\e,u_b^\e,v^\e)\,dx\,,
\eeq
thus allowing to obtain the convergence of $\|\Lambda^{\f12}Q\|_{L^2(\Omega_T)}$ towards zero, as $\e\to0$, with rate $\e^{\f12}$, under assumption \eqref{(H2)}, (see \eqref{evol Ep p=2 final}). This convergence result will be crucial in proving the convergence of the solution of the fast reaction system towards the solution of cross-diffusion system. 

Finally, when $p=1$ the total energy density \eqref{def hp} reduce to $h_1(u_a,u_b,v)=u_a+u_b$, so that $I^1_{\text{diff}}=I^1_{\text{fast}}=0$. It follows the uniform control of the densities $u_a^\e$, $u_b^\e$ in the Lebesgue spaces $L^{\infty}(0,\infty; L^1(\Omega))$ and $L^2(\Omega_T)$ obtained in Theorem \ref{th exist meso}. In order to estimates $I^p_{\text{rea}}$ and $I^p_{\text{diff}}$ in \eqref{evol energy} with $p>1$ and bootstrap the above $L^1$ and $L^2$ estimates to $L^p$ estimates, $p>2$, (see Lemma \ref{energy estimate}), we need instead a suitable analysis of the hessian matrix $\text{Hess}(h_p)$. Particular attention must be paid to the critical case $p=p_\beta$ (see\eqref{def palpha pbeta}), which requires to assume $A>0$ to control the $\text{Hess}(h_{a,p})$, (see \eqref{sub critical case}). This analysis is done in the rest of this section. 

\begin{remark}
The energy functionals \eqref{def  energy} are reminiscent of the functionals introduced in \cite{Desvillettes2015a,Brocchieri2021}. It is worth noticing that they are not the sum of functionals of the single densities $u_a$, $u_b$ and $v$. 
\end{remark}

\subsection{The gradient of the total energy density $h_p$} 
Let $p>1$. From  \eqref{def hp},\eqref{def theta omega},\eqref{def hap hbp bis},  the gradient $D h_p$ reads as 
\[
D h_p(u_a,u_b,v)=\Big(\pa_1 h_{a,p}(u_a,v),\pa_1 h_{b,p}(u_b,v),\pa_2 h_{a,p}(u_a,v)+\pa_2 h_{b,p}(u_a,v)\Big),
\]
where
\beq\label{def pa1ha pa1hb}
\begin{split}
\pa_{1}h_{a,p}(u_a,v)=\theta(u_a,v)^{\alpha(p-1)}u_a^{p-1},\quad
\pa_{1}h_{b,p}(u_b,v)=\omega(u_b,v)^{\beta(p-1)}u_b^{p-1},
\end{split}
\eeq
and 
\begin{equation}\label{def pa2ha pa2hb}
\begin{split}
\pa_2 h_{a,p}(u_a,v)&=\cfast\alpha(p-1)\int_{0}^{u_a}\theta(z,v)^{\alpha(p-1)-1}z^{p-1}dz,\\
\pa_2 h_{b,p}(u_b,v)&=\dfast\beta(p-1)\int_{0}^{u_b}\omega(z,v)^{\beta(p-1)-1}z^{p-1}dz.
\end{split}
\end{equation}

The derivatives $\pa_{1}h_{a,p}$ and $\pa_{1}h_{b,p}$ in \eqref{def pa1ha pa1hb} are well defined for all $u_a,u_b,v\ge0$ and all $A,B\ge0$. The same holds true for $\pa_2 h_{a,p}$ and $\pa_2 h_{b,p}$ since the integrals in \eqref{def pa2ha pa2hb} are finite. Indeed, by  \eqref{ineq theta omega} and \eqref{def q r}, we have  
\[
\begin{split}
0\le\theta(z,v)^{\alpha(p-1)-1}z^{p-1}\le \afast^{1-p}\,\theta(z,v)^{(1+\alpha)(p-1)-1}=\afast^{1-p}\,\theta(z,v)^{q(p)-2}\,,\\
0\le\omega(z,v)^{\beta(p-1)-1}z^{p-1}\le \bfast^{1-p}\,\omega(z,v)^{(1+\beta)(p-1)-1}=\bfast^{1-p}\,\omega(z,v)^{r(p)-2}\,,
\end{split}
\]
where $q(p),r(p)>1$ as $p>1$. Hence, by integration, we end up with 
\beq\label{est pa2 hap hbp}
\begin{split}
0\le\pa_2 h_{a,p}(u_a,v)\lesssim\afast^{-p}\,\theta(u_a,v)^{q(p)-1}\,,\\
0\le\pa_2 h_{b,p}(u_b,v)\lesssim\bfast^{-p}\,\omega(u_b,v)^{r(p)-1}\,.
\end{split}
\eeq
\subsection{The Hessian of the total energy density $h_p$}
Let $p>1$. The Hessian matrix of $h_p$ is 
\[
\text{Hess}(h_p)=\left(
\begin{array}{ccc}
\pa_{11}h_{a,p} & 0 & \pa_{12}h_{a,p}\\
0 & 0&  0\\
\pa_{21}h_{a,p} & 0 & \pa_{22}h_{a,p}
\end{array}
\right)
+\left(
\begin{array}{ccc}
0& 0 & 0\\
0 & \pa_{11}h_{b,p}&  \pa_{12}h_{b,p} \\
0 & \pa_{21}h_{b,p}  & \pa_{22}h_{b,p} 
\end{array}
\right),
\]
where
\begin{align}\label{hess hap}
\pa_{11}h_{a,p}(u_a,v)&=\afast\alpha(p-1)\theta(u_a,v)^{\alpha(p-1)-1}u_a^{p-1}+(p-1)\theta(u_a,v)^{\alpha(p-1)}u_a^{p-2},\notag\\
\pa_{12}h_{a,p}(u_a,v)&=\pa_{21}h_{a,p}(u_a,v)=\cfast\alpha(p-1)\theta(u_a,v)^{\alpha(p-1)-1}u_a^{p-1},\\
\pa_{22} h_{a,p}(u_a,v)&=\cfast^2\alpha(p-1)(\alpha(p-1)-1)\int_{0}^{u_a}\theta(z,v)^{\alpha(p-1)-2}z^{p-1}dz,\notag
\end{align}
and 
\begin{align}\label{hess hbp}
\pa_{11}h_{b,p}(u_b,v)&=\bfast\beta(p-1)\omega(u_b,v)^{\beta(p-1)-1}u_b^{p-1}+(p-1)\omega(u_b,v)^{\beta(p-1)}u_b^{p-2},\notag\\
\pa_{12}h_{b,p}(u_b,v)&=\pa_{21}h_{b,p}(u_b,v)=\dfast\beta(p-1)\omega(u_b,v)^{\beta(p-1)-1}u_b^{p-1},\\
\pa_{22} h_{b,p}(u_b,v)&=\dfast^2\beta(p-1)(\beta(p-1)-1)\int_{0}^{u_b}\omega(z,v)^{\beta(p-1)-2}z^{p-1}dz.\notag
\end{align}

\subsubsection{Hess$(h_{a,p})$}
From \eqref{hess hap}, we see that, for $p\in(1,2)$, we need a strictly positive density $u_a$ in order to define $u_a^{p-2}$ in $\pa_{11}h_{a,p}$. More precisely, this is necessary when $p=p_\alpha<2$ and $p=p_\beta<2$. The strict positivity of $u_a^\e$ is given by Theorem~\ref{th exist meso}.

Furthermore, the following cases will be considered to control the term appearing in both $\pa_{11}h_{a,p}$ and $\pa_{12}h_{a,p}$, i.e. $\theta(u_a,v)^{\alpha(p-1)-1}u_a^{p-1}$ . 
\begin{itemize}
\item[(a1)] If $1<p<p_\alpha$, then $\f1{p-1}-\alpha>1$ and, for all $u_a,v\ge0$, as $A>0$ (see \eqref{(H1)}), it holds
\beq\label{sub critical case}
0\le\theta(u_a,v)^{\alpha(p-1)-1}u_a^{p-1}=\left(\f{u_a}{\theta(u_a,v)^{\f1{p-1}-\alpha}}\right)^{p-1}\le C_A(p)\,.
\eeq
Estimate \eqref{sub critical case} will be crucial when $p=p_\beta$, since $p_\beta< p_\alpha$, if $\alpha<\beta$.
\item[(a2)] If $p\ge p_\alpha$,  then $\alpha(p-1)-1\ge-(p-1)$ and, for all $u_a,v\ge0$, by \eqref{ineq theta omega} it holds
\beq\label{critical+supercritical case}
0\le\theta(u_a,v)^{\alpha(p-1)-1}u_a^{p-1}\le\afast^{1-p} \theta(u_a,v)^{(1+\alpha)(p-1)-1}=\afast^{1-p} \theta(u_a,v)^{q(p)-2}\,,
\eeq
with $q(p)\ge q(p_\alpha)=2$, (see \eqref{q=r=2}).
\end{itemize}

Finally, as $A>0$,  $\pa_{22}h_{a,p}$ in \eqref{hess hap} is well defined. In order to simplify the computations, we will handle $\pa_{22}h_{a,p}$ only for $p\!>\!p_\alpha$ and the cases to be analysed are the~following.
\begin{itemize}
\item[(a3)] If  $p_\alpha<p<1+\f1\alpha$, then $-(p-1)<\alpha(p-1)-1<0$ and by \eqref{ineq theta omega}, we have
\[
0\le\theta(z,v)^{\alpha(p-1)-2}z^{p-1}\le \afast^{1-p}\, \theta(z,v)^{(1+\alpha)(p-1)-2}= \afast^{1-p}\, \theta(z,v)^{q(p)-3}\,.
\]
Hence
\beq\label{est pa22 hap}
\begin{split}
&|\pa_{22}h_{a,p}(u_a,v)|\lesssim \alpha(p-1)\,C(\alpha,p)\,\afast^{-p}\,\theta(u_a,v)^{q(p)-2}\,,\\
&C(\alpha,p):=\f{1-\alpha(p-1)}{(\alpha+1)(p-1)-1}>0\,.
\end{split}
\eeq
\item[(a4)] If $p\ge 1+\f1\alpha$, then $\alpha(p-1)-1\ge0$, so that $\pa_{22}h_{a,p}$ is positive and gives a negative term (see \eqref{Ipdiff 1}) in the evolution equation~\eqref{evol energy} that will be neglected. 
\end{itemize}

\subsubsection{Hess$(h_{b,p})$}\label{Hess(hbp)}
A similar analysis can be done for $\pa_{11}h_{b,p}$ and $\pa_{12}h_{b,p}$ in \eqref{hess hbp}. As before, we see that, for $p\in(1,2)$, we need a strictly positive density $u_b$, in order to define $u_b^{p-2}$. On the other hand, the critical case corresponding to the (a1) case above will not appear, since we will not consider any $p\in(1,p_\beta)$. This is one of the reason why, we do not need to assume $\phi$ strictly positive for all $\beta>0$ (i.e. $B>0$ for all $\beta>0$) as we do for $\psi$, (see \eqref{(H1)}). However, we need to make sure that $\pa_{22}h_{b,p}$ is well defined. To this end, we proceed as follows.

First, concerning $\pa_{11}h_{b,p}$ and $\pa_{12}h_{b,p}$, we observe, as in the (a2) case, that
\begin{itemize}
\item[(b1)] if $p\ge p_\beta$, then $\beta(p-1)-1\ge -(p-1)$ and, for all $u_b,v\ge0$,  it holds
\beq\label{critical+supercritical case 2}
0\le\omega(u_b,v)^{\beta(p-1)-1}u_b^{p-1}\le\bfast^{1-p}\,\omega(u_b,v)^{(1+\beta)(p-1)-1}=\bfast^{1-p}\,\omega(u_b,v)^{r(p)-2}\,,
\eeq
with  $r(p)\ge r(p_\beta)=2$,  (see \eqref{q=r=2}).
\end{itemize}
Next, as $B\ge0$ when $\beta\ge1$, $\pa_{22}h_{b,p}$ is not well defined when $p=p_\beta$, since the integrability of the function $\omega(z,v)^{\beta(p-1)-2}z^{p-1}$ in the neighbourhood of $z=0^+$ is not guaranteed, as it holds
\[
\omega(z,v)^{\beta(p_\beta-1)-2}z^{p_\beta-1}=\omega(z,v)^{-1-\f1{\beta+1}}z^{\f1{\beta+1}}\,.
\] 
Therefore, we will avoid this criticality and the following cases are considered. 
\begin{itemize}
\item[(b2)] If $p_\beta<p<1+\f1\beta$ then $-(p-1)<\beta(p-1)-1<0$ and we have by \eqref{ineq theta omega}
\[
0\le\omega(z,v)^{\beta(p-1)-2}z^{p-1}\le\bfast^{1-p}\omega(z,v)^{(1+\beta)(p-1)-2}=\bfast^{1-p}\omega(z,v)^{r(p)-3}\,.
\]
Hence, as $r(p)>r(p_\beta)=2$, $\pa_{22}h_{b,p}$ is well defined and 
\beq\label{est pa22 hbp}
\begin{split}
&|\pa_{22}h_{b,p}(u_b,v)|\lesssim\beta(p-1)\,C(\beta,p)\,\bfast^{-p}\omega(u_b,v)^{r(p)-2}\,,\\
&C(\beta,p):=\f{1-\beta(p-1)}{(\beta+1)(p-1)-1}>0\,.
\end{split}
\eeq
\item[(b3)] If $p\ge 1+\f1\beta$, then $\beta(p-1)-1\ge0$ and $\pa_{22}h_{b,p}$ is positive and gives a negative term (see \eqref{Ipdiff 1}) in the evolution equation \eqref{evol energy} that will be neglected. 
\end{itemize}
\section{The reaction contribution $I^p_{\text{rea}}$ to~$\Ep'(t)$ in \eqref{evol energy}}\label{Iprea}
In the sequel, it will be useful to employ the following elementary interpolation inequality: for all $x\ge0$, $C>0$ and $\gamma,\gamma_1\in\R$ such that $0<\gamma<\gamma_1$, it holds
\beq\label{elem. ineq}
x^\gamma\le C^{\gamma}+x^{\gamma_1}\,C^{\gamma-\gamma_1}\,.
\eeq

Let $p>1$. From \eqref{evol energy} we have that $I^p_{\text{rea}}=\int_\Omega D h_p\cdot {\mathcal F}\,dx$, i.e. 
\begin{align}\label{Ip rea}
I^{p}_{\text{rea}}&=\int_{\Omega}\pa_{1}h_{a,p}(u_a,v)f_a(u_a,u_b,v)\,dx+\int_{\Omega}\pa_{1}h_{b,p}(u_b,v)f_b(u_a,u_b,v)\,dx\notag\\
&\quad+\int_{\Omega}\pa_{2}h_{a,p}(u_a,v)f_v(u_a,u_b,v)\,dx+\int_{\Omega}\pa_{2}h_{b,p}(u_b,v)f_v(u_a,u_b,v)\,dx\notag\\
&:=J_1^{p}+J_2^{p}+J_3^{p}+J_4^{p}.
\end{align}

The competitive expression of the reaction functions $f_a$ and $f_b$ enable us to obtain estimates of $\|u_a\|_{L^{q(p)+1}(\Omega_T)}$ and $\|u_b\|_{L^{r(p)+1}(\Omega_T)}$ from $J_1^{p}$ and $J_2^{p}$, respectively (see \eqref{est J1p fin} and \eqref{est J2p}) and to absorbe some terms arising from the diffusion. Indeed, using \eqref{def pa1ha pa1hb}, the definitions of $f_a$ and $\theta$ in \eqref{def react functions}, \eqref{def theta omega}, assumption \eqref{(H3)} and neglecting the non-negative intra competition term $u_au_b$ in $f_a$, it holds
\begin{align}\label{est J1p}
J_1^{p}&=\int_{\Omega}\pa_{1}h_{a,p}(u_a,v)f_a(u_a,u_b,v)\,dx=\int_\Omega\theta(u_a,v)^{\alpha(p-1)}u_a^{p-1}\,f_a(u_a,u_b,v)\,dx\notag\\
&\le\eta_a\int_\Omega\theta(u_a,v)^{\alpha(p-1)}u_a^{p-1}\,u_a(1+A-(A+\afast u_a+\cfast v))\,dx\notag\\
&=\eta_a\int_\Omega u_a^{p}\left[\theta(u_a,v)^{\alpha(p-1)}(1+A-\theta(u_a,v))\right]\,dx.
\end{align}
Next, taking in \eqref{elem. ineq}, $\gamma_1=\gamma+1$ and $C>1$, we have the inequality
\beq\label{dis1}
x^\gamma(1-x)\le C^\gamma-(1-C^{-1})\,x^{\gamma+1}\,.
\eeq
Applying \eqref{dis1} with $x=\f{\theta(u_a,v)}{1+A}$, $\gamma=\alpha(p-1)$ and $C=1+A>1$, to the function in the square brackets in \eqref{est J1p}, we have
\[
\theta(u_a,v)^{\alpha(p-1)}(1+A-\theta(u_a,v))=(1+A)^{\gamma+1}\,x^\gamma(1-x)\le C^{2\gamma+1}-\f {A}{1+A}\,\theta(u_a,v)^{\gamma+1}\,.
\]
Then using \eqref{ineq theta omega} and \eqref{def q r},  we end up with
\beq\label{est J1p fin}
J_1^p\lesssim (1+A)^{2\alpha(p-1)+1}\,\|u_a\|_{L^p(\Omega)}^p-\afast^{\alpha(p-1)+1}\f {A}{1+A}\|u_a\|_{L^{q(p)+1}(\Omega)}^{q(p)+1}\,.
\eeq

Similarly, for $J_2^p$ we obtain that 
\begin{align*}
J_2^p&=\int_\Omega\omega(u_b,v)^{\beta(p-1)}u_b^{p-1}\,f_b(u_a,u_b,v)\,dx\notag\\
&\le\eta_b\int_\Omega u_b^{p}\left[\omega(u_b,v)^{\beta(p-1)}(1+B-\omega(u_b,v))\right]\,dx\,.
\end{align*}
However, as $B\ge0$ when $\beta\ge1$, we choose an arbitrary $\sigma>0$, we replace $B$ with $B\vee\sigma$ and, proceeding as before, we have
\beq\label{est J2p}
J_2^p\lesssim (1+B\vee\sigma)^{2\beta(p-1)+1}\|u_b\|_{L^p(\Omega)}^p-\bfast^{\beta(p-1)+1}\f{B\vee\sigma}{1+B\vee\sigma}
\|u_b\|_{L^{r(p)+1}(\Omega)}^{r(p)+1}\,.
\eeq

The terms $J_3^{p}$ and $J_4^{p}$ in \eqref{Ip rea} cannot give similar estimates since they contain the interaction between $\theta$ and $\omega$, through the reaction function $f_v$, so that \eqref{dis1} can not be applied. Therefore, we simply use \eqref{est pa2 hap hbp} and neglect all the negative terms in $f_v$, to end up with
\beq\label{est J3p}
J_3^p=\int_{\Omega}\pa_{2}h_{a,p}(u_a,v)f_v(u_a,u_b,v)\,dx\lesssim\afast^{-p}\int_\Omega v\,\theta(u_a,v)^{q(p)-1}\,dx\,,
\eeq
and
\beq\label{est J4p}
J_4^p=\int_{\Omega}\pa_{2}h_{b,p}(u_b,v)f_v(u_a,u_b,v)\,dx\lesssim\bfast^{-p}\int_\Omega v\,\omega(u_b,v)^{r(p)-1}\,dx\,.
\eeq

Finally, plugging \eqref{est J1p fin}--\eqref{est J4p} into \eqref{Ip rea},  we have for all $p>1$
\begin{align}\label{est Iprea}
I^{p}_{\text{rea}}&\lesssim (1+A)^{2\alpha(p-1)+1}\,\|u_a\|_{L^p(\Omega)}^p-\afast^{\alpha(p-1)+1}\f {A}{1+A}\|u_a\|_{L^{q(p)+1}(\Omega)}^{q(p)+1}\notag\\
&+(1+B\vee\sigma)^{2\beta(p-1)+1}\|u_b\|_{L^p(\Omega)}^p-\bfast^{\beta(p-1)+1}\f{B\vee\sigma}{1+B\vee\sigma}\|u_b\|_{L^{r(p)+1}(\Omega)}^{r(p)+1}\notag\\
&+ \afast^{-p}\| v\,\theta(u_a,v)^{q(p)-1}\|_{L^1(\Omega)}+\bfast^{-p}\| v\,\omega(u_b,v)^{r(p)-1}\|_{L^1(\Omega)}.
\end{align}
\section{The diffusion contribution $I^p_{\text{diff}}$ to~$\Ep'(t)$  in \eqref{evol energy}}\label{Section Ipdiff}
From \eqref{evol energy} we have 
\[
I^p_{\text{diff}}=d_a\int_{\Omega}\pa_1h_p\Delta u_a\,dx+d_b\int_\Omega\pa_2h_p\Delta u_b\,dx+d_v\int_\Omega\pa_3h_p\Delta v\,dx\,,
\]
and, by definition \eqref{def hp},
\beq\label{Ipdiff}
\begin{split}
I^p_{\text{diff}}&=d_a\int_{\Omega}\pa_1h_{a,p}\Delta u_a\,dx+d_b\int_\Omega\pa_1h_{b,p}\Delta u_b\,dx\\
&\qquad\qquad+d_v\int_\Omega\pa_2h_{a,p}\Delta v\,dx+d_v\int_\Omega\pa_2h_{b,p}\Delta v\,dx\,.
\end{split}
\eeq

As it is not possible to have a priori estimates on $\Delta u_a$ and $\Delta u_b$ uniform in $\e$, we have to apply Green's formula to the first and second integral in the right hand side of \eqref{Ipdiff}. Assumption  \eqref{(H1)} appears to be fundamental here, since we need to control $\pa_{11}h_{a,p}$ and $\pa_{12}h_{a,p}$, for $p<p_\alpha$ (see \eqref{hess hap} and \eqref{sub critical case}). 

On the other hand, as we do not have assumed the strict positivity of the transition function $\phi$ when $\beta\ge1$, we cannot bound $\pa_{22}h_{b,p}$ for $p=p_\beta$, when $\beta\ge1$, (see Subsection \ref{Hess(hbp)}). Hence, the third and forth integral in the right hand side of \eqref{Ipdiff} are left as they are for $p\le p_\alpha$ and $p\le p_\beta$, respectively, and the bound on $\Delta v$ used (see \eqref{max reg p=2}, \eqref{lap v}). Green's formula will be used for these two terms when $p>p_\alpha$ and $p>p_\beta$, respectively. 

To resume, using the boundary conditions \eqref{HNBC}, and the Heaviside functions 
\beq\label{def Heaviside}
\chi_\beta(p):=\left\{
\begin{split}
1\,,\qquad p>p_\beta\,,\\
0\,,\qquad p\le p_\beta\,,
\end{split}
\right.
\qquad 
\chi_\alpha(p):=\left\{
\begin{split}
1\,,\qquad p>p_\alpha\,,\\
0\,,\qquad p\le p_\alpha\,,
\end{split}
\right.
\eeq
$I^p_{\text{diff}}$ rewrites as
\beq\label{Ipdiff 1}
\begin{split}
I^p_{\text{diff}}=&-d_a\int_\Omega\pa_{11}h_{a,p}|\nabla u_a|^2\,dx-(d_a+\chi_\alpha(p)d_v)\int_\Omega\pa_{12}h_{a,p}\nabla u_a\cdot\nabla v\, dx\\
&-d_b\int_\Omega\pa_{11}h_{b,p}|\nabla u_b|^2\, dx-(d_b+\chi_\beta(p)d_v)\int_\Omega\pa_{12}h_{b,p}\nabla u_b\cdot\nabla v\, dx\\
&-d_v\chi_\alpha(p)\int_\Omega\pa_{22}h_{a,p}|\nabla v|^2\, dx
-d_v\chi_\beta(p)\int_\Omega\pa_{22}h_{b,p}|\nabla v|^2\, dx\\
&+d_v(1-\chi_\alpha(p))\int_\Omega\pa_2h_{a,p}\Delta v\,dx
+d_v(1-\chi_\beta(p))\int_\Omega\pa_2h_{b,p}\Delta v\,dx\\
:=& K_1^p+K_2^p+K_3^p+K_4^p+K_5^p+K_6^p\,.
\end{split}
\eeq

In the rest of the section, we will estimate each of the $K_i^p$ terms above. \\

\noindent{\bf Estimate of $K_1^p$.} From \eqref{hess hap}, $K_1^p$ reads as
\beq\label{def K1p}
\begin{split}
K_1^p=&-d_a\int_\Omega\pa_{11}h_{a,p}|\nabla u_a|^2\,dx-(d_a+\chi_\alpha(p)d_v)\int_\Omega\pa_{12}h_{a,p}\nabla u_a\cdot\nabla v\, dx\\
=&-d_a\,\afast\alpha(p-1)\int_\Omega\theta(u_a,v)^{\alpha(p-1)-1}u_a^{p-1}|\nabla u_a|^2\, dx\\
&-d_a(p-1)\int_\Omega\theta(u_a,v)^{\alpha(p-1)}u_a^{p-2}|\nabla u_a|^2\, dx\\
&-(d_a+\chi_\alpha(p)d_v)\,\cfast\alpha(p-1)\int_\Omega \theta(u_a,v)^{\alpha(p-1)-1}u_a^{p-1}    \nabla u_a\cdot\nabla v\, dx\,.
\end{split}
\eeq
Then, by Young's inequality applied to the third integral in the right hand side of~\eqref{def K1p}, we have
\beq\label{est K1p}
\begin{split}
K_1^p\lesssim&-\alpha(p-1)\int_\Omega\theta(u_a,v)^{\alpha(p-1)-1}u_a^{p-1}|\nabla u_a|^2\, dx\\
&-d_a(p-1)\int_\Omega\theta(u_a,v)^{\alpha(p-1)}u_a^{p-2}|\nabla u_a|^2\, dx\\
&+\alpha(p-1)\int_\Omega \theta(u_a,v)^{\alpha(p-1)-1}u_a^{p-1}|\nabla v|^2\,dx\,.
\end{split}
\eeq
Next, we use \eqref{ineq theta omega} in the second integral in the right hand side of \eqref{est K1p} to have 
\[
\theta(u_a,v)^{\alpha(p-1)}u_a^{p-2}\ge\afast^{\alpha(p-1)}\,u_a^{(\alpha+1)(p-1)-1}=\afast^{\alpha(p-1)}\,u_a^{q(p)-2}\,.
\] 
Hence, neglecting the first integral, we obtain
\beq\label{est K1p bis}
\begin{split}
K_1^p\lesssim&-d_a(p-1)\afast^{\alpha(p-1)}\int_\Omega u_a^{q(p)-2}|\nabla u_a|^2\, dx\\
&+\alpha(p-1)\int_\Omega \theta(u_a,v)^{\alpha(p-1)-1}u_a^{p-1}|\nabla v|^2\,dx\,.
\end{split}
\eeq

Finally, for the the second integral in the right hand side of \eqref{est K1p bis}, we proceed according the value of $p$\\ 
$(i)$ if $1<p<p_\alpha$, we employ \eqref{sub critical case} to obtain
\beq\label{K1p 1}
K_1^p\lesssim-d_a\f{4(p-1)}{q(p)^2}\,\afast^{\alpha(p-1)}\|\nabla u_a^{q(p)/2}\|_{L^2(\Omega)}^2+ \alpha(p-1)C_A(p)\|\nabla v\|_{L^2(\Omega)}^2\,;
\eeq
$(ii)$ if $p\ge p_\alpha$, we employ \eqref{critical+supercritical case} and we have
\beq\label{K1p 2}
K_1^p\lesssim -d_a\f{4(p-1)}{q(p)^2}\,\afast^{\alpha(p-1)}\|\nabla u_a^{q(p)/2}\|_{L^2(\Omega)}^2+\alpha(p-1)\,\afast^{1-p}\|\theta(u_a,v)^{q(p)/2-1}\,\nabla v\|_{L^2(\Omega)}^2\,.
\eeq

\noindent{\bf Estimate of $K_2^p$.} From \eqref{hess hbp}, $K_2^p$ reads as
\[
\begin{split}
K_2^p=&-d_b\int_\Omega\pa_{11}h_{b,p}|\nabla u_b|^2\, dx-(d_b+\chi_\beta(p)d_v)\int_\Omega\pa_{12}h_{b,p}\nabla u_b\cdot\nabla v\, dx\\
=&-d_b\,\bfast\beta(p-1)\int_\Omega\omega(u_b,v)^{\beta(p-1)-1}u_b^{p-1}|\nabla u_b|^2\, dx\\
&-d_b(p-1)\int_\Omega\omega(u_b,v)^{\beta(p-1)}u_b^{p-2}|\nabla u_b|^2\, dx\\
&-(d_b+\chi_\beta(p)d_v)\,\dfast\beta(p-1)\int_\Omega \omega(u_b,v)^{\beta(p-1)-1}u_b^{p-1} \nabla u_b\cdot\nabla v\, dx\,.
\end{split}
\]
Again by Young's inequality we have
\[
\begin{split}
K_2^p\lesssim&-\beta(p-1)\int_\Omega\omega(u_b,v)^{\beta(p-1)-1}u_b^{p-1}|\nabla u_b|^2\, dx\\
&-d_b(p-1)\int_\Omega\omega(u_b,v)^{\beta(p-1)}u_b^{p-2}|\nabla u_b|^2\, dx\\
&+\beta(p-1)\int_\Omega \omega(u_b,v)^{\beta(p-1)-1}u_b^{p-1}|\nabla v|^2\, dx\,.
\end{split}
\]
Then, as before, we neglect the first integral, we use \eqref{ineq theta omega} in the second integral and \eqref{critical+supercritical case 2} in the third one, to end up with, for all $p\ge p_\beta$, 
\beq\label{K2p}
K_2^p\lesssim -d_b\f{4(p-1)}{r(p)^2}\,\bfast^{\beta(p-1)}\|\nabla u_b^{r(p)/2}\|_{L^2(\Omega)}^2+\beta(p-1)\,\bfast^{1-p}\|\omega(u_b,v)^{r(p)/2-1}\,\nabla v\|_{L^2(\Omega)}^2\,.
\eeq

\noindent{\bf Estimate of $K_3^p$ and $K_4^p$.} It is sufficient to estimate $K_3^p$ for $p_\alpha< p<1+\alpha^{-1}$, since out of this range of $p$, $K_3^p$ is either zero or  nonpositive (see \eqref{hess hap}) and can be neglected.  From \eqref{est pa22 hap} it holds 
\beq\label{K3p}
\begin{split}
0&<K_3^p=-d_v\chi_\alpha(p)\int_\Omega\pa_{22}h_{a,p}|\nabla v|^2\, dx\\
&\lesssim \alpha(p-1)\,C(\alpha,p)\,\afast^{-p}\|\theta(u_a,v)^{q(p)/2-1}\nabla v\|_{L^2(\Omega)}^2\,,\quad p_\alpha< p<1+\alpha^{-1}.
\end{split}
\eeq

Similarly, $K_4^p$ is strictly positive for $p_\beta< p<1+\beta^{-1}$ (see \eqref{hess hbp}), and in this range, by \eqref{est pa22 hbp}, it holds
\beq\label{K4p}
\begin{split}
0&<K_4^p=-d_v\chi_\beta(p)\int_\Omega\pa_{22}h_{b,p}|\nabla v|^2\, dx\\
&\lesssim \beta(p-1)\,C(\beta,p)\,\bfast^{-p}\|\omega(u_b,v)^{r(p)/2-1}\nabla v\|_{L^2(\Omega)}^2\,,\quad p_\beta< p<1+\beta^{-1}.
\end{split}
\eeq

\noindent{\bf Estimate of $K_5^p$ and $K_6^p$.}  The term $K_5^p$ is not trivial for $p\le p_\alpha$ and it will be employed only for $p=p_\beta$ and $p=p_\alpha$. Thus, using \eqref{est pa2 hap hbp} and \eqref{def q r}, we obtain
\beq\label{K5pbeta}
\begin{split}
0<K_5^{p_\beta}&=d_v\int_\Omega\pa_2h_{a,p_\beta}\Delta v\,dx\lesssim \afast^{-p_\beta} \int_\Omega \theta(u_a,v)^{\f{\alpha+1}{\beta+1}}|\Delta v|\,dx\\
&\lesssim \|\theta(u_a,v)^{q(p_\beta)-1}\|_{L^2(\Omega)}\|\Delta v\|_{L^2(\Omega)}\,,
\end{split}
\eeq
and
\beq\label{K5palpha}
0<K_5^{p_\alpha}\lesssim\|\theta(u_a,v)\|_{L^2(\Omega)}\|\Delta v\|_{L^2(\Omega)}\,.
\eeq

Finally, the term $K_6^p$ is not trivial for $p\le p_\beta$ and it will be employed only for $p=p_\beta$. Using \eqref{est pa2 hap hbp}, we have  
\beq\label{K6pbeta}
0< K_6^{p_\beta}=d_v\int_\Omega\pa_2h_{b,p_\beta}\Delta v\,dx\lesssim \|\omega(u_b,v)\|_{L^2(\Omega)}\|\Delta v\|_{L^2(\Omega)}.
\eeq
\section{Energy estimates: proof of Lemma \ref{energy estimate}}\label{sect proof energy lemma}
This section is devoted to the proof of Lemma \ref{energy estimate}, based on the computations obtained in Sections \ref{energy} to \ref{Section Ipdiff}, on the maximal regularity \eqref{max reg} giving \eqref{est gradv} and on a bootstrap argument. We recall that $\alpha\le\beta$, so that $p_\beta\le p_\alpha$, and that $q(p_\alpha)=r(p_\beta)=2$. Therefore, we estimate $\mathcal{E}_p$ along the  trajectories of the solution $(u_a^\e,u_b^\e,v^\e)$, starting with $p={p_\beta}$, then $p={p_\alpha}$ and finally $p\in[2,p^0_{\alpha,\beta})$ (see \eqref{def palphabeta}), using the differential equation below (see \eqref{evol energy} and~\eqref{Ipdiff 1})
\beq\label{est Ep 0}
\mathcal{E}'_{p}(t)= I^{p}_{\text{diff}}+I^{p}_{\text{rea}}+I^{p}_{\text{fast}}=\sum_{i=1}^6 K_i^{p}+I^{p}_{\text{rea}}+I^{p}_{\text{fast}}\,.
\eeq
%
\subsection{Estimates from $\mathcal{E}_{p_\beta}$, $\alpha<\beta$}\label{Epbeta}
Assume $\alpha<\beta$. Taking into account that $K_3^{p_\beta}=K_4^{p_\beta}=0$ (see \eqref{K3p},\eqref{K4p}), from estimates \eqref{K1p 1}, \eqref{K2p} where we use $r(p_\beta)=2$, \eqref{K5pbeta} and \eqref{K6pbeta}, we have 
\beq\label{est Ipbetadiff}
\begin{split}
I^{p_\beta}_{\text{diff}}\lesssim&-\|\nabla u_a^{q(p_\beta)/2}\|_{L^2(\Omega)}^2-\|\nabla u_b\|_{L^2(\Omega)}^2+\|\nabla v\|_{L^2(\Omega)}^2\\
&+\|\theta(u_a,v)^{q(p_\beta)-1}\|_{L^2(\Omega)}\|\Delta v\|_{L^2(\Omega)}+\|\omega(u_b,v)\|_{L^2(\Omega)}\|\Delta v\|_{L^2(\Omega)}\,.
\end{split}
\eeq
Using \eqref{Ipfats neg}, plugging \eqref{est Ipbetadiff} and estimate \eqref{est Iprea} of $I^{p_\beta}_{\text{rea}}$ into \eqref{est Ep 0}, rearranging the terms and integrating the obtained  inequality over $(0,T)$, we get 
\beq\label{est Epbeta}
\begin{split}
&\mathcal{E}_{p_\beta}(T)-\mathcal{E}_{p_\beta}(0)+\|\nabla u_a^{q(p_\beta)/2}\|_{L^2(\Omega_T)}^2+\|\nabla u_b\|_{L^2(\Omega_T)}^2+\|u_a\|_{L^{q(p_\beta)+1}(\Omega_T)}^{q(p_\beta)+1}+\|u_b\|_{L^3(\Omega_T)}^3\\
&\lesssim\ \|\nabla v\|_{L^2(\Omega_T)}^2+\|\Delta v\|^2_{L^2(\Omega_T)}+\|\theta(u_a,v)^{q(p_\beta)-1}\|^2_{L^2(\Omega_T)}+\|\omega(u_b,v)\|^2_{L^2(\Omega_T)}\\
&+\|u_a\|_{L^{p_\beta}(\Omega_T)}^{p_\beta}+\|u_b\|_{L^{p_\beta}(\Omega_T)}^{p_\beta}
+\| v\,\theta(u_a,v)^{q(p_\beta)-1}\|_{L^1(\Omega_T)}+\|v\,\omega(u_b,v)\|_{L^1(\Omega_T)}\,.
\end{split}
\eeq
As $p_\beta<2$ and $q(p_\beta)-1=\f{\alpha+1}{\beta+1}<1$, recalling that $\theta$ and $\omega$ are affine functions, the estimates obtained in Theorem~\ref{th exist meso}  allow us to control all the terms in the right hand side of \eqref{est Epbeta}. Hence, using \eqref{est from E}, we end up with 
\begin{align}\label{est from E 1}
&\|u_a\|^{q(p_\beta)}_{L^\infty(0,T;L^{q(p_\beta)}(\Omega))}+\|u_b\|^2_{L^\infty(0,T;L^{2}(\Omega))}+\|\nabla u_a^{q(p_\beta)/2}\|_{L^2(\Omega_T)}^2+\|\nabla u_b\|_{L^2(\Omega_T)}^2\notag\\
&\qquad+\|u_a\|_{L^{q(p_\beta)+1}(\Omega_T)}^{q(p_\beta)+1}+\|u_b\|_{L^3(\Omega_T)}^3\lesssim \mathcal{E}_{p_\beta}(0)+C(|\Omega|,T)\,.
\end{align} 
\subsection{Estimates from $\mathcal{E}_{p_\alpha}$, $\alpha<\beta$}\label{Epalpha}
Assume $\alpha<\beta$. Taking into account that $K_3^{p_\alpha}=K_6^{p_\alpha}=0$, using \eqref{K1p 2} with $q(p_\alpha)=2$, \eqref{K2p}, \eqref{K4p}, \eqref{K5palpha}, we have
\beq\label{est Ipalphadiff}
\begin{split}
I^{p_\alpha}_{\text{diff}}\lesssim&-\|\nabla u_a\|_{L^2(\Omega)}^2-\|\nabla u_b^{r(p_\alpha)/2}\|_{L^2(\Omega)}^2+\|\nabla v\|_{L^2(\Omega)}^2\\
&+\|\omega(u_b,v)^{r(p_\alpha)/2-1}\,\nabla v\|_{L^2(\Omega)}^2+\|\theta(u_a,v)\|_{L^2(\Omega)}\|\Delta v\|_{L^2(\Omega)}\,.
\end{split}
\eeq
Using \eqref{Ipfats neg}, plugging \eqref{est Ipalphadiff} and estimate \eqref{est Iprea} of $I^{p_\alpha}_{\text{rea}}$ into \eqref{est Ep 0}, rearranging the terms and integrating the obtained  inequality over $(0,T)$, we get
\beq\label{est Epalpha}
\begin{split}
&\mathcal{E}_{p_\alpha}(T)-\mathcal{E}_{p_\alpha}(0)+\|\nabla u_a\|_{L^2(\Omega_T)}^2+\|\nabla u_b^{r(p_\alpha)/2}\|_{L^2(\Omega_T)}^2+\|u_a\|_{L^{3}(\Omega_T)}^{3}+\|u_b\|_{L^{r(p_\alpha)+1}(\Omega_T)}^{r(p_\alpha)+1}\\
&\lesssim\ \|\nabla v\|_{L^2(\Omega_T)}^2+\|\omega(u_b,v)^{r(p_\alpha)/2-1}\nabla v\|_{L^2(\Omega_T)}^2
+\|\theta(u_a,v)\|^2_{L^2(\Omega_T)}+\|\Delta v\|^2_{L^2(\Omega_T)}\\
&+\|u_a\|_{L^{p_\alpha}(\Omega_T)}^{p_\alpha}+\|u_b\|_{L^{p_\alpha}(\Omega_T)}^{p_\alpha}
+\|v\,\theta(u_a,v)\|_{L^1(\Omega_T)}+\|v\,\omega(u_b,v)^{r(p_\alpha)-1}\|_{L^1(\Omega_T)}\,.
\end{split}
\eeq

It is worth noticing that despite $p_\alpha<2$, $r(p_\alpha)$ can be large without any restriction on $\beta-\alpha$, since from \eqref{def q r},\eqref{def palpha pbeta} and $\alpha<\beta$, it holds
\beq\label{rpalpha}
r(p_\alpha)=\f{\beta+1}{\alpha+1}+1=\f{\beta-\alpha}{\alpha+1}+2>2\,.
\eeq
Hence, in order to obtain new a priori estimates on $u_a$ and $u_b$ from  \eqref{est Epalpha},  we need to get rid of the terms
\beq\label{def I J}
I:=\|\omega(u_b,v)^{r(p_\alpha)/2-1}\nabla v\|^2_{L^2(\Omega_T)}\quad\text{and}\quad
J:=\|v\,\omega(u_b,v)^{r(p_\alpha)-1}\|_{L^1(\Omega_T)}\,.
\eeq

Let $\delta\in(0,1)$. Applying Young's inequality into $I$ in \eqref{def I J} we have
\beq\label{I=I1+I2}
I\lesssim {\delta}\int_{\Omega_T}|\nabla v|^6\,dx\,dt+\delta^{-\f12}\int_{\Omega_T}\omega(u_b,v)^{\f32(r(p_\alpha)-2)}\,dx\,dt=:I_1+I_2\,.
\eeq
Then, by \eqref{est gradv}  with $p=3$,
\beq\label{I1 alpha<beta}
I_1={\delta}\|\nabla v\|^6_{L^6(\Omega_T)}\lesssim {\delta}(1+T)+{\delta}\,\|u_a\|_{L^{3}(\Omega_T)}^{3}+{\delta}\,\|u_b\|_{L^{3}(\Omega_T)}^{3}\,.
\eeq
On the other hand, by \eqref{rpalpha} and assumption \eqref{(H2)},  
\[
0<\f32(r(p_\alpha)-2)=\f32\f{\beta-\alpha}{\alpha+1}<\f{\beta-\alpha}{\alpha+1}+3=r(p_\alpha)+1\,.
\]
Using \eqref{elem. ineq} with  $\gamma=\f32(r(p_\alpha)-2)$, $\gamma_1=r(p_\alpha)+1$ and $C=C(\delta)>0$ such that $\delta^{-\f12}\,C(\delta)^{\gamma-\gamma_1}=\delta$, 
we obtain 
\[
I_2=\delta^{-\f12}\!\int_{\Omega_T} \omega(u_b,v)^{\f32(r(p_\alpha)-2)}\,dxdt
\le \delta C(\delta)^{r(p_\alpha)+1}|\Omega_T|+\delta\!\int_{\Omega_T}\omega(u_b,v)^{r(p_\alpha)+1}dxdt.
\]
Finally, recalling definition \eqref{def theta omega} of $\omega$, by Jensen's inequality and the boundedness of $v$ (see \eqref{bound v}), we have
\beq\label{I2 alpha<beta 3}
I_2\lesssim \delta\,C(\delta)^{r(p_\alpha)+1}|\Omega_T|+\delta|\Omega_T|+\delta\,\| u_b\|^{r(p_\alpha)+1}_{L^{r(p_\alpha)+1}(\Omega_T)}\,.
\eeq
Plugging \eqref{I1 alpha<beta} and \eqref{I2 alpha<beta 3} into \eqref{I=I1+I2}, we obtain
\beq\label{I first}
I\lesssim {\delta}(1+C(\delta,|\Omega|,T))+{\delta}\,\|u_a\|_{L^{3}(\Omega_T)}^{3}+{\delta}\,\|u_b\|_{L^{3}(\Omega_T)}^{3}
+\delta\,\| u_b\|^{r(p_\alpha)+1}_{L^{r(p_\alpha)+1}(\Omega_T)}\,.
\eeq

Next, by \eqref{rpalpha}, $r(p_\alpha)-1>1$. Hence, using again Jensen's inequality and the boundedness of $v$ in $J$ defined in \eqref{def I J}, we have
\beq\label{est J 0}
J=\int_{\Omega_T}v\,\omega(u_b,v)^{r(p_\alpha)-1}\,dx\,dt\lesssim \int_{\Omega_T} u_b^{r(p_\alpha)-1}\,dx\,dt+|\Omega_T|\,.
\eeq
Using \eqref{elem. ineq} with $\gamma\!=\!r(p_\alpha)-1$, $\gamma_1\!=\!r(p_\alpha)+1$ and $C\!=\!\delta^{-\f12}$ so that $C^{\gamma-\gamma_1}\!=\!\delta$, we end up with
\beq\label{est J}
J\lesssim \delta^{\f{1-r(p_\alpha)}2}|\Omega_T|+\delta\,\| u_b\|^{r(p_\alpha)+1}_{L^{r(p_\alpha)+1}(\Omega_T)}+|\Omega_T|\,.
\eeq

To conclude, \eqref{I first} and \eqref{est J} imply that there exists $C(\delta,|\Omega|,T)>0$ such that
\beq\label{est I+J}
I+J\lesssim C(\delta,|\Omega|,T)+{\delta}\,\|u_a\|_{L^{3}(\Omega_T)}^{3}+\delta\,\|u_b\|_{L^{3}(\Omega_T)}^{3}+\delta\,\| u_b\|^{r(p_\alpha)+1}_{L^{r(p_\alpha)+1}(\Omega_T)}\,.
\eeq
Choosing $\delta$ small enough, plugging \eqref{est I+J} into \eqref{est Epalpha} and rearranging the terms, we obtain
\beq\label{est Epalpha final}
\begin{split}
\mathcal{E}_{p_\alpha}(T)&-\mathcal{E}_{p_\alpha}(0)+\|\nabla u_a\|_{L^2(\Omega_T)}^2+\|\nabla u_b^{r(p_\alpha)/2}\|_{L^2(\Omega_T)}^2\\
&\qquad\qquad+(1-\delta)\|u_a\|_{L^{3}(\Omega_T)}^{3}+(1-\delta)\|u_b\|_{L^{r(p_\alpha)+1}(\Omega_T)}^{r(p_\alpha)+1}\\
\lesssim&\ C(\delta,|\Omega|,T)+\|\nabla v\|_{L^2(\Omega_T)}^2+\|\theta(u_a,v)\|^2_{L^2(\Omega_T)}+\|\Delta v\|^2_{L^2(\Omega_T)}\\
&+\|u_a\|_{L^{p_\alpha}(\Omega_T)}^{p_\alpha}+\|u_b\|_{L^{p_\alpha}(\Omega_T)}^{p_\alpha}+\|v\,\theta(u_a,v)\|_{L^1(\Omega_T)}+\delta\,\|u_b\|_{L^{3}(\Omega_T)}^{3}\,.
\end{split}
\eeq
Recalling that $p_\alpha<2$, the estimates obtained in Theorem~\ref{th exist meso} plus the estimate of $\|u_b\|_{L^{3}(\Omega_T)}$ obtained in \eqref{est from E 1}, will allow us to control all the terms in the right hand side of \eqref{est Epalpha final}, so that, using \eqref{est from E}, we get
\begin{align}\label{est from E 2}
&\|u_a\|^2_{L^\infty(0,T;L^{2}(\Omega))}+\|u_b\|^{r(p_\alpha)}_{L^\infty(0,T;L^{r(p_\alpha)}(\Omega))}+\|\nabla u_a\|_{L^2(\Omega_T)}^2+\|\nabla u_b^{r(p_\alpha)/2}\|_{L^2(\Omega_T)}^2\notag \\
&\qquad+\|u_a\|_{L^{3}(\Omega_T)}^{3}+\|u_b\|_{L^{r(p_\alpha)+1}(\Omega_T)}^{r(p_\alpha)+1}\lesssim \mathcal{E}_{p_\alpha}(0)+C(|\Omega|,T)\,.
\end{align}
\subsection{Estimates from $\mathcal{E}_{p_\alpha}=\mathcal{E}_{p_\beta}$, $\alpha=\beta$}\label{Ep alpha=Beta}
If $\alpha=\beta$, then $p_\alpha=p_\beta$, $K_3^{p_\alpha}=K_4^{p_\alpha}=0$ and we can use \eqref{K1p 2}, \eqref{K2p}, \eqref{K5palpha}, \eqref{K6pbeta}, to obtain
\beq\label{est Ipalpha=beta diff}
\begin{split}
I^{p_\beta}_{\text{diff}}
&\lesssim-\|\nabla u_a\|_{L^2(\Omega)}^2-\|\nabla u_b\|_{L^2(\Omega)}^2+\|\nabla v\|_{L^2(\Omega)}^2\\
&\quad+\|\theta(u_a,v)\|_{L^2(\Omega)}\|\Delta v\|_{L^2(\Omega)}+\|\omega(u_b,v)\|_{L^2(\Omega)}\|\Delta v\|_{L^2(\Omega)}\,.
\end{split}
\eeq
Employing again \eqref{Ipfats neg} and estimate \eqref{est Iprea} for $I^{p_\beta}_{\text{rea}}=I^{p_\alpha}_{\text{rea}}$, plugging \eqref{est Ipalpha=beta diff} into \eqref{est Ep 0}  and integrating the obtained inequality over $(0,T)$, we get 
\begin{align}\label{energy est alpha=beta}
&\mathcal{E}_{p_\beta=p_\alpha}(T)-\mathcal{E}_{p_\beta=p_\alpha}(0)+\|\nabla u_a\|_{L^2(\Omega_T)}^2+\|\nabla u_b\|_{L^2(\Omega_T)}^2+\|u_a\|_{L^{3}(\Omega_T)}^{3}+\|u_b\|_{L^3(\Omega_T)}^3\notag\\
&\quad\lesssim \|\nabla v\|_{L^2(\Omega_T)}^2+\|\Delta v\|_{L^2(\Omega_T)}^2+\|\theta(u_a,v)\|_{L^2(\Omega_T)}^2+\|\omega(u_b,v)\|_{L^2(\Omega_T)}^2\\
&\qquad+\|u_a\|_{L^{p_\alpha}(\Omega_T)}^{p_\alpha}+\|u_b\|_{L^{p_\alpha}(\Omega_T)}^{p_\alpha}+\|v(\theta(u_a,v)+\omega(u_b,v))\|_{L^1(\Omega_T)}\notag\,.
\end{align}
Then, we see that all the terms in the right hand side of \eqref{energy est alpha=beta} are controlled by the estimates obtained in Theorem~\ref{th exist meso}, and we get 
\beq\label{est from E 3}
\begin{split}
\|u_a\|^2_{L^\infty(0,T;L^{2}(\Omega))}&+\|u_b\|^2_{L^\infty(0,T;L^{2}(\Omega))}+\|\nabla u_a\|_{L^2(\Omega_T)}^2+\|\nabla u_b\|_{L^2(\Omega_T)}^2\\
&+\|u_a\|_{L^{3}(\Omega_T)}^{3}+\|u_b\|_{L^{3}(\Omega_T)}^{3}\lesssim \mathcal{E}_{p_\beta=p_\alpha}(0)+C(|\Omega|,T)\,.
\end{split}
\eeq 
\subsection{Estimates from $\mathcal{E}_{p}$, $p\in[2,p^0_{\alpha,\beta})$}\label{Energy p>=2)}
To begin with, we fix $p$ in $[2,p^0_{\alpha,\beta})$, where $p^0_{\alpha,\beta}$ is defined below and $p^0_{\alpha,\beta}\!>\!2$ by~\eqref{(H2)}
\beq\label{def palphabeta}
p^0_{\alpha,\beta}:=\left\{
\begin{split}
&1+\f4{\beta-3\alpha-2}\,,\qquad\text{if }2(\alpha+1)<\beta-\alpha<2(\alpha+3)\,,\\
&+\infty\qquad\qquad\qquad\quad\text{if }0\le\beta-\alpha\le2(\alpha+1)\,.
\end{split}
\right.
\eeq
Observing that (because of definition \eqref{def Heaviside}) the terms $K_5^p$ and $K_6^p$ in \eqref{Ipdiff 1} both vanish for $p\ge2$, \eqref{Ipdiff 1} reads as $I^{p}_{\text{diff}}=K_1^p+K_2^p+K_3^p+K_4^p$. Moreover, the term $K_3^p$ in \eqref{K3p} (respectively $K_4^p$ in \eqref{K4p}) gives a positive contribution to $I^{p}_{\text{diff}}$ if and only if $p\in(p_\alpha,1+\f1\alpha)$ (respectively $p\in(p_\beta,1+\f1\beta$)). If $\alpha<1$ (respectively $\beta<1$), there are $p\!\in\![2,1+\f1\alpha)$, (respectively $p\!\in\![2,1+\f1\beta)$), and in that case we use the decreasing character of the constant $C(\alpha,p)$ in \eqref{K3p}, defined in \eqref{est pa22 hap}, to obtain 
\[
0<C(\alpha,p)\le C(\alpha,2)=(1-\alpha)/\alpha\,.
\] 
Hence, the upper bound \eqref{K3p} of $K_3^p$ can be absorbed by the upper bound \eqref{K1p 2} of $K_1^p$ (respectively $0<C(\beta,p)\le(1-\beta)/\beta$ and the upper bound \eqref{K4p} of $K_4^p$ can be absorbed by the upper bound \eqref{K2p} of $K_2^p$). Therefore, by \eqref{K1p 2},\eqref{K2p}, it holds 
\begin{align}
&I^{p}_{\text{diff}}
\lesssim-\f{4(p-1)}{q(p)^2}\,\afast^{\alpha(p-1)}\|\nabla u_a^{\f{q(p)}2}\|_{L^2(\Omega)}^2+(p-1)\,\afast^{1-p}\|\theta(u_a,v)^{\f{q(p)}2-1}\,\nabla v\|_{L^2(\Omega)}^2\notag\\
&-\f{4(p-1)}{r(p)^2}\,\bfast^{\beta(p-1)}\|\nabla u_b^{\f{r(p)}2}\|_{L^2(\Omega)}^2+(p-1)\,\bfast^{1-p}\|\omega(u_b,v)^{\f{r(p)}2-1}\,\nabla v\|_{L^2(\Omega)}^2\label{est Ipdiff}.
\end{align}
Plugging \eqref{est Ipdiff} and the estimate \eqref{est Iprea} of $I^{p}_{\text{rea}}$ into \eqref{est Ep 0}, rearranging the terms and integrating the obtained inequality over $(0,T)$, we end up with
\[
\begin{split}
&\mathcal{E}_{p}(T)-\mathcal{E}_{p}(0)+\f{4(p-1)}{q(p)^2}\,\afast^{\alpha(p-1)}\|\nabla u_a^{\f{q(p)}2}\|_{L^2(\Omega_T)}^2
+\f{4(p-1)}{r(p)^2}\,\bfast^{\beta(p-1)}\|\nabla u_b^{\f{r(p)}2}\|_{L^2(\Omega_T)}^2\\
&+\afast^{\alpha(p-1)+1}\f {A}{1+A}\|u_a\|_{L^{q(p)+1}(\Omega_T)}^{q(p)+1}
+\bfast^{\beta(p-1)+1}\f{B\vee\sigma}{1+B\vee\sigma}\|u_b\|_{L^{r(p)+1}(\Omega_T)}^{r(p)+1}\\
&-\int_0^TI^{p}_{\text{fast}}(t)\,dt\\
&\qquad\lesssim(p-1)\,\afast^{1-p}\|\theta(u_a,v)^{\f{q(p)}2-1}\,\nabla v\|_{L^2(\Omega_T)}^2\\
&\qquad+(p-1)\,\bfast^{1-p}\|\omega(u_b,v)^{\f{r(p)}2-1}\,\nabla v\|_{L^2(\Omega_T)}^2\\
&\qquad +(1+A)^{2\alpha(p-1)+1}\,\|u_a\|_{L^p(\Omega_T)}^p +(1+B\vee\sigma)^{2\beta(p-1)+1}\|u_b\|_{L^p(\Omega_T)}^p\\
&\qquad+ \afast^{-p}\|v\,\theta(u_a,v)^{q(p)-1}\|_{L^1(\Omega_T)}+\bfast^{-p}\|v\,\omega(u_b,v)^{r(p)-1}\|_{L^1(\Omega_T)}\,.
\end{split}
\]
Next, we observe that, by \eqref{est from E} and the above inequality, for all $p\ge2$, there exists $\mathcal{C}_{p}(\alpha,\beta,A,B,\afast,\bfast)>0$ such that 
\begin{align}\label{evol Ep p>=2}
&\|u_a^\e\|^{q(p)}_{L^\infty(0,T;L^{q(p)}(\Omega))}+\|u_b^\e\|^{r(p)}_{L^\infty(0,T;L^{r(p)}(\Omega))}+\|\nabla u_a^{\f{q(p)}2}\|_{L^2(\Omega_T)}^2+\|\nabla u_b^{\f{r(p)}2}\|_{L^2(\Omega_T)}^2\notag\\
&\qquad+\|u_a\|_{L^{q(p)+1}(\Omega_T)}^{q(p)+1}+\|u_b\|_{L^{r(p)+1}(\Omega_T)}^{r(p)+1}-\int_0^TI^{p}_{\text{fast}}(t)\,dt\notag\\
&\le \mathcal{C}_{p}\|\theta(u_a,v)^{\f{q(p)}2-1}\,\nabla v\|_{L^2(\Omega_T)}^2+ \mathcal{C}_{p}\|\omega(u_b,v)^{\f{r(p)}2-1}\,\nabla v\|_{L^2(\Omega_T)}^2\\
&\qquad+\mathcal{C}_{p}\|u_a\|_{L^p(\Omega_T)}^p+\mathcal{C}_{p}\|u_b\|_{L^p(\Omega_T)}^p\notag\\
&\qquad+\mathcal{C}_{p}\|v\,\theta(u_a,v)^{q(p)-1}\|_{L^1(\Omega_T)}+\mathcal{C}_{p}\|v\,\omega(u_b,v)^{r(p)-1}\|_{L^1(\Omega_T)}+\mathcal{C}_{p}\mathcal{E}_{p}(0)\notag\\[1.5ex]
&:=Z_1^p+Z_2^p+Z_3^p+Z_4^p+Z_5^p+Z_6^p+\mathcal{C}_{p}\mathcal{E}_{p}(0)\notag\,.
\end{align}

We proceed estimating $Z_1^p,Z_2^p$ and $Z_5^p+Z_6^p$ in such a way to be all absorbed by $\|u_a\|_{L^{q(p)+1}(\Omega_T)}^{q(p)+1}$ and $\|u_b\|_{L^{r(p)+1}(\Omega_T)}^{r(p)+1}$ in the left hand side of \eqref{evol Ep p>=2}. As the computations will be in the line of the computations carried out in the case $p=p_\alpha$ with $\alpha<\beta$ (see Subsection \ref{Epalpha}), redundant details will be omitted. Hereafter, the constant $\mathcal{C}_{p}$ will change from line to line and it may depends also on $|\Omega|$,  $T$ and the constants $C_p^\MR,C_p^\ADN,C_p^\GN$ in \eqref{max reg v},\eqref{ADN},\eqref{GN} respectively. 
\medskip

\noindent{\bf Estimate of $Z_1^p$.} Let $\delta\!>\!0$ to be chosen later and note that,  if $p\!\ge\!2$, $q(p)-2\ge\alpha>0$. Applying Young's inequality, we have
\beq\label{est Z1p 1}
\begin{split}
Z_1^p&=\mathcal{C}_{p}\int_{\Omega_T}\theta(u_a,v)^{q(p)-2}\,|\nabla v|^2\,dxdt\\
&\le\delta\,\mathcal{C}_{p}\int_{\Omega_T}\,|\nabla v|^{2(q(p)+1)}\,dxdt+\delta^{-\f1{q(p)}}\,\mathcal{C}_{p}\int_{\Omega_T}\theta(u_a,v)^{(q(p)-2)\f{q(p)+1}{q(p)}}\,dxdt.
\end{split}
\eeq
By \eqref{est gradv} and $q(p)\le r(p)$, we obtain
\beq\label{est Z1p max reg}
\delta\int_{\Omega_T}\,|\nabla v|^{2(q(p)+1)}\,dx\,dt\le\delta\,\mathcal{C}_{p}\left(1+\|u_a\|^{q(p)+1}_{L^{q(p)+1}(\Omega_T)}+\|u_b\|^{r(p)+1}_{L^{r(p)+1}(\Omega_T)}\right)\,.
\eeq
Next, using \eqref{elem. ineq} with $\gamma=(q(p)-2)\f{q(p)+1}{q(p)}$, $\gamma_1=q(p)+1$ and $C=\delta^{-\f12}$ so that 
$\delta^{-\f1{q(p)}}\,C^{\gamma-\gamma_1}=\delta$, we get
\beq\label{est Z1p 2}
\begin{split}
\delta^{-\f1{q(p)}}&\int_{\Omega_T}\theta(u_a,v)^{(q(p)-2)\f{q(p)+1}{q(p)}}\,dx\,dt\le
\delta^{-\f{q(p)-1}2}|\Omega_T|+\delta\int_{\Omega_T}\theta(u_a,v)^{q(p)+1}\,dx\,dt\\
&\le\delta^{-\f{q(p)-1}2}|\Omega_T|+\delta\,\mathcal{C}_{p}\left(1+\|u_a\|^{q(p)+1}_{L^{q(p)+1}(\Omega_T)}\right)\,.
\end{split}
\eeq
Plugging \eqref{est Z1p max reg},\eqref{est Z1p 2} into \eqref{est Z1p 1}, we end up with
\beq\label{est Z1p 3}
Z_1^p\le \mathcal{C}_p\delta(1+\delta^{-\f{q(p)+1}2})+\delta\,\mathcal{C}_p\left(\|u_a\|^{q(p)+1}_{L^{q(p)+1}(\Omega_T)}+\|u_b\|^{r(p)+1}_{L^{r(p)+1}(\Omega_T)}\right)\,.
\eeq

\noindent{\bf Estimate of $Z_2^p$.} Let us observe that $r(p)-2\ge\beta>0$, if $p\ge2$. Hence, proceeding as above,
\begin{align}\label{est Z2p alpha<=beta 1}
Z_2^p&=\mathcal{C}_p\int_{\Omega_T}\omega(u_b,v)^{r(p)-2}\,|\nabla v|^2\,dxdt\notag\\
&\le\delta\,\mathcal{C}_p\int_{\Omega_T}\,|\nabla v|^{2(q(p)+1)}\,dxdt+\delta^{-\f1{q(p)}}\mathcal{C}_p\int_{\Omega_T}\omega(u_b,v)^{(r(p)-2)\f{q(p)+1}{q(p)}}\,dxdt\notag\\
&\le\delta\,\mathcal{C}_p\left(1+\|u_a\|^{q(p)+1}_{L^{q(p)+1}(\Omega_T)}+\|u_b\|^{r(p)+1}_{L^{r(p)+1}(\Omega_T)}\right)\notag\\
&\quad+\delta^{-\f1{q(p)}}\mathcal{C}_p\int_{\Omega_T}\omega(u_b,v)^{(r(p)-2)\f{q(p)+1}{q(p)}}\,dxdt\,.
\end{align}
If $\alpha=\beta$, then $r(p)=q(p)$ and we can proceed as in \eqref{est Z1p 2} to obtain from~\eqref{est Z2p alpha<=beta 1}
\beq\label{est Z2p alpha=beta}
Z_2^p\le \delta\,\mathcal{C}_p(1+\delta^{-\f{r(p)+1}2})+\delta\,\mathcal{C}_p\left(\|u_a\|^{q(p)+1}_{L^{q(p)+1}(\Omega_T)}+\|u_b\|^{r(p)+1}_{L^{r(p)+1}(\Omega_T)}\right)\,.
\eeq
If $\alpha<\beta$, using definitions \eqref{def q r}, as $p\in[2,p^0_{\alpha,\beta})$, it holds
\beq\label{critical inequality}
0<(r(p)-2)\f{q(p)+1}{q(p)}< r(p)+1\,.
\eeq
Therefore, we apply once again \eqref{elem. ineq} with $\gamma=(r(p)-2)\f{q(p)+1}{q(p)}$, $\gamma_1=r(p)+1$ and $C=C(\delta)>0$ such that 
$\delta^{-\f{1}{q(p)}}\,C(\delta)^{\gamma-\gamma_1}=\delta$, to get
\beq\label{est Z2p alpha<beta 1}
\begin{split}
\delta^{-\f1{q(p)}}&\int_{\Omega_T}\omega(u_b,v)^{(r(p)-2)\f{q(p)+1}{q(p)}}\,dxdt
\le\delta\,C(\delta)^{r(p)+1}|\Omega_T|+\delta \int_{\Omega_T}\omega(u_b,v)^{r(p)+1}\,dxdt\\
&\le\delta\,C(\delta)^{r(p)+1}|\Omega_T|+\delta\, \mathcal{C}_p\left(1+\|u_b\|^{r(p)+1}_{L^{r(p)+1}(\Omega_T)}\right)\,.
\end{split}
\eeq
Plugging \eqref{est Z2p alpha<beta 1} into \eqref{est Z2p alpha<=beta 1} and taking into account that $C(\delta)^{r(p)+1}=\delta^{-\f{(q(p)+1)(r(p)+1)}{2+3q(p)-r(p)}}$, we obtain 
\beq\label{est Z2p alpha<beta}
Z_2^p\le \delta\,\mathcal{C}_p(1+\delta^{-\f{(q(p)+1)(r(p)+1)}{2+3q(p)-r(p)}})+\delta\,\mathcal{C}_p\left(\|u_a\|^{q(p)+1}_{L^{q(p)+1}(\Omega_T)}+\|u_b\|^{r(p)+1}_{L^{r(p)+1}(\Omega_T)}\right)\,.
\eeq
Note that \eqref{est Z2p alpha<beta} becomes \eqref{est Z2p alpha=beta} when $r(p)=q(p)$ (i.e. $\alpha=\beta$).
\medskip

\noindent{\bf Estimate of $Z_5^p+Z_6^p$.} Proceeding as in \eqref{est J 0},\eqref{est J}, we have  
\beq\label{est Z5p+Z6p}
\begin{split}
Z_5^p+Z_6^p&=\mathcal{C}_p\int_{\Omega_T} v\,\theta(u_a,v)^{q(p)-1}\,dx\,dt+\mathcal{C}_p\int_{\Omega_T} v\,\omega(u_b,v)^{r(p)-1}\,dx\,dt\\
&\le\mathcal{C}_p\int_{\Omega_T} u_a^{^{q(p)-1}}\,dxdt+\mathcal{C}_p\int_{\Omega_T} u_b^{^{r(p)-1}}\,dxdt+\mathcal{C}_p|\Omega_T|\\
&\le\mathcal{C}_p\delta (\delta^{-\f{q(p)+1}2}+ \delta^{-\f{r(p)+1}2})\\
&\qquad+\delta\,\mathcal{C}_p(\|u_a\|^{q(p)+1}_{L^{q(p)+1}(\Omega_T)}+\|u_b\|^{r(p)+1}_{L^{r(p)+1}(\Omega_T)})+\mathcal{C}_p|\Omega_T|\,.
\end{split}
\eeq

\noindent{\bf Final energy estimate.} 
Plugg \eqref{est Z1p 3}, \eqref{est Z2p alpha<beta}, \eqref{est Z5p+Z6p} into \eqref{evol Ep p>=2} and, for $\delta'\in(0,1)$, set $\delta$ so that $3\delta\,\mathcal{C}_p=1-\delta'$. 
Then, $\delta^{-1}\lesssim\mathcal{C}_p$ and it holds
\begin{align}\label{evol Ep p>=2 final}
&\|u_a\|^{q(p)}_{L^\infty(0,T;L^{q(p)}(\Omega))}+\|u_b\|^{r(p)}_{L^\infty(0,T;L^{r(p)}(\Omega))}
+\|\nabla u_a^{q(p)/2}\|_{L^2(\Omega_T)}^2+ \|\nabla u_b^{r(p)/2}\|_{L^2(\Omega_T)}^2\notag\\
&\qquad\qquad+\delta'\left(\|u_a\|_{L^{q(p)+1}(\Omega_T)}^{q(p)+1}+\|u_b\|_{L^{r(p)+1}(\Omega_T)}^{r(p)+1}\right)-\int_0^TI^{p}_{\text{fast}}(t)\,dt\notag\\
&\quad\le \mathcal{C}_p(\|u_a\|_{L^p(\Omega_T)}^p +\|u_b\|_{L^p(\Omega_T)}^p+\mathcal{E}_{p}(0)+|\Omega_T|)\\
&\qquad\qquad+ (1+\mathcal{C}_p^{\f{q(p)+1}2}+ \mathcal{C}_p^{\f{r(p)+1}2}+\mathcal{C}_p^{\f{(q(p)+1)(r(p)+1)}{2+3q(p)-r(p)}})\notag\,.
\end{align}
\subsection{Bootstrapping and end of the proof}
We are now ready to prove Lemma \ref{energy estimate}. 

First, recall that $q(p_\alpha)=r(p_\beta)=2$. So, if $\alpha<\beta$, estimates \eqref{est from E 1} and  \eqref{est from E 2} imply \eqref{basic energy estimate lemma 2.2} by interpolation. If $\alpha=\beta$, \eqref{est from E 3} is exactly \eqref{basic energy estimate lemma 2.2}.

Next, we have proved that, if $p\in[2,p^0_{\alpha,\beta})$ (see \eqref{def palphabeta}), \eqref{evol Ep p>=2 final} holds true and gives a bound on $\|u_a\|_{L^{q(p)}(\Omega_T)}$ and $\|u_b\|_{L^{r(p)}(\Omega_T)}$ as soon as we have a bound on the $L^{p}(\Omega_T)$ norm of $u_a$ and $u_b$. Hence, recalling that $q(p)\le r(p)$ and starting from exponent $2=q(p_\alpha)$ and the $L^{2}(\Omega_T)$ estimates of $u_a$ and $u_b$ in \eqref{result p=1}, we can bootstrap an $L^{q(p)}(\Omega_T)$ bound of $u_a$ and $u_b$ to an $L^{q(q(p))}(\Omega_T)$ bound of $u_a$ and $u_b$ until  $q(p)<p^0_{\alpha,\beta}$. The two cases below have to be considered.
\begin{itemize}
\item[(i)] $\beta-\alpha\in[0,2(\alpha+1)]$. Then, $p^0_{\alpha,\beta}=+\infty$ and estimate \eqref{evol Ep p>=2 final} implies \eqref{energy estimate 1 lemma 2.2}. 
\item[(ii)] $\beta-\alpha\in(2(\alpha+1),2(\alpha+3))=\cup_nI_n$, (see \eqref{def In}). Then, $p^0_{\alpha,\beta}\in(2,+\infty)$. In order to set up the bootstrap procedure, we denote
\[
p^n_{\alpha,\beta}:=1+\f4{(\alpha+1)^n(\beta-3\alpha-2)}\,,\qquad n\ge1\,,
\]
and we observe that, by \eqref{def q r}, it holds 
\[
p^0_{\alpha,\beta}=(\overbrace{q\circ\cdots\circ q}^{n-times})(p^n_{\alpha,\beta})=:q_n(p^n_{\alpha,\beta})\,,\qquad n\ge0\,.
\] 
Let $n_{\alpha,\beta}\ge0$ be the largest integer such that $\beta-\alpha\in I_{n_{\alpha,\beta}}$. We have that $2<p^{n_{\alpha,\beta}}_{\alpha,\beta}$ and $q_{n_{\alpha,\beta}}(2)<q_{n_{\alpha,\beta}}(p^{n_{\alpha,\beta}}_{\alpha,\beta})=p^0_{\alpha,\beta}$. Therefore, starting from \eqref{result p=1}, we can bootstrap \eqref{evol Ep p>=2 final} till $q_{n_{\alpha,\beta}}(2)=(\alpha+1)^{n_{\alpha,\beta}}+1$ to get estimate \eqref{energy estimate 2 lemma 2.2}. 
\end{itemize}

Finally, taking $p=2$ in \eqref{evol Ep p>=2 final}, by \eqref{Ifreap2} we obtain  \eqref{evol Ep p=2 final}. The maximal regularity estimate \eqref{lap v} and the gradient estimate \eqref{gradient v} follow by \eqref{max reg} and \eqref{est gradv} respectively, using \eqref{energy estimate 1 lemma 2.2} or  \eqref{energy estimate 2 lemma 2.2}, according to the value of $\beta-\alpha$. For that, it is worth recall that when $\beta-\alpha\in(2(\alpha+1),2(\alpha+3))$ and $p=(\alpha+1)^{n_{\alpha,\beta}}+1$, then $q(p)+1=2+(\alpha+1)^{n_{\alpha,\beta}+1}$ . 

\section{Existence for the cross-diffusion system}\label{sect proof existence macro}
This section is devoted to the proofs of Lemma \ref{converg (ua,ub)} and Theorem \ref{thm existence}. The latter follows by compactness arguments based on the previous estimates on the unique global strict solution $(u_a^\e,u_b^\e,v^\e)$ of \eqref{meso system}--\eqref{def phi psi}. The key point is the identification of the limit (as $\e\to0$) of the densities pair $(u^\e_a,u^\e_b)$ with the unique solution of the nonlinear system \eqref{nonlinearsys} corresponding to the limit of the densities pair $(u^\e\!=u_a^\e+u_b^\e,v^\e)$. This is the object of Lemma \ref{converg (ua,ub)}. 

\subsection{Proof of Lemma \ref{converg (ua,ub)}}
By estimates \eqref{evol Ep p=2 final}, \eqref{result p=1} there exists $C(T)>0$ such that, for all $\e>0$,
\[
\|(\Lambda^{1/2}Q)(u_a^\e,u_b^\e,v^\e)\|_{L^2(\Omega_T)}\le\sqrt{\e}\, C(T)\,.
\]
As $\Lambda(u_a^\e,u_b^\e,v^\e)\ge A^\alpha>0$, the latter implies 
\begin{equation}\label{conv lambda Q}
\|Q(u_a^\e,u_b^\e,v^\e)\|_{L^2(\Omega_T)}\le A^{-\alpha/2}\sqrt{\e}\, C(T)\,.
\end{equation}

Now, following \cite{Brocchieri2021}, let us define $\caQ(u,u_b,v)\!:=\!Q(u-u_b,u_b,v)$, for $(u,u_b,v)\!\in\R^3_+$ such that $u_b\le u$. Note that, by \eqref{def Q}, \eqref{def phi psi}, for all $(u_a,u_b,v)\!\in\R^3_+$,
\begin{align}
&\pa_1Q(u_a,u_b,v)=-\psi(\afast u_a+\cfast v)/\Lambda(u_a,u_b,v)\label{pa1Q}\\
&\qquad-\afast(u_a+u_b)\phi(\bfast u_b+\dfast v)\psi'(\afast u_a+\cfast v)/\Lambda^{2}(u_a,u_b,v)\notag\\
&\pa_2Q(u_a,u_b,v)=\phi(\bfast u_b+\dfast v)/\Lambda(u_a,u_b,v)\label{pa2Q}\\
&\qquad+\bfast(u_a+u_b)\phi'(\bfast u_b+\dfast v)\psi(\afast u_a+\cfast v)/\Lambda^{2}(u_a,u_b,v)\notag\,.
\end{align}
Hence, 
\beq\label{-pa1Q+pa2Q}
\pa_2\caQ(u,u_b,v)=-\pa_1Q(u-u_b,u_b,v)+\pa_2Q(u-u_b,u_b,v)\ge1\,.
\eeq
Recalling that $(u_a^*(u^\e,v^\e),u_b^*(u^\e,v^\e))$ is the unique solution of system \eqref{nonlinearsys} corresponding to $(u^\e,v^\e)$, it holds $\caQ(u^\e,u_b^*(u^\e,v^\e),v^\e)=0$. Therefore, for some intermediate value $\xi\!\in\![0,u^\e]$ between $u_b^\e\!\in\![0,u^\e]$ and $u_b^*(u^\e,v^\e)\!\in\![0,u^\e]$, it follows
\[
Q(u_a^\e,u_b^\e,v^\e)=\caQ(u^\e,u_b^\e,v^\e)-\caQ(u^\e,u_b^*(u^\e,v^\e),v^\e)=\pa_2\caQ(u^\e,\xi,v^\e)(u_b^\e-u_b^*(u^\e,v^\e))\,.
\]
Using \eqref{-pa1Q+pa2Q} and the previous identity, we have
\beq\label{ube-ube*}
|u_b^\e-u_b^*(u^\e,v^\e)|\le \pa_2\caQ(u^\e,\xi,v^\e)|u_b^\e-u_b^*(u^\e,v^\e)|=|Q(u_a^\e,u_b^\e,v^\e)|\,,
\eeq
and \eqref{conv slow manifold} follows by \eqref{conv lambda Q}. The Lemma is proved. 

\subsection{Proof of Theorem \ref{thm existence}}
In the sequel $p$ will satisfies~\eqref{def p thm existence macro}.

\bigskip
                                                                                                                                                                                                                                                                                                                                           \noindent{\bf Convergence of $(v_\e)_\e$.}  Thanks to estimates \eqref{bound v}, \eqref{lap v}, \eqref{gradient v}, for all $T>0$, the sequence $(v^\e)_\e$ is bounded in $W^{1,p}(\Omega_T)\cap L^p(0,T;W^{2,p}(\Omega))\cap L^{\infty}((0,\infty)\times\Omega)$. By \eqref{energy estimate 1 lemma 2.2},\eqref{energy estimate 2 lemma 2.2}, the sequence $(f_v(u_a^\e,u_b^\e,v^\e))_\e$ is bounded in $L^p(\Omega_T)$, (see also \eqref{def react functions}). Owing to the Aubin-Lions's Lemma \cite{Aubin1987} and standard weak compactness arguments, it follows that, for any $T>0$, there exists a subsequence of $(v^\e)_\e$ (still denoted $v^\e$) and $v\in W^{1,p}(\Omega_T)\cap L^p(0,T;W^{2,p}(\Omega))\cap L^{\infty}((0,T)\times\Omega)$ such that,
\beq\label{conv of v}
\begin{split}
&v^\e\to v\,,\quad\text{in } L^p(0,T;W^{1,p}(\Omega))\text{ and a.e. in } \Omega_T\,,\quad\text{as }\e\to0\,,\\
& v^\e \rightharpoonup  v\, \,\quad\text{ in } W^{1,p}(\Omega_T)\text{ and in }L^p(0,T;W^{2,p}(\Omega))\,, \quad\text{as }\e\to0\,.
\end{split}
\eeq
Moreover, $v$ is nonnegative as $v^\e$ is nonnegative. 
\medskip

\noindent{\bf Convergence of $(u^\e)_\e\!=\!(u_a^\e+u_b^\e)_\e$.}
By estimates \eqref{basic energy estimate lemma 2.2}--\eqref{energy estimate 2 lemma 2.2}, (see also Remark~\ref{remark grad u in L2}), for all $T>0$, the sequences $(u^\e_a)_\e,\,(u^\e_b)_\e$ are bounded in $L^2(0,T;H^1(\Omega))\cap L^{p}(\Omega_T)\cap L^\infty(0,T;L^{p-1}(\Omega))$. Hence, $(f_u(u^\e_a,u^\e_b,v^\e))_\e$ is bounded in $L^{\f p2}(\Omega_T)$ (because of the quadratic terms, see \eqref{def fu}, \eqref{def react functions}). We denote $s=\f p2\wedge 2$ and we continue according the value of $\beta-\alpha$. 

If $\beta-\alpha\in[0,2(\alpha+1)]$, we can choose $p=4$ (see \eqref{def p thm existence macro}), which yields $s=2$ and $(f_u(u^\e_a,u^\e_b,v^\e))_\e$ is bounded in $L^{2}(\Omega_T)$. Therefore, the equation satisfied by $u^\e$ (see \eqref{meso system ue ve}) implies that $(\pa_t u^\e)_{\e}$ is bounded in $L^{2}(0,T;H^{-1}(\Omega))$. 

If $\beta-\alpha\in(2(\alpha+1),2(\alpha+3))$, then $p>3$ (see \eqref{def p thm existence macro}), $s\in(\f32,2]$ and we have to argue according the space dimension. If $N\le 6$, by Sobolev's embedding theorem $L^s(\Omega)\subset H^{-1}(\Omega)$ with continuous embedding, so that $(\pa_t u^\e)_{\e}$ is bounded in $L^{s}(0,T;H^{-1}(\Omega))$. If $N>6$, by the assumption $(\alpha+1)^{n_{\alpha,\beta}+1}\ge2$, we have that $p\ge4$, which yields $s=2$ and again $(\pa_t u^\e)_{\e}$ is bounded in $L^{2}(0,T;H^{-1}(\Omega))$. 

As above, for any $T>0$, it holds the existence of a subsequence of $(u^\e)_\e$ (still denoted $u^\e$) and $u\in L^2(0,T;H^1(\Omega))\cap L^p(\Omega_T)\cap L^\infty(0,T;L^{p-1}(\Omega))$, such that
\beq\label{conv of u}
\begin{split}
& u^\e=u_a^\e+u_b^\e\to u\,,\quad\text{in } L^2(\Omega_T)\text{ and a.e. in } \Omega_T\,,\quad\text{as }\e\to0\,,\\
& u^\e\,\rightharpoonup u\,\qquad\text{ in } L^2(0,T;H^1(\Omega))\text{ and in } L^p(\Omega_T)\,,\quad\text{as }\e\to0\,,\\
& \pa_tu^\e\rightharpoonup\pa_t u\,\qquad\text{ in } L^s(0,T;H^{-1}(\Omega))\,,\quad\text{as }\e\to0\,. 
\end{split}
\eeq
Furthermore, $u$ is nonnegative as $u^\e$ is positive.
\medskip 

\noindent{\bf Convergence of $(u^\e_a)_\e,(u^\e_b)_\e$.} 
By the boundedness of the sequences $(u^\e_a)_\e,\,(u^\e_b)_\e$ quoted above, for any $T>0$, there exists $u_a,u_b\in L^2(0,T;H^1(\Omega))\cap L^{p}(\Omega_T)\cap L^\infty(0,T;L^{p-1}(\Omega))$ such that, for subsequences (still denoted $u^\e_a$ and $u^\e_b$)
\beq\label{weak conv ua ub}
u^\e_a\rightharpoonup u_a\,,\quad u^\e_b\rightharpoonup u_b \,\quad\text{ in } L^2(0,T;H^1(\Omega))\text{ and in } L^p(\Omega_T)\,,\quad\text{as }\e\to0\,.
\eeq
On the other hand, using \eqref{conv slow manifold}, there exists subsequences such that
\[
Q(u_a^\e,u_b^\e,v^\e)\to 0\,,\quad |u_b^\e-u_b^*(u^\e,v^\e)|=|u_a^\e-u_a^*(u^\e,v^\e)|\to 0\,,\quad \text{ a.e. in }\Omega_T\,.
\]
By the continuity of the map $(u_a^*,u_b^*)$ (see Lemma \ref{property cross-diff}) and the a.e. convergences of $(u^\e,v^\e)$ towards $(u,v)$ obtained above, it follows that 
\beq\label{conv solu NLS}
(u_a^\e,u_b^\e)\to(u_a^*(u,v),u_b^*(u,v))\,,\qquad\text{a.e. in }\Omega_T\,,\quad\text{as }\e\to0\,,
\eeq
Therefore, $(u_a^*(u,v),u_b^*(u,v))=(u_a,u_b)$ a.e. in $\Omega_T$, and 
\[
u_a^*(u,v),u_b^*(u,v)\in L^2(0,T;H^1(\Omega))\cap L^{p}(\Omega_T)\cap L^\infty(0,T;L^{p-1}(\Omega))\,.
\]

\noindent{\bf Diagonal extraction.} Since $T>0$ is arbitrarily large, we can apply the diagonal extraction argument. It follows that there exists a subsequence $(\e_k)_k$ and a pair of nonnegative measurable functions $(u,v):(0,\infty)\times\Omega\to\R^2_+$ satisfying $(i), (ii)$ in the statement of Theorem \ref{thm existence} and
\beq
u^{\e_k}\to u\,,\ v^{\e_k}\to v\,,\ (u_a^{\e_k},u_b^{\e_k})\to(u_a^*(u,v),u_b^*(u,v))\,,\ \text{ a.e. in } (0,\infty)\times\Omega\,,\text{ as }k\to\infty\,,
\eeq
and also the convergence in \eqref{conv of v}, \eqref{conv of u}, \eqref{weak conv ua ub}, thus giving $(iii), (iv)$. 
\medskip

\noindent{\bf Conclusion.} It remains to show that $(u,v)$ is a global (weak, strong) solution according to \eqref{weak eq u}. To begin with, we consider the weak formulation of the equation for $u^\e$ in \eqref{meso system ue ve} with test functions as in \eqref{weak eq u}, i.e.
\begin{align}\label{weak eq uae+ube}
-\int_0^T\!\!\int_\Omega u^\e\,\partial_{t}w\,dxdt
&-\int_{\Omega}u^{\init}\,w(0)\,dx+\int_0^T\!\!\int_\Omega\nabla(d_au^\e_a+d_bu^\e_b)\cdot\nabla w\,dxdt\notag\\
&=\int_0^T\!\!\int_\Omega f_u(u_a^\e, u_b^\e,v^\e)\,w\,dxdt\,.
\end{align}
It is worth noticing that the term in the right hand side of \eqref{weak eq uae+ube} is well defined. Indeed, this is clear if $\beta-\alpha\in[0,2(\alpha+1)]$ or if $\beta-\alpha\in(2(\alpha+1),2(\alpha+3))$ and $N>6$, since then $f_u(u_a^\e, u_b^\e,v^\e)\in L^2(\Omega_T)$, as observed above. On the other hand, if $\beta-\alpha\in(2(\alpha+1),2(\alpha+3))$ and $2<N\le6$, the Sobolev's embedding $H^1(\Omega)\subset L^{\f{2N}{N-2}}(\Omega)$ implies that the quadratic terms in $f_u(u_a^\e, u_b^\e,v^\e)$ belong to $L^1(0,T;L^{\f{N}{N-2}}(\Omega))$, while $w\in L^\infty(0,T;L^{\f{2N}{N-2}}(\Omega))$, and we are able to conclude. If $N=1,2$, similar arguments give us the claim. 

Hence, thanks to the above convergence properties of $u_a^{\e_k},u_b^{\e_k},v^{\e_k}$,  the boundedness of $(f_u(u_a^\e, u_b^\e,v^\e))_\e$ in $L^{\f p2}(\Omega_T)$, for all $T>0$, the convergence of $(f_u(u_a^{\e_k},u_b^{\e_k},v^{\e_k}))_k$ towards $f_u(u_a^*(u,v),u_b^*(u,v),v)$ a.e. in $(0,\infty)\times\Omega$ as $k\to\infty$, we can pass to the limit $k\to\infty$ in \eqref{weak eq uae+ube} and the equation for $u$ in \eqref{weak eq u} holds true. 

Furthermore, $u\in W^{1,s}(0,T;H^{-1}(\Omega))$ and $W^{1,s}(0,T;H^{-1}(\Omega))$ is continuously embedded in $C^0([0,T];H^{-1}(\Omega))$. Hence, the operator $w\to w(0)$ is weakly sequentially continuous from $W^{1,s}(0,T;H^{-1}(\Omega))$ weak to $H^{-1}(\Omega)$ weak. As, for all $T\!>\!0$, $u^{\e_k}\rightharpoonup u$ in $W^{1,s}(0,T;H^{-1}(\Omega))$  and $u^{\e_k}(0)\!=\!u^\init$, we have $u(0)\!=\!u^\init$ a.e. in $\Omega$. 

Next, concerning the $v$ component of the solution, recall that $(f_v(u_a^\e,u_b^\e,v^\e))_\e$ is bounded in $L^p(\Omega_T)$, for all $T>0$, and that $f_v(u_a^{\e_k},u_b^{\e_k},v^{\e_k})\to f_v(u_a^*(u,v),u_b^*(u,v),v)$ a.e. in $(0,\infty)\times\Omega$ as $k\to\infty$. Hence, the previous convergence properties of $(v^{\e_k})_k$, applied to the equation
\[
\partial_tv^{\e_k}-d_v\Delta v^{\e_k}=f_v(u_a^{\e_k},u_b^{\e_k},\,v^{\e_k})
\] 
yield that
\[
\pa_tv=d_v\Delta v+f_v(u_a^*(u,v),u_b^*(u,v),v)\,,\quad\text{in }{\cal D}'((0,\infty)\times\Omega)\,,
\]
and then that $v$ satisfies the above equation in $L^p(\Omega_T)$, for all $T>0$, by the $L^p(\Omega_T)$-integrability of each term in the equation. As $W^{1,p}(\Omega_T)\subset W^{1,p}(0,T;L^p(\Omega))\subset C^0([0,T];L^p(\Omega))$ with continuous embeddings, the operator $w\to w(0)$ is weakly sequentially continuous from $W^{1,p}(\Omega_T)$ weak to $L^p(\Omega)$ weak. As, for all $T>0$, $v^{\e_k}\rightharpoonup v$ in $W^{1,p}(\Omega_T)$ and $v^{\e_k}(0)=v^\init$, we have $v(0)=v^\init$ a.e. in $\Omega$. Finally, the regularity of $v$ implies that $v$ satisfies homogeneous Neumann boundary condition in the sense of traces. Theorem \ref{thm existence} is proved. 
\section{Uniqueness for the cross-diffusion system}\label{sect proof uniqueness macro}
This section is devoted to the proof of Theorem \ref{Thm uniqueness}. More precisely, we show that, for all $T>0$, the function 
\beq\label{def philambda}
\phi_\lambda(\tau):=\|u_1-u_2\|^2_{L^2(\Omega_\tau)}+\lambda\|v_1-v_2\|^2_{L^2(\Omega_\tau)}\,,\qquad \tau\in[0,T]\,,
\eeq
with $\lambda>0$ large enough, satisfies an integral inequality that gives \eqref{est stability} by Gronwall's lemma. As a consequence, solutions of \eqref{weak eq u}, whose components are both bounded, are unique. The key tool are test functions introduced by Oleinik (see \cite{Vazquez1992} and the references therein, see also \cite{LionsMagenes1972}). 

Let us denote
\[
f_u^*(u,v):=f_u(u_a^*(u,v),u_b^*(u,v),v)\,,\qquad f_v^*(u,v):=f_v(u_a^*(u,v),u_b^*(u,v),v)\,.
\]
The weak formulation of \eqref{macro system}, \eqref{structure}, \eqref{nonlinearsys} for both components $u,v$, writes as
\begin{align}
-\iint_{\Omega_T} u\,\pa_tw_1\,dxdt&-\int_\Omega u^\init\,w_1(0)\,dx\notag\\
&+\iint_{\Omega_T}\nabla A(u,v)\cdot\nabla w_1\,dxdt=\iint_{\Omega_T}f_u^*(u,v)\,w_1\,dxdt\,,\label{weak formulation u}\\
-\iint_{\Omega_T} v\,\pa_tw_2\,dxdt&-\int_\Omega v^\init\,w_2(0)\,dx\notag\\
&+\iint_{\Omega_T}\nabla v\cdot\nabla w_2\,dxdt=\iint_{\Omega_T}f_v^*(u,v)\,w_2\,dxdt\label{weak formulation v}\,,
\end{align}
for all $T>0$ and all $w_1,w_2\in C^1([0,T];H^1(\Omega))$ such that $w(T)=0$. 

We may assume without loss of generality that $d_b>d_a$, so that the function $A(u,v)$ in \eqref{structure}  rewrites as 
\beq\label{rewrite A}
A(u,v)=d_a u_a^*(u,v)+d_b u_b^*(u,v)=d_a\,u+(d_b-d_a)\,u_b^*(u,v)\,,
\eeq
and, using \eqref{partial u ua* ub*}, \eqref{partial v ua* ub*}, it holds
\begin{align}
&\pa_1A(u,v)=d_a+(d_b-d_a)\,\pa_1u_b^*(u,v)\in (d_a,d_b)\label{pa1A}\\
&\pa_2A(u,v)=(d_b-d_a)\,\pa_2u_b^*(u,v)\in \left(-(d_b-d_a)\f{\dfast}{\bfast},(d_b-d_a)\f{\cfast}{\afast}\right)\,. \label{pa2A}
\end{align}

Let $(u_i,v_i)$, $i=1,2$, be two solutions with initial data $(u_i^\init,v_i^\init)$ and let us denote $A_i:=A(u_i,v_i)$, $i=1,2$. As the functions $A_i$ belong to $L^2(0,T;H^1(\Omega))$, the function
\beq\label{test function 1}
w_1^\tau(t,x)=
\begin{cases}
\int_t^\tau(A_1(s,x)-A_2(s,x))\,ds\,,&\text{if } 0\le t\le\tau\,,\\
\ 0\,,&\text{if } \tau\le t\le T\,,
\end{cases}
\eeq
belongs to $H^1([0,T];H^1(\Omega))$ and $w_1^\tau(T)=0$, for all $\tau\in[0,T]$. Therefore, by the density of $C^1([0,T];H^1(\Omega))$ in $H^1([0,T];H^1(\Omega))$, we can use $w_1^\tau$ as test functions. Hence, testing the equation satisfied by $u_1-u_2$ against $w_1^\tau$ , we obtain (see \eqref{weak formulation u})
\begin{align}
\iint_{\Omega_\tau}&(u_1-u_2)(t,x)(A_1-A_2)(t,x)\,dxdt-\int_\Omega(u_1^\init-u_2^\init)(x)\int_0^\tau(A_1-A_2)(s,x)\,ds\,dx\notag\\
&\qquad+\iint_{\Omega_\tau}\nabla(A_1-A_2)(t,x)\cdot\int_t^\tau\nabla(A_1-A_2)(s,x)\,ds\,dxdt\label{weak form u1-u2}\\
&= \iint_{\Omega_\tau}(f^*_u(u_1,v_1)-f^*_u(u_2,v_2))(t,x)\int_t^\tau(A_1-A_2)(s,x)\,ds\,dxdt\notag\,.
\end{align}

It is convenient to split the first term in the left hand side of \eqref{weak form u1-u2} as 
\begin{align}
\iint_{\Omega_\tau}(u_1-u_2)(t,x)(A_1&-A_2)(t,x)\,dxdt\notag\\
&=\iint_{\Omega_\tau}(u_1-u_2)(t,x)(A(u_1,v_1)-A(u_2,v_1))(t,x)\,dxdt\notag\\
&\ +\iint_{\Omega_\tau}(u_1-u_2)(t,x)(A(u_2,v_1)-A(u_2,v_2))(t,x)\,dxdt\notag\\
&=: I_1+I_2\label{first term}\,.
\end{align}
Indeed, by \eqref{pa1A}, $A$ is increasing in $u$ with $\pa_1A(u,v)$ lower bounded, so that $I_1$ is positive and lower bounded
\beq\label{first term 1}
I_1\ge d_a\|u_1-u_2\|^2_{L^2(\Omega_\tau)}\,.
\eeq
On the other hand, denoting $c_2:=\!(d_b-d_a)\!\left(\f{\cfast}{\afast}\vee\f{\dfast}{\bfast}\right)$, by \eqref{pa2A} we obtain for $I_2$
\beq\label{first term 2}
|I_2|\le \sigma\, d_a\|u_1-u_2\|^2_{L^2(\Omega_\tau)}+\f{c_2^2}{4\sigma\,d_a}\|v_1-v_2\|^2_{L^2(\Omega_\tau)}\,,
\eeq
with $\sigma>0$ to be chosen later.  

Using \eqref{pa1A}, \eqref{pa2A} again, the second term in the left hand side of \eqref{weak form u1-u2} can be estimated as following, for all $\tau\in[0,T]$, 
\begin{align}
&\int_\Omega(u_1^\init-u_2^\init)(x)\int_0^\tau(A_1-A_2)(s,x)\,ds\,dx\le 
\f{\sigma\, d_a}{2d_b^2}\|A_1-A_2\|^2_{L^2(\Omega_\tau)}\notag\\
&\qquad\qquad\qquad+\f {Td_b^2}{2\sigma\, d_a}\|u_1^\init-u_2^\init\|^2_{L^2(\Omega)}\notag\\
&\le \sigma\, d_a\|u_1-u_2\|^2_{L^2(\Omega_\tau)}+\sigma\f{c_2^2d_a}{d_b^2}\|v_1-v_2\|^2_{L^2(\Omega_\tau)}
+\f {Td_b^2}{2\sigma\, d_a}\|u_1^\init-u_2^\init\|^2_{L^2(\Omega)}\,.\label{second term}
\end{align}

Denoting $y_i(t,x)=\int_t^\tau\pa_{x_i}(A_1-A_2)(s,x)\,ds$, the third term in the left hand side of \eqref{weak form u1-u2} turn out to be positive and it can be neglected, since it writes as
\beq\label{positive term}
\begin{split}
\iint_{\Omega_\tau}\nabla(A_1-A_2)(t,x)&\cdot\int_t^\tau\nabla(A_1-A_2)(s,x)\,ds\,dxdt\\
&=-\f12\int_\Omega\int_0^\tau\pa_t\left(\sum_iy_i^2(t,x)\right)\,dtdx\\
&=\f12\int_\Omega\sum_i\left(\int_0^\tau\pa_{x_i}(A_1-A_2)(s,x)\,ds\right)^2\,dx\,.
\end{split}
\eeq

Finally, for the term in the right hand side of \eqref{weak form u1-u2} we have
\beq\label{react term for u 1}
\begin{split}
\iint_{\Omega_\tau}&(f^*_u(u_1,v_1)-f^*_u(u_2,v_2))(t,x)\int_t^\tau(A_1-A_2)(s,x)\,ds\,dxdt\\
&\le\int_0^\tau ds\int_0^s dt\int_\Omega dx|(f^*_u(u_1,v_1)-f^*_u(u_2,v_2))(t,x)||(A_1-A_2)(s,x)|
\end{split}
\eeq
and 
\begin{align}
&\int_\Omega |(f^*_u(u_1,v_1)-f^*_u(u_2,v_2))(t,x)||(A_1-A_2)(s,x)|dx\notag\\
&\le \int_\Omega |(f^*_u(u_1,v_1)-f^*_u(u_2,v_2))(t,x)|\left(d_b|(u_1-u_2)(s,x)|+c_2|(v_1-v_2)(s,x)|\right)dx\notag\\
&\le \f{\sigma\,d_b^2}{T}\|(u_1-u_2)(s)\|^2_{L^2(\Omega)}+\f{\sigma\,c_2^2}{T}\|(v_1-v_2)(s)\|^2_{L^2(\Omega)}\notag\\
&\qquad\qquad +\f{T}{2\sigma}\|(f^*_u(u_1,v_1)-f^*_u(u_2,v_2))(t)\|^2_{L^2(\Omega)}\label{react term for u 2}\,.
\end{align}
Plugging \eqref{react term for u 2} into \eqref{react term for u 1}, we end up with 
\beq\label{react term for u 3}
\begin{split}
\iint_{\Omega_\tau}&(f^*_u(u_1,v_1)-f^*_u(u_2,v_2))(t,x)\int_t^\tau(A_1-A_2)(s,x)\,ds\,dxdt\\
&\le \sigma\, d_b^2\|u_1-u_2\|^2_{L^2(\Omega_\tau)}+\sigma\,c_2^2\|v_1-v_2\|^2_{L^2(\Omega_\tau)}\\
&\qquad\qquad +\f{T}{2\sigma}\int_0^\tau\int_0^s\|(f^*_u(u_1,v_1)-f^*_u(u_2,v_2))(t)\|^2_{L^2(\Omega)}\,dtds\,.
\end{split}
\eeq
Now, gathering \eqref{first term}--\eqref{positive term} and \eqref{react term for u 3}, we have that there exists $C_1\!=\!C_1(\sigma,d_a,d_b,c_2)\!>\!0$ such that, for all $T>0$, it holds
\begin{align}
(d_a&-\sigma(2d_a+d_b^2))\|u_1-u_2\|^2_{L^2(\Omega_\tau)}\le C_1\,\|v_1-v_2\|^2_{L^2(\Omega_\tau)}
+\f {Td_b^2}{2\sigma\, d_a}\|u_1^\init-u_2^\init\|^2_{L^2(\Omega)}\notag\\
&+\f{T}{2\sigma}\int_0^\tau\int_0^s\|(f^*_u(u_1,v_1)-f^*_u(u_2,v_2))(t)\|^2_{L^2(\Omega)}\,dtds\label{u1-u2}\,.
\end{align}

The same type of computations can be performed for the equation satisfied by $v_1-v_2$ using the test functions 
\beq\label{test function 2}
w_2^\tau(t,x)=
\begin{cases}
\int_t^\tau(v_1(s,x)-v_2(s,x))\,ds\,,&\text{if } 0\le t\le\tau\,,\\
\ 0\,,&\text{if } \tau\le t\le T\,.
\end{cases}
\eeq
 Indeed, from \eqref{weak formulation v} we have
\begin{align}
\iint_{\Omega_\tau}(v_1-&v_2)^2(t,x)\,dxdt-\int_\Omega(v_1^\init-v_2^\init)(x)\int_0^\tau(v_1-v_2)(s,x)\,ds\,dx\notag\\
&+\iint_{\Omega_\tau}\nabla(v_1-v_2)(t,x)\cdot\int_t^\tau\nabla(v_1-v_2)(s,x)\,ds\,dxdt\label{weak form v1-v2}\\
&=\iint_{\Omega_\tau}(f^*_v(u_1,v_1)-f^*_v(u_2,v_2))(t,x)\int_t^\tau(v_1-v_2)(s,x)\,ds\,dxdt\notag\,.
\end{align}
The first term in the left hand side of \eqref{weak form v1-v2} is left as it is, while all the other terms are estimated similarly as before, to obtain 
\beq\label{v1-v2}
\begin{split}
(1-2\sigma)\|v_1&-v_2\|^2_{L^2(\Omega_\tau)}\le\f T{4\sigma}\|v_1^\init-v_2^\init\|^2_{L^2(\Omega)}\\
&+\f T{4\sigma}\int_0^\tau\int_0^s\|(f^*_v(u_1,v_1)-f^*_v(u_2,v_2))(t)\|^2_{L^2(\Omega)}\,dtds\,.
\end{split}
\eeq

It remains to choose $\sigma>0$ small enough to have $d_a-\sigma(2d_a+d_b^2)>0$ in \eqref{u1-u2} and then $\lambda>0$ large enough (depending on $C_1(\sigma,d_a,d_b,c_2)$), so that, adding \eqref{u1-u2}, \eqref{v1-v2}, the function $\phi_\lambda(\tau)$ in \eqref{def philambda} satisfies, for all $\tau\in[0,T]$, 
\beq\label{int eq philambda}
\begin{split}
\phi_\lambda(\tau)\le&\ T\,C_2\left(\|u_1^\init-u_2^\init\|^2_{L^2(\Omega)}+\lambda\, \|v_1^\init-v_2^\init\|^2_{L^2(\Omega)}\right)\\
&+T\,C_3\left(\int_0^\tau\int_0^s\|(f^*_u(u_1,v_1)-f^*_u(u_2,v_2))(t)\|^2_{L^2(\Omega)}\,dtds\right.\\
&\qquad\qquad+\left.\lambda\int_0^\tau\int_0^s\|(f^*_v(u_1,v_1)-f^*_v(u_2,v_2))(t)\|^2_{L^2(\Omega)}\,dtds\right)\,,
\end{split}
\eeq
where $C_2,C_3$ are positive constants depending only on $\sigma,d_a,d_b,c_2$. 

Finally, as $f_u(u_a,u_b,v), f_v(u_a,u_b,v)$ in \eqref{def fu}, \eqref{def react functions} are locally Lipschitz continuous (due to the quadratic terms) and recalling that $u_a^*,u_b^*$  are $C^1(\R_+^2)$ with bounded gradient, $f^*_u(u,v), f^*_v(u,v)$ are also locally Lipschitz continuous. Therefore, there exists a positive constant $C_4$, depending on $\sigma,d_a,d_b,c_2,\lambda$ and 
$\|u_i\|_{L^\infty((0,T)\times\Omega)}$, $\|v_i\|_{L^\infty((0,\infty)\times\Omega)}$, such that \eqref{int eq philambda} gives us, for all $\tau\in[0,T]$,  
\[
\phi_\lambda(\tau)\le\ T\,C_2\left(\|u_1^\init-u_2^\init\|^2_{L^2(\Omega)}+\lambda\, \|v_1^\init-v_2^\init\|^2_{L^2(\Omega)}\right)+T\,C_4\int_0^\tau \phi_\lambda(s)\,ds\,.
\]
Gronwall's lemma implies, for all $\tau\in[0,T]$,
\[
\|u_1-u_2\|^2_{L^2(\Omega_\tau)}+\lambda\|v_1-v_2\|^2_{L^2(\Omega_\tau)}\le T\,C_2\,e^{T^2C_4}\left(\|u_1^\init-u_2^\init\|^2_{L^2(\Omega)}+\lambda\, \|v_1^\init-v_2^\init\|^2_{L^2(\Omega)}\right)
\]
and \eqref{est stability} follows. 
\section{The rate of convergence: proof of Theorem \ref{thm rate conv}}\label{sect proof rate of conv}
Let  $(u,v)$ be a nonnegative global classical solution of \eqref{macro system}--\eqref{nonlinearsys} satisfying \eqref{bound assumpt2}, $(u_a^\e,u_b^\e,v^\e)$ be the unique nonnegative global strict solution of \eqref{meso system}--\eqref{def phi psi} and $u^\e=u_a^\e+u_b^\e$. We denote
\beq\label{def Ue,Ve,We}
U^\e:=u^\e-u\,,\qquad V^\e:=v^\e-v\,,\qquad w:=u_b^*(u,v)\,,\qquad W^\e:=u_b^\e-w\,,
\eeq
and $\delta:=d_b-d_a>0$ so that, by \eqref{def Ue,Ve,We}, \eqref{structure},
\beq\label{A for rate of conv}
\begin{split}
d_a\,u_a^\e+d_b\,u_b^\e&=d_a\,u^\e+\delta\,u_b^\e=d_a(U^\e+u)+\delta(W^\e+w)\,,\\
A(u,v)&=d_a u_a^*(u,v)+d_b u_b^*(u,v)=d_a\,u+\delta\,w\,.
\end{split}
\eeq
It is worth noticing that by \eqref{A for rate of conv} and the definition of $w$ in \eqref{def Ue,Ve,We}, it holds
\[
\nabla A(u,v)= (d_a+\delta\,\pa_1u_b^*(u,v))\nabla u+ \delta\,\pa_2u_b^*(u,v)\nabla v\,.
\]
As observed in \eqref{pa1A}, $d_a+\delta\,\pa_1u_b^*(u,v)\!\in\!(d_a,d_b)$. Therefore, the homogeneous Neumann boundary conditions satisfied by $A(u,v)$ and $v$ give homogeneous Neumann boundary conditions for $u$ and consequently for $w$, $U^\e$ and $W^\e$. 

Furthermore, under the constraint $u=u_a+u_b$, with $u_a,u_b\in\R_+$, the reaction functions $f_u=f_a+f_b,f_b,f_v$ and $Q$, defined in \eqref{def react functions}, \eqref{def Q}, can be written as
\beq\label{new reactions}
\begin{split}
&f_u(u_a,u_b,v)=f_u(u-u_b,u_b,v)=:\caFu(u,u_b,v)\\
&f_v(u_a,u_b,v)=f_v(u-u_b,u_b,v)=:\caFv(u,u_b,v)\\
&f_b(u_a,u_b,v)=f_b(u-u_b,u_b,v)=:\caFw(u,u_b,v)\\
&Q(u_a,u_b,v)=Q(u-u_b,u_b,v)=:\caQ(u,u_b,v)\,.
\end{split}
\eeq
Then, it is easily seen that the triplet $(U^\e,V^\e,W^\e)$ satisfies over $(0,T)\times\Omega$ and for all $T>0$, the fast reaction-diffusion system below
\beq\label{U,V,W system}
\left\{
\begin{split}
&\partial_tU^\e-\Delta(d_aU^\e+\delta\,W^\e)=\caFu(U^\e+u,W^\e+w,V^\e+v)-\caFu(u,w,v)\\[1.4ex]
&\partial_tV^\e-d_v\,\Delta V^\e=\caFv(U^\e+u,W^\e+w,V^\e+v)-\caFv(u,w,v)\\[1.4ex]
&\partial_tW^\e-d_b\,\Delta W^\e=\caFw(U^\e+u,W^\e+w,V^\e+v)\\ 
&\qquad\qquad\qquad\qquad-\e^{-1}\caQ(U^\e+u,W^\e+w,V^\e+v)-(\partial_tw-d_b\,\Delta\,w)
\end{split}
\right.
\eeq
together with the initial and boundary conditions
\beq\label{BC+IC}
\begin{cases}
\nabla U^\e\cdot \vec{n}=\nabla V^\e\cdot \vec{n}=\nabla W^\e\cdot \vec{n}=0\,,\qquad&\In\,(0,T)\times\pa\Omega,\\
U^\e(0)=V^\e(0)=0\,,\quad W^\e(0)=u_b^\init-u_b^*(u^\init,v^\init),&\In\,\Omega\,.
\end{cases}
\eeq

Theorem \ref{thm rate conv} will be proved estimating (term by term) the time evolution of the functional
\beq
\begin{split}\label{def L}
&\caL(t):=\f{\gamma_1}2\| U^\e(t)\|^2_{L^2(\Omega)}+\f{\gamma_2}2\| V^\e(t)\|^2_{L^2(\Omega)}+\f{\e\gamma_3}2\| W^\e(t)\|^2_{L^2(\Omega)}\\
&\qquad\qquad\qquad\qquad\qquad\qquad+\f{\e\delta}2\|\nabla W^\e(t)\|^2_{L^2(\Omega)}+E(t)
\end{split}
\eeq
where
\begin{align}
&E(t):=\!-\!\int_\Omega \big[\caP(U^\e+u,w,V^\e+v)-\caP(u,w,v)\label{def E}\\
&\quad\qquad\qquad\qquad\qquad\qquad\qquad-\pa_1\caP(u,w,v)\,U^\e-\pa_3\caP(u,w,v)\,V^\e\big]dx\notag\\
&\caP(x,y,z):=\int_0^x\caQ(\xi,y,z)\,d\xi\,,\qquad (x,y,z)\in\R_+^3,\label{def caP}
\end{align}
and $\gamma_1,\gamma_2,\gamma_3$ are strictly positive constants (independent of $\e$) to be chosen later.  The functional \eqref{def L}--\eqref{def caP} is inspired by\! \cite{Iida2006}. However we require minor properties for $\caQ$, than in \cite{Iida2006}. The key tool employed to handle the functional is simply Taylor's~formula. 
\medskip

\noindent{\bf Step 1. Preliminaries.} 
To begin with, note that by assumptions \eqref{bound assumpt1}--\eqref{bound assumpt2} and estimate \eqref{bound v}, for all $T>0$, there exists ${\cal M}_T>0$, such that $u,v,u^\e,v^\e\in[0,{\cal M}_T]$, for all $\e\in(0,1)$. As $0\le u_b^\e\le u^\e$ and $0\le w=u_b^*(u,v)\le u$, it follows that $u_b^\e,w\in[0,{\cal M}_T]$, for all $\e\in(0,1)$. Therefore, in the sequel we can invoque the boundedness over $[0,{\cal M}_T]^3$ of the reaction functions in \eqref{new reactions} and their derivatives and of $\pa_iQ, \pa_{ij}Q,\pa_{ijk}Q$, $i,j,k\in\{1,2,3\}$. Indeed, as we have assumed $A,B\!>\!0$, for all $\alpha,\beta$, the function $Q$ belongs to $C^\infty(\R^3_+)$.  In particular, from \eqref{pa1Q} it is easily seen that $\pa_1Q\!<\!0$ and that there exists ${\cal C}(\alpha,\beta,A,B,{\cal M}_T)\!>\!0$ such that
\beq\label{inf pa1Q}
\caK_0:=\inf_{[0,{\cal M}_T]^3}|\pa_1Q(u_a,u_b,v)|\ge\inf_{[0,{\cal M}_T]^3}\psi(\afast u_a+\cfast v)/\Lambda(u_a,u_b,v)\ge {\cal C}>0\,.
\eeq
On the other hand, using the inequalities
\[
\psi/\Lambda\le1,\quad\phi/\Lambda\le1,\quad\afast u_a\psi'/\Lambda\le\alpha\psi/\Lambda\le\alpha,\quad
\psi'/\Lambda\le\psi'/\psi\le\alpha A^{-1},
\]
it holds
\beq\label{sup pa1Q}
\caK_1:=\sup_{[0,{\cal M}_T]^3}|\pa_1Q(u_a,u_b,v)|\le 1+\alpha+\afast {\cal M}_T\,\alpha A^{-1}\,.
\eeq
Furthermore, by \eqref{pa1Q}, \eqref{pa2Q}, we have
\beq\label{inf -pa1Q+pa2Q}
\inf_{\R_+^3}(-\pa_1Q(u_a,u_b,v)+\pa_2Q(u_a,u_b,v))\ge1\,.
\eeq

By \eqref{bound assumpt2} again, for all $T>0$, $\pa_tu,\pa_tv$, $\nabla u,\nabla v$, $\Delta u,\Delta v$, $\nabla\pa_tu, \nabla\pa_tv$, $\nabla\Delta u, \nabla\Delta v$ are bounded over $[0,T]\!\times\!\overline\Omega$. As a consequence, $\pa_tw,\Delta w$, $\nabla\pa_tw,\nabla\Delta w$ are also bounded over $[0,T]\!\times\!\overline\Omega$. Indeed, $w\!=\!u_b^*(u,v)$ and the gradient of the map $u_b^*:\R^2_+\mapsto\R^2_+$ is given by \eqref{grad ub*}, \eqref{grad q}, with $q$ defined in \eqref{def Q}. As $A,B>0$, $q\in C^\infty(\R^3_+)$. In particular $\pa_2q-\pa_1q\ge A^\alpha+B^\beta>0$. Therefore, $\pa_iu_b^*, \partial_{ij}u_b^*, \partial_{ijk}u_b^*$, $i,j,k\in\{1,2\}$, are locally bounded and this is sufficient since $u,v\in[0,{\cal M}_T]$. 

Hereafter, the $\xi_i$, $i=1,2,3$, appearing when applying Taylor's formula, belong all to $[0,{\cal M}_T]$. The constants in the estimates will change from line to line and only the dependence on $d_a,d_b,d_v,\delta,\caK_0,\caK_1$ is kept explicit. 
\medskip

\noindent{\bf Step 2. The evolution equation of $\| U^\e\|^2_{L^2(\Omega)}$ and $\| V^\e\|^2_{L^2(\Omega)}$.} We have
\[
\begin{split}
\f12\f d{dt}\| U^\e\|^2_{L^2(\Omega)}=&-d_a\| \nabla U^\e\|^2_{L^2(\Omega)}-\delta\int_\Omega \nabla W^\e\cdot\nabla U^\e\, dx\\
&+\int_\Omega [\caFu(U^\e+u,W^\e+w,V^\e+v)-\caFu(u,w,v)]\,U^\e\,dx\\
=&-d_a\| \nabla U^\e\|^2_{L^2(\Omega)}-\delta\int_\Omega \nabla W^\e\cdot\nabla U^\e\, dx\\
&+\int_\Omega \big(D\caFu(\xi_1,\xi_2,\xi_3)\cdot (U^\e,W^\e,V^\e)\big)U^\e\,dx\,,
\end{split}
\] 
and
\[
\begin{split}
\f12\f d{dt}\| V^\e\|^2_{L^2(\Omega)}=&-d_v\| \nabla V^\e\|^2_{L^2(\Omega)}\\
&+\int_\Omega [\caFv(U^\e+u,W^\e+w,V^\e+v)-\caFv(u,w,v)]\,V^\e\,dx\\
=&-d_v\| \nabla V^\e\|^2_{L^2(\Omega)}+\int_\Omega \big(D\caFv(\xi_1,\xi_2,\xi_3)\cdot (U^\e,W^\e,V^\e)\big)V^\e\,dx\,.
\end{split}
\]
Hence, the local boundedness of $D \caFu, D \caFv$ and ad hoc Young's inequality, give

\beq\label{L2 norm U}
\begin{split}
\f12\f d{dt}\| U^\e\|^2_{L^2(\Omega)}\le& -\f{d_a}2\| \nabla U^\e\|^2_{L^2(\Omega)}+\f{\delta^2}{2d_a}\|\nabla W^\e\|^2_{L^2(\Omega)}\\
&+C(\|U^\e\|_{L^2(\Omega)}+\|W^\e\|_{L^2(\Omega)}+\|V^\e\|_{L^2(\Omega)})\|U^\e\|_{L^2(\Omega)}\\
\le&-\f{d_a}2\| \nabla U^\e\|^2_{L^2(\Omega)}+\f{\delta^2}{2d_a}\|\nabla W^\e\|^2_{L^2(\Omega)}\\
&+C(\|U^\e\|^2_{L^2(\Omega)}+\|V^\e\|^2_{L^2(\Omega)})+\f1{4}\| W^\e\|^2_{L^2(\Omega)}\,,
\end{split}
\eeq
and similarly
\beq\label{L2 norm V}
\begin{split}
\f12\f d{dt}\| V^\e\|^2_{L^2(\Omega)}
\le& -d_v\| \nabla V^\e\|^2_{L^2(\Omega)}+C(\|U^\e\|^2_{L^2(\Omega)}+\|V^\e\|^2_{L^2(\Omega)})+\f1{4}\| W^\e\|^2_{L^2(\Omega)}\,.
\end{split}
\eeq

\noindent{\bf Step 3. The evolution equation of $\| W^\e\|^2_{L^2(\Omega)}$.}  We have
\begin{align}
&\f12\f d{dt}\| W^\e\|^2_{L^2(\Omega)}=-d_b\| \nabla W^\e\|^2_{L^2(\Omega)}+I_1+I_2\,,\label{est W}\\
&I_1:=\int_\Omega \caFw(U^\e+u,W^\e+w,V^\e+v)\,W^\e\,dx-\int_\Omega (\partial_tw-d_b\,\Delta\,w)\,W^\e\,dx\,,\notag\\
&I_2:=-\f1\e\int_\Omega\caQ(U^\e+u,W^\e+w,V^\e+v)\,W^\e\,dx\notag\,.
\end{align}
By the local boundedness of $\caFw$ and the boundedness of $\partial_tw,\Delta\,w$
\beq\label{est I1}
I_1\le C(1+d_b)\|W^\e\|_{L^1(\Omega)}\,.
\eeq
Next, observing that $\caQ(u,w,v)=Q(u_a^*(u,v),u_b^*(u,v),v)=0$ because $(u,v)$ is a solution of \eqref{macro system}--\eqref{nonlinearsys}, we write $I_2$ as
\[
\begin{split}
I_2=&-\f1\e\int_\Omega[\caQ(U^\e+u,W^\e+w,V^\e+v)-\caQ(U^\e+u,w,V^\e+v)]\,W^\e\,dx\\
&-\f1\e\int_\Omega[\caQ(U^\e+u,w,V^\e+v)-\caQ(u,w,v)]\,W^\e\,dx\,.
\end{split}
\]
The latter allows us to use again Taylor's formula to obtain
\[
\begin{split}
I_2=&-\f1\e\int_\Omega\pa_2\caQ(U^\e+u,\xi_1,V^\e+v)(W^\e)^2dx\\
&-\f1\e\int_\Omega [\pa_1\caQ(\xi_2,w,\xi_3)U^\e+\pa_3\caQ(\xi_2,w,\xi_3)V^\e]W^\e dx\,.
\end{split}
\]
As by the definition of $\caQ$ in \eqref{new reactions} it holds that $\pa_2\caQ=-\pa_1Q+\pa_2Q$, using \eqref{inf -pa1Q+pa2Q} together with the local boundedness of $D Q$, we get
\beq\label{est I2}
I_2\le -\f1\e\,\|W^\e\|^2_{L^2(\Omega)}+\f C\e(\|U^\e\|_{L^2(\Omega)}+\|V^\e\|_{L^2(\Omega)})\|W^\e\|_{L^2(\Omega)}\,.
\eeq
Plugging \eqref{est I1},\eqref{est I2} into \eqref{est W}, we end up with the estimate 
\[
\begin{split}
\f12\f d{dt}\| W^\e\|^2_{L^2(\Omega)}\le& -d_b\| \nabla W^\e\|^2_{L^2(\Omega)}+C(1+d_b)\|W^\e\|_{L^1(\Omega)}\\
&-\f1\e\,\|W^\e\|^2_{L^2(\Omega)}+\f C\e(\|U^\e\|_{L^2(\Omega)}+\|V^\e\|_{L^2(\Omega)})\|W^\e\|_{L^2(\Omega)}\,.
\end{split}
\]
Finally, we multiply the above inequality by $\e$ and use Young's inequality to~get 
\beq\label{est1 W}
\begin{split}
&\f\e2\f d{dt}\| W^\e\|^2_{L^2(\Omega)}\le -d_b\e\| \nabla W^\e\|^2_{L^2(\Omega)}+C(1+d_b)^2{\e^2}-\f12\|W^\e\|^2_{L^2(\Omega)}\\
&\qquad\qquad\qquad\qquad+C(\|U^\e\|_{L^2(\Omega)}+\|V^\e\|_{L^2(\Omega)})\|W^\e\|_{L^2(\Omega)}\\
&\le-d_b\e\| \nabla W^\e\|^2_{L^2(\Omega)}+C(1+d_b)^2{\e^2}-\f14\|W^\e\|^2_{L^2(\Omega)}+ C(\|U^\e\|^2_{L^2(\Omega)}+\|V^\e\|^2_{L^2(\Omega)}).
\end{split}
\eeq

\medskip

\noindent{\bf Step 4. The evolution equation of $\| \nabla W^\e\|^2_{L^2(\Omega)}$.}  Multiplying the equation for $W^\e$ in \eqref{U,V,W system} by $-\Delta W^\e$, we have
\begin{align}
\f12\f d{dt}\|\nabla W^\e\|^2_{L^2(\Omega)}&=-d_b\|\Delta W^\e\|^2_{L^2(\Omega)}+J_1+J_2+J_3\,,\label{est1 grad W}\\
J_1&:=-\int_\Omega \caFw(U^\e+u,W^\e+w,V^\e+v)\,\Delta W^\e\,dx\notag\\
J_2&:=\int_\Omega (\partial_tw-d_b\,\Delta\,w)\,\Delta W^\e\,dx\notag\\
J_3&:=\f1\e\int_\Omega\caQ(U^\e+u,W^\e+w,V^\e+v)\,\Delta W^\e\,dx\notag\,.
\end{align}
Using Taylor's formula into $J_1$ and the local boundedness of $D \caFw$, we have
\[
\begin{split}
J_1&=-\int_\Omega \caFw(u,w,v)\,\Delta W^\e\,dx-\int_\Omega (D\caFw(\xi_1,\xi_2,\xi_3)\cdot (U^\e,W^\e,V^\e))\,\Delta W^\e\,dx\\
&\le C\|\nabla W^\e\|_{L^1(\Omega)}+C(\|U^\e\|_{L^2(\Omega)}+\|W^\e\|_{L^2(\Omega)}+\|V^\e\|_{L^2(\Omega)})\|\Delta W^\e\|_{L^2(\Omega)}\,.
\end{split}
\]
Moreover, by the boundedness of $\nabla\partial_tw,\nabla\Delta\,w$, we obtain for $J_2$
\[
J_2\le C(1+d_b)\|\nabla W^\e\|_{L^1(\Omega)}\,.
\]
Therefore, 
\beq\label{est J1+J2}
\begin{split}
J_1&+J_2\le C(1+d_b)\|\nabla W^\e\|_{L^1(\Omega)}\\
&\qquad +C(\|U^\e\|_{L^2(\Omega)}+\|W^\e\|_{L^2(\Omega)}+\|V^\e\|_{L^2(\Omega)})\|\Delta W^\e\|_{L^2(\Omega)}\\
\le& C(1+d_b)^2\e+\f1{4\e}\|\nabla W^\e\|^2_{L^2(\Omega)}\\
&+C\,d_b^{-1}(\|U^\e\|^2_{L^2(\Omega)}+\|W^\e\|^2_{L^2(\Omega)}+\|V^\e\|^2_{L^2(\Omega)})+\f{d_b}2\|\Delta W^\e\|^2_{L^2(\Omega)}\,.
\end{split}
\eeq

The term $J_3$ is the more challenging and we proceed as for $I_2$ in Step 4, i.e. 
\begin{align}\label{def J32}
J_3=&\f1\e\int_\Omega [\caQ(U^\e+u,W^\e+w,V^\e+v)-\caQ(U^\e+u,w,V^\e+v)]\,\Delta W^\e\,dx\notag\\
&+\f1\e\int_\Omega[\caQ(U^\e+u,w,V^\e+v)-\caQ(u,w,v)]\,\Delta W^\e\,dx\notag\\
=:&J_3^1+J_3^2\,.
\end{align}
For $J_3^1$ we write
\[
J_3^1=-\f1\e\int_\Omega \nabla [\caQ(U^\e+u,W^\e+w,V^\e+v)-\caQ(U^\e+u,w,V^\e+v)]\cdot\nabla W^\e\,dx\,,
\]
so that 
\[
\begin{split}
J_3^1=&-\f1\e\int_\Omega [\pa_1\caQ(U^\e+u,W^\e+w,V^\e+v)-\pa_1\caQ(U^\e+u,w,V^\e+v)]\nabla(U^\e+u)\cdot\nabla W^\e\,dx\\
&-\f1\e\int_\Omega \pa_2\caQ(U^\e+u,W^\e+w,V^\e+v)\,|\nabla W^\e|^2\,dx\\
&-\f1\e\int_\Omega [\pa_2\caQ(U^\e+u,W^\e+w,V^\e+v)-\pa_2\caQ(U^\e+u,w,V^\e+v)]\nabla w\cdot\nabla W^\e\,dx\\
&-\f1\e\int_\Omega [\pa_3\caQ(U^\e+u,W^\e+w,V^\e+v)-\pa_3\caQ(U^\e+u,w,V^\e+v)]\nabla(V^\e+v)\cdot\nabla W^\e\,dx\,.
\end{split}
\]
Moreover, using $\pa_2\caQ=-\pa_1Q+\pa_2Q$ and \eqref{inf -pa1Q+pa2Q}, and rearranging the terms
\[
\begin{split}
J_3^1\le&-\f1\e\int_\Omega [\pa_1\caQ(U^\e+u,W^\e+w,V^\e+v)-\pa_1\caQ(U^\e+u,w,V^\e+v)]\nabla U^\e\cdot\nabla W^\e\,dx\\
&-\f1\e\int_\Omega \pa_{12}\caQ(U^\e+u,\xi_1,V^\e+v)\,W^\e (\nabla u\cdot\nabla W^\e)\,dx\\
&-\f1\e\|\nabla W^\e\|^2_{L^2(\Omega)}\\
&-\f1\e\int_\Omega \pa_{22}\caQ(U^\e+u,\xi_2,V^\e+v)\,W^\e(\nabla w\cdot\nabla W^\e)\,dx\\
&-\f1\e\int_\Omega [\pa_3\caQ(U^\e+u,W^\e+w,V^\e+v)-\pa_3\caQ(U^\e+u,w,V^\e+v)]\nabla V^\e\cdot\nabla W^\e\,dx\\
&-\f1\e\int_\Omega \pa_{32}\caQ(U^\e+u,\xi_3,V^\e+v)\,W^\e(\nabla v\cdot\nabla W^\e)\,dx\,.
\end{split}
\]
Next, using \eqref{sup pa1Q}, the boundedness of $\nabla u,\nabla v,\nabla w$ and the local boundedness of $DQ,D^2Q$, we obtain 
\beq\label{est J31}
\begin{split}
J_3^1\le& \f{2\caK_1}\e\|\nabla U^\e\|_{L^2(\Omega)}\,\|\nabla W^\e\|_{L^2(\Omega)}+\f C\e\| W^\e\|_{L^2(\Omega)}\,\|\nabla W^\e\|_{L^2(\Omega)}\\
&-\f1\e\|\nabla W^\e\|^2_{L^2(\Omega)}+\f C\e\|\nabla V^\e\|_{L^2(\Omega)}\,\|\nabla W^\e\|_{L^2(\Omega)}\\
\le& \f{6\caK_1^2}{\e}\|\nabla U^\e\|^2_{L^2(\Omega)}-\f1{2\e}\|\nabla W^\e\|^2_{L^2(\Omega)}+\f C{\e}\| W^\e\|^2_{L^2(\Omega)}+\f C{\e}\|\nabla V^\e\|^2_{L^2(\Omega)}\,.
\end{split}
\eeq

The term $J_3^2=\e^{-1}\int_\Omega[\caQ(U^\e+u,w,V^\e+v)-\caQ(u,w,v)]\,\Delta W^\e\,dx$ will be absorbed by the evolution equation of $E(t)$. Therefore, we let it as it is. Plugging into \eqref{est1 grad W} multiplied by $\e$, the estimates \eqref{est J1+J2}, \eqref{est J31} and the definition of $J_3^2$, we have
\beq\label{est2 grad W}
\begin{split}
\f{\e}2\f d{dt}&\|\nabla W^\e\|^2_{L^2(\Omega)}\le-\f{d_b\,\e}2\|\Delta W^\e\|^2_{L^2(\Omega)}+C(1+d_b)^2\e^2-\f14\|\nabla W^\e\|^2_{L^2(\Omega)}\\
&+6\caK_1^2\|\nabla U^\e\|^2_{L^2(\Omega)}+C\|\nabla V^\e\|^2_{L^2(\Omega)}\\
&+Cd_b^{-1}\e(\|U^\e\|^2_{L^2(\Omega)}+\|V^\e\|^2_{L^2(\Omega)})+C(d_b^{-1}\e+1)\|W^\e\|^2_{L^2(\Omega)}\\
&+\int_\Omega[\caQ(U^\e+u,w,V^\e+v)-\caQ(u,w,v)]\,\Delta W^\e\,dx\,.
\end{split}
\eeq

\noindent{\bf Step 5. The evolution equation of $E(t)$ defined in \eqref{def E}.} Let us write
\beq\label{eq E(t)}
\f d{dt}E(t)=\sum_{i=1}^5\,L_i(t)\,,
\eeq
with
\[
\begin{split}
&L_1:=-\int_\Omega[\pa_1\caP(U^\e+u,w,V^\e+v)-\pa_1\caP(u,w,v)]\,\pa_tU^\e\,dx\\
&L_2:=-\int_\Omega[\pa_1\caP(U^\e+u,w,V^\e+v)-\pa_1\caP(u,w,v)\\
&\qquad\qquad\qquad-\pa_{11}\caP(u,w,v)\,U^\e-\pa_{31}\caP(u,w,v)\,V^\e]\,\pa_tu\,dx\\
&L_3:=-\int_\Omega[\pa_2\caP(U^\e+u,w,V^\e+v)-\pa_2\caP(u,w,v)\\
&\qquad\qquad\qquad-\pa_{12}\caP(u,w,v)\,U^\e-\pa_{32}\caP(u,w,v)\,V^\e]\,\pa_tw\,dx
\end{split}
\]
\[
\begin{split}
&L_4:=-\int_\Omega[\pa_3\caP(U^\e+u,w,V^\e+v)-\pa_3\caP(u,w,v)]\,\pa_tV^\e\,dx\\
&L_5:=-\int_\Omega[\pa_3\caP(U^\e+u,w,V^\e+v)-\pa_3\caP(u,w,v)\\
&\qquad\qquad\qquad-\pa_{13}\caP(u,w,v)\,U^\e-\pa_{33}\caP(u,w,v)\,V^\e]\,\pa_tv\,dx\,.
\end{split}
\]

The terms $L_2,L_3,L_5$ are easily controlled using second order Taylor's formula applied to $\pa_1\caP=\caQ,\pa_2\caP,\pa_3\caP$ respectively, to obtain
\beq\label{est L2,L3,L5}
\begin{split}
&L_2\le C\int_\Omega((U^\e)^2+(V^\e)^2)|\pa_t u|\, dx\le C(\| U^\e\|^2_{L^2(\Omega)}+\| V^\e\|^2_{L^2(\Omega)})\,,\\
&L_3\le C\int_\Omega((U^\e)^2+(V^\e)^2)|\pa_t w|\, dx\le C(\| U^\e\|^2_{L^2(\Omega)}+\| V^\e\|^2_{L^2(\Omega)})\,,\\
&L_5\le C\int_\Omega((U^\e)^2+(V^\e)^2)|\pa_t v|\, dx\le C(\| U^\e\|^2_{L^2(\Omega)}+\| V^\e\|^2_{L^2(\Omega)})\,.
\end{split}
\eeq

Next, by \eqref{def caP} and the equation for $U^\e$ in \eqref{U,V,W system}, we have for $L_1$
\beq\label{def L1}
\begin{split}
L_1=&-\int_\Omega[\caQ(U^\e+u,w,V^\e+v)-\caQ(u,w,v)][d_a\Delta U^\e+\delta\Delta W^\e\\
&\qquad\qquad+\caFu(U^\e+u,W^\e+w,V^\e+v)-\caFu(u,w,v)]\,dx\\
=&:L_1^1+L_1^2+L_1^3\,.
\end{split}
\eeq

The $L_1^2$ term is fundamental since it allows us to get rid of the term $J_3^2$ defined in \eqref{def J32}. Indeed, $L_1^2$ reads as
\beq\label{est L12}
L_1^2=-\delta \int_\Omega[\caQ(U^\e+u,w,V^\e+v)-\caQ(u,w,v)]\,\Delta W^\e\,dx=-\e\,\delta\,J_3^2\,.
\eeq

The control of $L_1^3$ follows simply as
\begin{align}
L_1^3=&-\int_\Omega[\caQ(U^\e+u,w,V^\e+v)-\caQ(u,w,v)][\caFu(U^\e+u,W^\e+w,V^\e+v)-\caFu(u,w,v)]\,dx\notag\\
&\le C\int_\Omega(|U^\e|+|V^\e|)(|U^\e|+|V^\e|+|W^\e|)\,dx\label{est L13}\\
&\le C(\| U^\e\|^2_{L^2(\Omega)}+\| V^\e\|^2_{L^2(\Omega)})+\f18\| W^\e\|^2_{L^2(\Omega)}\notag\,.
\end{align}

The control of $L_1^1$ follows the same computations done for $J_3^1$ defined in \eqref{def J32}. First we write
\[
\begin{split}
L_1^1=&\ d_a\int_\Omega \nabla [\caQ(U^\e+u,w,V^\e+v)-\caQ(u,w,v)]\cdot\nabla U^\e\,dx\\
=&\ d_a\int_\Omega\pa_1\caQ(U^\e+u,w,V^\e+v)\,|\nabla U^\e|^2\,dx\\
&+d_a\int_\Omega [\pa_1\caQ(U^\e+u,w,V^\e+v)-\pa_1\caQ(u,w,v)]\,\nabla u\cdot\nabla U^\e\,dx\\
&+d_a\int_\Omega [\pa_2\caQ(U^\e+u,w,V^\e+v)-\pa_2\caQ(u,w,v)]\,\nabla w\cdot\nabla U^\e\,dx\\
&+d_a\int_\Omega\pa_3\caQ(U^\e+u,w,V^\e+v)\,\nabla V^\e\cdot\nabla U^\e\,dx\\
&+d_a\int_\Omega [\pa_3\caQ(U^\e+u,w,V^\e+v)-\pa_3\caQ(u,w,v)]\,\nabla v\cdot\nabla U^\e\,dx\,.
\end{split}
\]
Next, since $\pa_1\caQ=\pa_1Q$, using \eqref{inf pa1Q} and skipping few details, we have
\[
\begin{split}
L_1^1\le&-d_a\,\caK_0\|\nabla U^\e\|^2_{L^2(\Omega)}
+d_a\,C\int_\Omega(|U^\e|+|V^\e|)(|\nabla u|+|\nabla w|+|\nabla v|)\,|\nabla U^\e|\,dx\\
&+d_a\,C\|\nabla U^\e\|_{L^2(\Omega)}\|\nabla V^\e\|_{L^2(\Omega)}\,,
\end{split}
\]
and by the boundedness of $\nabla u$, $\nabla w$, $\nabla v$ and ad hoc Young's inequalities, we obtain
\beq\label{est L11}
L_1^1\le-\f{d_a\caK_0}2\,\|\nabla U^\e\|^2_{L^2(\Omega)}+\f C{\caK_0}(\| U^\e\|^2_{L^2(\Omega)}+\| V^\e\|^2_{L^2(\Omega)})+C\f {d_a}{\caK_0}\|\nabla V^\e\|^2_{L^2(\Omega)}\,.
\eeq
Plugging \eqref{est L12}, \eqref{est L13}, \eqref{est L11} into \eqref{def L1} we end up with
\beq\label{est L1}
\begin{split}
L_1\le&-\f{d_a\caK_0}2\,\|\nabla U^\e\|^2_{L^2(\Omega)}+C\f {d_a}{\caK_0}\,\|\nabla V^\e\|^2_{L^2(\Omega)}\\
&+C(1+\caK_0^{-1})(\| U^\e\|^2_{L^2(\Omega)}+\| V^\e\|^2_{L^2(\Omega)})+\f18\| W^\e\|^2_{L^2(\Omega)}\\
&-\delta \int_\Omega[\caQ(U^\e+u,w,V^\e+v)-\caQ(u,w,v)]\,\Delta W^\e\,dx\,.
\end{split}
\eeq

We consider now the term $L_4$ and use the equation for $V^\e$ in \eqref{U,V,W system}
\beq\label{def L4}
\begin{split}
L_4=&-\int_\Omega[\pa_3\caP(U^\e+u,w,V^\e+v)-\pa_3\caP(u,w,v)][d_v\,\Delta V^\e\\
&\qquad+\caFv(U^\e+u,W^\e+w,V^\e+v)-\caFv(u,w,v)]\,dx\\
=&:L_4^1+L_4^2\,.
\end{split}
\eeq
For the $L_4^1$ term we have 
\[
\begin{split}
L_4^1=&d_v\int_\Omega\nabla[\pa_3\caP(U^\e+u,w,V^\e+v)-\pa_3\caP(u,w,v)]\cdot\nabla V^\e\,dx\\
=&d_v\int_\Omega\pa_{31}\caP(U^\e+u,w,V^\e+v)\nabla U^\e\cdot\nabla V^\e\,dx\\
&+d_v\int_\Omega[\pa_{31}\caP(U^\e+u,w,V^\e+v)-\pa_{31}\caP(u,w,v)]\nabla u\cdot\nabla V^\e\,dx\\
&+d_v\int_\Omega[\pa_{32}\caP(U^\e+u,w,V^\e+v)-\pa_{32}\caP(u,w,v)]\nabla w\cdot\nabla V^\e\,dx\\
&+d_v\int_\Omega\pa_{33}\caP(U^\e+u,w,V^\e+v)\,|\nabla V^\e|^2\,dx\\
&+d_v\int_\Omega[\pa_{33}\caP(U^\e+u,w,V^\e+v)-\pa_{33}\caP(u,w,v)]\nabla v\cdot\nabla V^\e\,dx\,,
\end{split}
\]
so that
\beq\label{est L41}
\begin{split}
L_4^1\le&\f{d_a\caK_0}4\,\|\nabla U^\e\|^2_{L^2(\Omega)}+C\f{d_v^2}{d_a\caK_0}\|\nabla V^\e\|^2_{L^2(\Omega)}\\
&+Cd_v\int_\Omega(|U^\e|+|V^\e|)(|\nabla u|+|\nabla w|+|\nabla v|)\,|\nabla V^\e|\,dx+Cd_v\|\nabla V^\e\|^2_{L^2(\Omega)}\\
\le&\f{d_a\caK_0}4\,\|\nabla U^\e\|^2_{L^2(\Omega)}+Cd_v(1+d_v(d_a\caK_0)^{-1})\|\nabla V^\e\|^2_{L^2(\Omega)}\\
&+Cd_v(\| U^\e\|^2_{L^2(\Omega)}+\| V^\e\|^2_{L^2(\Omega)}).
\end{split}
\eeq
For the $L_4^2$ term in \eqref{def L4}, skipping again few details, we have 
\beq\label{est L42}
\begin{split}
L_4^2&\le C\int_\Omega(|U^\e|+|V^\e|)(|U^\e|+|V^\e|+|W^\e|)\,dx\\
&\le C(\| U^\e\|^2_{L^2(\Omega)}+\| V^\e\|^2_{L^2(\Omega)})+\f18\| W^\e\|^2_{L^2(\Omega)}\,.
\end{split}
\eeq
Plugging \eqref{est L41}, \eqref{est L42} into \eqref{def L4} we obtain
\beq\label{est L4}
\begin{split}
L_4\le&\f{d_a\caK_0}4\,\|\nabla U^\e\|^2_{L^2(\Omega)}+Cd_v(1+d_v(d_a\caK_0)^{-1})\|\nabla V^\e\|^2_{L^2(\Omega)}\\
&\qquad+C(1+d_v)(\| U^\e\|^2_{L^2(\Omega)}+\| V^\e\|^2_{L^2(\Omega)})+\f18\| W^\e\|^2_{L^2(\Omega)}\,.
\end{split}
\eeq

Gathering \eqref{est L2,L3,L5}, \eqref{est L1}, \eqref{est L4} in \eqref{eq E(t)}, we end up with the estimate of $E'(t)$
\begin{align}
\f d{dt}E(t)\le&-\f{d_a\caK_0}4\,\|\nabla U^\e\|^2_{L^2(\Omega)}+C(d_v+d_v^2(d_a\caK_0)^{-1}+d_a\caK_0^{-1})\|\nabla V^\e\|^2_{L^2(\Omega)}\notag\\
&+C(1+d_v+\caK_0^{-1})(\| U^\e\|^2_{L^2(\Omega)}+\| V^\e\|^2_{L^2(\Omega)})+\f14\| W^\e\|^2_{L^2(\Omega)}\notag\\
&-\delta \int_\Omega[\caQ(U^\e+u,w,V^\e+v)-\caQ(u,w,v)]\,\Delta W^\e\,dx\label{est E}\,.
\end{align}

\noindent{\bf Step 6. A positive lower bound.} Let $T>0$ and $\gamma_1,\gamma_2\in\R_+^2$. We claim that, if $\gamma_2>0$ is large enough, there exists $c_0=c_0(\gamma_1,\caK_0,T)>0$ such that
\beq\label{lower bound E}
E(t)+\f{\gamma_1}2\| U^\e\|^2_{L^2(\Omega)}+\f{\gamma_2}2\| V^\e\|^2_{L^2(\Omega)}\ge c_0(\| U^\e\|^2_{L^2(\Omega)}+\| V^\e\|^2_{L^2(\Omega)})\,.
\eeq
Indeed, from \eqref{def E} and a second order Taylor's formula, we have
\[
E(t)\!=\!-\int_\Omega [\pa_{11}\caP(\xi_1,w,\xi_2)\,(U^\e)^2+2\pa_{13}\caP(\xi_1,w,\xi_2)U^\e V^\e+\pa_{33}\caP(\xi_1,w,\xi_2)\,(V^\e)^2]dx.
\]
Moreover, from \eqref{new reactions}, \eqref{def caP} and \eqref{inf pa1Q}, we have
\beq\label{lower bound}
-\pa_{11}\caP(\xi_1,w,\xi_2)=-\pa_{1}\caQ(\xi_1,w,\xi_2)=-\pa_{1}Q(\xi_1,w,\xi_2)\ge\caK_0>0\,.
\eeq
Therefore, using \eqref{lower bound} and the local boundedness of $DQ$ and $D^2Q$, we obtain 
\[
E(t)\ge \caK_0\| U^\e\|^2_{L^2(\Omega)}-C_1\| U^\e\|_{L^2(\Omega)}\| V^\e\|_{L^2(\Omega)}-C_2\| V^\e\|^2_{L^2(\Omega)}\,,
\]
and, for $\theta>0$ to be chosen,
\[
\begin{split}
E(t)+\f{\gamma_1}2\| U^\e\|^2_{L^2(\Omega)}+\f{\gamma_2}2\| V^\e\|^2_{L^2(\Omega)}\ge& (\caK_0+\f{\gamma_1}2-\theta)\| U^\e\|^2_{L^2(\Omega)}\\
&+(\f{\gamma_2}2-C_2-\f{C_1^2}{4\theta})\| V^\e\|^2_{L^2(\Omega)}\,.
\end{split}
\]
In order to obtain \eqref{lower bound E} it is sufficient to choose $\theta>0$ small enough so that $\caK_0+\f{\gamma_1}2-\theta>0$ and then $\gamma_2>0$ large enough so that $\f{\gamma_2}2-C_2-\f{C_1^2}{4\theta}>0$. 
\medskip

\noindent{\bf Step 7. End of the proof.} Recall definition \eqref{def L} of $\caL(t)$.  By \eqref{L2 norm U}, \eqref{L2 norm V}, \eqref{est1 W},\eqref{est2 grad W} and \eqref{est E}, neglecting useless negative terms, we have
\begin{align}
\f d{dt}\caL(t)\le& c_1(\gamma_1)\| \nabla U^\e\|^2_{L^2(\Omega)}+c_2(\gamma_2)\| \nabla V^\e\|^2_{L^2(\Omega)}+c_3(\gamma_1)\| \nabla W^\e\|^2_{L^2(\Omega)}\notag\\
&+c_4(\gamma_1,\gamma_2,\gamma_3)(\| U^\e\|^2_{L^2(\Omega)}+\| V^\e\|^2_{L^2(\Omega)})\notag\\
&+c_5(\gamma_1,\gamma_2,\gamma_3)\|W^\e\|^2_{L^2(\Omega)}+C(1+d_b)^2(\gamma_3+\delta)\e^2\label{ed L(t)}\,,
\end{align}
where
\beq\label{constants}
\begin{split}
&c_1(\gamma_1)=-\f{d_a}4\,\caK_0-\f{d_a}2\,\gamma_1+{6\delta}\,\caK_1^2\,,\\
&c_2(\gamma_2)=-d_v\,\gamma_2+C(\delta+d_v+d_v^2(d_a\caK_0)^{-1}+d_a\caK_0^{-1})\,,\\
&c_3(\gamma_1)=\f{\delta^2}{2d_a}\,\gamma_1-\f\delta4\,,\\
&c_5(\gamma_1,\gamma_2,\gamma_3)=\f14(1+\gamma_1+\gamma_2-\gamma_3)+\delta\,C(1+\e\,d_b^{-1})\,.
\end{split}
\eeq
We need now to determine $\gamma_1,\gamma_2,\gamma_3\in\R_+$ (independent of $\e$) such that the constant $c_1,c_2,c_3,c_5$ are negative. We will see that this is possible provided that the ratio $d_b/d_a$ is small enough (see  \eqref{hyp rate of convergence}). Since the constant $c_4$ is positive for all $\gamma_1,\gamma_2,\gamma_3\!\in\!\R_+$, we do not need the explicit formula. 

Note that $c_1$ and $c_3$ in \eqref{constants} are both negative if and only~if
\beq\label{condition gamma1}
\f{12\caK_1^2\delta}{d_a}-\f{\caK_0}2<\gamma_1<\f{d_a}{2\delta}\,.
\eeq
If $\f{12\caK_1^2\delta}{d_a}-\f{\caK_0}2\le0$, i.e. $\f{d_b}{d_a}\le 1+\f{\caK_0}{24\caK_1^2}$, then it is sufficient to choose $\gamma_1\in(0,\f{d_a}{2\delta})$ to have $c_1,c_3$ negative. On the other hand,  if $\f{12\caK_1^2\delta}{d_a}-\f{\caK_0}2>0$, then it is easily seen that $\f{12\caK_1^2\delta}{d_a}-\f{\caK_0}2<\f{d_a}{2\delta}$ if $\f{d_b}{d_a}<1+x_+$, where $x_+$ is the positive root of the polynomial function $(12\caK_1^2\,x^2-\f{\caK_0}2\,x-\f12)$ and satisfies $x_+>\f{\caK_0}{24\caK_1^2}$. Therefore, it is again possible to find $\gamma_1>0$ satisfying \eqref{condition gamma1} and giving negative $c_1$ and $c_3$. 

For this $\gamma_1\!>\!0$, we choose $\gamma_2\!>\!0$ large enough so that at the same time \eqref{lower bound E} holds true and $c_2$ is negative, and finally $\gamma_3\!>\!0$ large enough so that  $c_5$ is negative as well.

Now, by \eqref{ed L(t)} and \eqref{lower bound E}, $\caL(t)$ satisfies
\[
\begin{split}
\f d{dt}\caL(t)\le&c_4\,(\| U^\e\|^2_{L^2(\Omega)}+\| V^\e\|^2_{L^2(\Omega)})+C(1+d_b)^2\e^2(\gamma_3+\delta)\\
\le&c_4\,c_0^{-1}\,(E(t)+\f{\gamma_1}2\| U^\e\|^2_{L^2(\Omega)}+\f{\gamma_2}2\| V^\e\|^2_{L^2(\Omega)})+C(1+d_b)^2\e^2(\gamma_3+\delta)\\
\le&c_4\,c_0^{-1}\,\caL(t)+C(1+d_b)^2(\gamma_3+\delta)\e^2\,.
\end{split}
\]
Integrating the differential inequality above over $(0,t)$, we obtain that there exists $C(T)>0$ such that, for all $t\in(0,T)$, it holds
\[
E(t)+\f{\gamma_1}2\| U^\e(t)\|^2_{L^2(\Omega)}+\f{\gamma_2}2\| V^\e(t)\|^2_{L^2(\Omega)}\le\caL(t)\le C(T)(\e^2+\caL(0))\,.
\]
Observing that (see \eqref{def L}, \eqref{def Ue,Ve,We} and \eqref{initial layer})
\beq\label{L(0)}
\caL(0)=\f{\e\gamma_3}2\| W^\e(0)\|^2_{L^2(\Omega)}+\f{\e\delta}2\|\nabla W^\e(0)\|^2_{L^2(\Omega)}=C(T)\,\e\,\e_\init^2\,,
\eeq
by \eqref{lower bound E} again, the above inequality implies
\beq\label{U+V}
\| U^\e(t)\|^2_{L^2(\Omega)}+\| V^\e(t)\|^2_{L^2(\Omega)}\le C(T)(\e^2+\e\,\e_\init^2)\,,\quad t\in(0,T)\,.
\eeq
Next, we plug \eqref{U+V} into \eqref{est1 W}, we neglect useless negative terms and we obtain
\[
\f d{dt}\| W^\e(t)\|^2_{L^2(\Omega)}\le C(T)(\e+\e_\init^2)-\f1{2\e}\| W^\e(t)\|^2_{L^2(\Omega)}\,,
\]
i.e. 
\beq\label{est W final}
\| W^\e(t)\|^2_{L^2(\Omega)}\le \e_\init^2\,\mathrm{e}^{-\f1{2\e}t}+C(T)(\e^2+\e\,\e_\init^2)\,,\quad t\in(0,T)\,.
\eeq
Finally, we plug \eqref{U+V}, \eqref{est W final} into \eqref{ed L(t)} to obtain
\beq\label{final est L}
\begin{split}
\f d{dt}\caL(t)-c_1\| \nabla U^\e(t)\|^2_{L^2(\Omega)}&-c_2\| \nabla V^\e(t)\|^2_{L^2(\Omega)}-c_3\| \nabla W^\e(t)\|^2_{L^2(\Omega)}\\
&\le \e_\init^2\,\mathrm{e}^{-\f1{2\e}t}+C(T)(\e^2+\e\,\e_\init^2)\,.
\end{split}
\eeq

We are now able to conclude. Indeed, \eqref{u-ue v-ve} follows by  \eqref{U+V} and \eqref{final est L} integrated over $(0,t)$ and taking into account the positivity of $\caL(t)$, the negativity of $c_1,c_2,c_3$ and \eqref{L(0)}. \eqref{uae-ua* ube-ub*} follows by \eqref{est W final}, \eqref{final est L} and $u_a^\e-u_a^*(u,v)=(u^\e-u)-W^\e$. 
\appendix
\section{Proof of Theorem \ref{th exist meso}}\label{appendix A}
Throughout the proof, we omit the $\e$ superscript for the sake of clarity. 

First, we truncate the reaction functions in \eqref{meso system} in order to obtain globally Lipschitz functions in $L^p(\Omega)$. To do this, fix $M>0$. Let us denote $U=(u_a,u_b,v)$, $|U|=\max({|u_a|,|u_b|,|v|})$ and let 
\[
f_\nu ^M(U):=\left\{
\begin{split}
&f_\nu (U)\,,\qquad\quad\quad\ \ \text{if }|U|\le M\,,\\
&f_\nu \left(M\f U{|U|}\right)\,,\quad\quad\text{if }|U|> M\,,
\end{split}
\right.
\]
for $\nu  =a, b, v$, where $f_\nu $ is defined by \eqref{def react functions}. Moreover, let
\[
\psi^M(x):=\left\{
\begin{split}
& A^\alpha\,,\qquad\text{if }x<0\,,\\
&\psi(x)\,,\quad\ \text{if }x\in[0,M]\,,\\
&\psi(M)\,,\quad\text{if }x>M\,,
\end{split}
\right.
\quad
\phi^M(x):=\left\{
\begin{split}
& B^\beta\,,\qquad\text{if }x<0\,,\\
&\phi(x)\,,\quad\ \text{if }x\in[0,M]\,,\\
&\phi(M)\,,\quad\text{if }x>M\,,
\end{split}
\right.
\]
where $\psi,\phi$ are the transiction functions defined by \eqref{def phi psi}. 

Next, we let $\Lambda^M_a, \Lambda^M_b$ be the relative satisfaction measures defined by \eqref{def rel sat measures} with $\psi,\phi$ replaced by $\psi^M,\phi^M$ and also $Q^M(U)=\Lambda^M_bg^M(u_b)-\Lambda^M_ag^M(u_a)$, where
\[
g^M(x):=
\left\{
\begin{split}
&x\,,\quad\qquad\ \text{if }|x|\le M\,,\\
&M\f{x}{|x|}\,,\quad\ \,\text{if }|x|> M\,.
\end{split}
\right.
\]
Finally, we define the nonlinear mapping $F^M:\R^3\to\R^3$ as
\[
F^M(U):=(f_a^M(U)+\e^{-1}Q^M(U), f_b^M(U)-\e^{-1}Q^M(U), f_v^M)\,.
\]
Note that the functions $f_\nu^M,\psi^M,\phi^M$ are globally Lipschitz and bounded. Therefore, $F^M$ is globally Lipschitz, bounded and  
\[
F^M(U):=(f_a(U)+\e^{-1}Q(U), f_b(U)-\e^{-1}Q(U), f_v(U))\,,\quad\text{if}\quad |U|\le M\,.
\]

Given $p\in(1,+\infty)$, we consider the operator $A_p$ on $X_p:=(L^p(\Omega))^3$ defined by
\beq\label{def operator A}
\left\{
\begin{split}
&D(A_p)=D_p^3\qquad\text{(see \eqref{def Dp})}\\
&A_pU=(d_a\Delta u_a,d_b\Delta u_b,d_v\Delta v)\quad\text{for }U\in D(A_p)\,,
\end{split}
\right.
\eeq
and the abstract initial value problem
\beq\label{abstract IVP}
U'(t)=A_pU(t)+F^M(U(t))\,,\quad t>0\,,\quad U(0)=U^\init:=(u_a^\init,u_b^\init,v^\init)\,.
\eeq

We will solve \eqref{def operator A},\eqref{abstract IVP} and then we will get rid of the truncation. The main ingredient is that $A_p\!:\!D(A_p)\!\subset \!X_p\to X_p$ is a sectorial operator (\cite{Lunardi1995}, Theorem~3.1.3). Hence it generates in $X_p$ an analytic semigroup denoted $(e^{tA_p})_{t\ge0}$ (\cite{Lunardi1995}, Chapter~2). Moreover, $A_p$ is closed so that $D(A_p)$, endowed with the graph norm, is a Banach space. $D(A_p)$ being also dense in $X_p$, the semigroup is strongly continuous, i.e. $\lim_{t\to0}e^{tA_p}U=U$, for all $U\in X_p$. Furthermore, there exists $K_p\!>\!0$ and $\omega_p\!\in\!\R$ such that (see \cite{Lunardi1995}, Proposition 2.1.1)
\beq\label{bound etA}
\|e^{tA}\|_{L(X_p)}\le K_pe^{\omega_p\,t}\,,\qquad\forall\ t\ge0\,.
\eeq

\textit{First step: well-posedness of \eqref{def operator A},\eqref{abstract IVP}}.
Let $\|\cdot\|_p$ denote the usual norm in $X_p$. We start proving that \eqref{def operator A},\eqref{abstract IVP} has a unique mild solution, i.e. a unique function $U\in C^0([0,\infty),X_p)$ such that 
\beq\label{def mild sol}
U(t)=e^{tA}U^\init+\int_0^te^{(t-s)A}F^M(U(s))\,ds\,,\qquad\forall\ t\ge0\,.
\eeq

It is easily seen that $F^M$ maps $X_p$ into $X_p$ and 
\beq\label{Lip FM}
\|F^M(U)\|_p\le \|F^M(0)\|_p\!+\!L_M\|U\|_p=L_M\|U\|_p\,,\quad\forall\ U\in X_p\,.
\eeq
Therefore, \eqref{def mild sol} makes sense since, by assumption \eqref{hyp ID}, $U^\init\in X_p$ and, for all $U\in C^0([0,\infty),X_p)$ and all $t>0$, $F^M(U(\cdot))\in L^1((0,t);X_p)$. Moreover, the Lipschitz property of $F^M$ together with \eqref{bound etA} and Gronwall's Lemma gives us the uniqueness of \eqref{def mild sol}. The same ingredients give us the continuous dependence of $U$ with respect to $U^\init$. Therefore, it remains to prove the existence of $U$ and that $U$ belongs to $C^1([0,\infty);X_p)\cap C^0([0,\infty);D_p(A))$ for all $p\in(1,+\infty)$. 

Let $\theta>0$ be such that $\omega_p+\theta>0$. The existence is proved using the contraction mapping principle in the space 
\[
E:=\{U\in C^0([0,\infty),X_p) : \|U\|_E=\sup_{t\ge0}e^{-(\omega_p+\theta)t}\|U(t)\|_{p}<\infty\}\,,
\]
that is a Banach space when endowed with the norm $\|U\|_E$. Hence, given $U\in E$, we set
\[
\Phi(U)(t)=e^{tA}U^\init+\int_0^te^{(t-s)A}F^M(U(s))\,ds\,,\quad \forall t\ge0\,.
\]
We claim that $\Phi$ maps $E$ into $E$ and it is a contraction provided that $\theta>K_pL_M$. Indeed, it is clear that $\Phi(U)\in C([0,\infty),X_p)$. Moreover, using \eqref{bound etA} and \eqref{Lip FM} or the Lipschitz property of $F^M$, we obtain that, for all $U,V\in E$, 
\[
\|\Phi(U)\|_E\le K_p\|U^\init\|_p+K_p\theta^{-1}L_M\|U\|_E\,, 
\]
and
\[
\|\Phi(U)-\Phi(V)\|_E\le K_p\theta^{-1}L_M\|U-V\|_E\,. 
\]

Now, for all $T\in(0,\infty)$, the mild solution $U$ is Lipschitz $[0,T]\to X_p$.
Indeed, on the one hand, using \eqref{bound etA}, the Lipschitz property of $F^M$ and Gronwall's Lemma again, we have
\beq\label{U Lip 1}
\|U(t)\|_p\le K\,e^{(K_pL_M+\omega_p)t}\|U^\init\|_p\,,\qquad \forall t\ge0\,.
\eeq
On the other hand, as $U(t+h)$, $t,h>0$, is a mild solution of \eqref{abstract IVP} with initial data $U(h)\in X_p$, proceeding as above, we obtain
\beq\label{Lip U(t) basic}
\|U(t+h)-U(t)\|_p\le K\,e^{(K_pL_M+\omega_p)t}\|U(h)-U^\init\|_p\,,\qquad\forall\ t\ge0\,.
\eeq
Next, by \cite{Lunardi1995} Propositions 2.1.1, 2.1.4, and since $U^\init\in D_p(A)$ (see \eqref{hyp ID}),  it holds
\[
e^{hA}U^\init-U^\init=A\int_0^h e^{sA}U^\init\,ds=\int_0^h A\,e^{sA}U^\init\,ds=\int_0^h e^{sA}\,AU^\init\,ds\,.
\]
Therefore, from \eqref{bound etA}-\eqref{Lip FM} and the above equality, we get 
\beq\label{U Lip 2}
\|U(h)-U^\init\|_p\le K_pe^{\omega_p h}h\left(\|AU^\init\|_p+L_M\sup_{0\le s\le h}\|U(s)\|_p\right)\,,
\eeq
and the Lipschitz property follows from \eqref{U Lip 1}, \eqref{Lip U(t) basic}, \eqref{U Lip 2}. 

As a consequence, $F^M(U(\cdot))$ is Lipschitz $[0,T]\to X_p$. Therefore, taking also into account that $F^M$ is bounded, and applying Theorem 4.3.1 (ii) and Lemma 4.1.6 in \cite{Lunardi1995}, we have that $U\!\in C^1([0,\infty);X_p)\cap C^0([0,\infty);D(A_p))$. 

Finally, note that, since $\Omega$ is bounded, if $p\ge q$, $X_p\subset X_q$ (with continuous embedding), $D(A_p)\subset D(A_q)$ and $A_pU=A_qU$ if $U\in D(A_p)$. Therefore, by assumption \eqref{hyp ID}, the above time regularity holds true for all $p\in(1,+\infty)$. We can drop the subscript $p$ in the sequel. 
\medskip

\textit{Second step: well-posedness of \eqref{meso system}--\eqref{def phi psi}.}
Let $U\!=\!(u_a,u_b,v)$ be the unique solution of \eqref{def operator A},\eqref{abstract IVP}. Then, $u_a,u_b>0$ and $v\ge0$ on $(0,+\infty)\times\Omega$, since the semi-group $(e^{tA})_{t\ge0}$ is strongly positive, the initial data are non-negative, with $u_a^\init,u_b^\init$ not identically zero, and the nonlinear mapping $F^M$ is quasi-positive (see \cite{Haraux2017}). 

Multiplying the equation for $u_a$ in \eqref{abstract IVP} by $u_a^{p-1}$, $p>1$, using the positivity of the solution  in $f_a^M$ and $Q^M$, the fact that the satisfaction measure $\Lambda_b^M$ lies in $[0,1)$ and the Young inequality, we get
\[
\begin{split}
\f1p\f d{dt}\int_\Omega u_a^p\,dx&\le\eta_a\int_\Omega u_a^p\,dx+\f1\e\int_\Omega u_b\, u_a^{p-1}\,dx\\
&\le\overline\eta\int_\Omega u_a^p\,dx+\f1{\e p}\int_\Omega u_b^p\,dx+\f1\e(1-\f1p)\int_\Omega u_a^p\,dx\,.
\end{split}
\]
Similarly, for $u_b$ we have 
\[
\f1p\f d{dt}\int_\Omega u_b^p\,dx\le\overline\eta\int_\Omega u_b^p\,dx+\f1{\e p}\int_\Omega u_a^p\,dx+\f1\e(1-\f1p)\int_\Omega u_b^p\,dx\,,
\]
so that, for all $t>0$, 
\[
\f d{dt}(\|u_a(t)\|^p_{L^p(\Omega)}+\|u_b(t)\|^p_{L^p(\Omega)})
\le p(\overline\eta+\f1\e)(\|u_a(t)\|^p_{L^p(\Omega)}+\|u_b(t)\|^p_{L^p(\Omega)})\,.
\]
Integrating the above differential inequality over $(0,T)$, we end up with 
\[
\sup_{0\le t\le T}(\|u_a\|^p_{L^p(\Omega)}+\|u_b\|^p_{L^p(\Omega)})\le(\|u_a^\init\|^p_{L^p(\Omega)}+\|u_b^\init\|^p_{L^p(\Omega)})^{\f1p}\,e^{(\overline\eta+\f1\e)\, T}\,,
\]
implying, as $p\to\infty$ and for all $T>0$, 
\[
\|u_a\|_{L^{\infty}(0,T;L^\infty(\Omega))},\|u_b\|_{L^{\infty}(0,T;L^\infty(\Omega))}\le
(\|u_a^\init\|_{L^\infty(\Omega)}\vee\|u_b^\init\|_{L^\infty(\Omega)})\,e^{(\overline\eta+\f1\e)\, T}:=M_1.
\]
Similarly, we prove that
\[
\|v\|_{L^{\infty}(0,T;L^\infty(\Omega))}\le \|v^\init\|_{L^\infty(\Omega)}\,e^{\eta_v T}:=M_2.
\]

We now fix $T>0$ and we let $M_T=M_1\vee M_2$. It follows from the above estimates that $|U(t,x)|\le M_T$, for all $(t,x)\in[0,T]\times\Omega$. Therefore, $F^{M_T}(U):=(f_a(U)+\e^{-1}Q(U), f_b(U)-\e^{-1}Q(U), f_v(U))$, on $[0,T]\times\Omega$. Thus, we see that $U$ satisfies \eqref{meso system} on $[0,T]\times\Omega$. We denote this solution by $U^T$, since at this stage it might depend on $T$. We claim that in fact it does not. Indeed, let $0< T< S$. On $[0,T]$, both $U^T$ and $U^S$ are bounded. Therefore, both $U^T$ and $U^S$ satisfy the same equation \eqref{abstract IVP} provided $M$ is chosen sufficiently large. By uniqueness for the equation \eqref{abstract IVP}, it follows that $U^T=U^S$ on $[0,T]\times\Omega$. Thus we obtain a solution $U$ of \eqref{meso system} on $[0,\infty)\times\Omega$. 
\medskip

\textit{Third step: estimates uniform in $\e$}.
It is easily seen from \eqref{meso system}, \eqref{def react functions} (using the positivity of the solution and \eqref{notations 1}) that $u_a+u_b$ satisfies
\beq\label{norm(ua+ub)L1}
\f d{dt}\|u_a+u_b\|_{L^1(\Omega)}\le\overline\eta\|u_a+u_b\|_{L^1(\Omega)}-a\eta_a\| u_a\|_{L^2(\Omega)}^2-b\eta_b\|u_b\|_{L^2(\Omega)}^2\,,
\eeq
so that $y(t)=\|(u_a+u_b)(t)\|_{L^1(\Omega)}$ satisfies
\[
y'(t)\le\overline\eta\,y(t)-\f\eta{2|\Omega|}\,y^2(t)\,,\qquad \forall\ t\ge0\,.
\]
Integrating the above differential inequality, we obtain \eqref{LinftyL1}. Furthermore, integrating \eqref{norm(ua+ub)L1} over $(0,T)$, we obtain
\[
a\eta_a\| u_a\|_{L^2(0,T;L^2(\Omega))}^2+b\eta_b\|u_b\|_{L^2(0,T;L^2(\Omega))}^2
\le \| u^{\init}_a+u^{\init}_b\|_{L^1(\Omega)}+\overline\eta\int_0^{T}y(t)\,dt\,.
\]
Using \eqref{LinftyL1} in the above inequality, we get \eqref{result p=1}. 

Next, multiplying the equation for $v$ in  \eqref{meso system} by $v^{p-1}$, $p>1$, and using again the positivity of the solution in $f_v$, we get
\beq\label{eq vp}
\f1p\f d{dt}\int_\Omega v^p\,dx\le\eta_v\int_\Omega v^p\,dx-r_v\int_\Omega v^{p+1}\,dx\,.
\eeq
Plugging into \eqref{eq vp} the H\"older inequality
\[
\int_\Omega v^p\,dx\le \Big(\int_\Omega v^{p+1}\,dx\Big)^{\f p{p+1}}|\Omega|^\f1{p+1},
\]
we have
\beq\label{ineq vp}
\begin{split}
\f1p\f d{dt}\int_\Omega v^p\,dx\le
\eta_v\,\int_\Omega v^p\,dx-\f{r_v}{|\Omega|^{1/p}}\,\Big(\int_\Omega v^{p}\,dx\Big)^{1+\f1p}\,.
\end{split}
\eeq
Integrating \eqref{ineq vp} over $(0,t)$, $t>0$, we obtain
\beq\label{norm(v)Lp}
\|v(t)\|_{L^p(\Omega)}\le \max\Big\{\|v^\init\|_{L^p(\Omega)},\f{\eta_v}{r_v}|\Omega|^{\f1p}\Big\},
\eeq
implying \eqref{bound v} as $p\to\infty$.
\medskip

\noindent\textit{Fourth step: maximal regularity and further estimates}.
Let $p=2$. Then, $(e^{tA})_{t\ge0}$ is a semigroup of contraction, i.e. $\|e^{tA}\Phi\|_{X_2}\le\|\Phi\|_{X_2}$, for all $\Phi\in X_2$ and all $t\ge0$ (\cite{CazenaveHaraux1998}, Theorem 3.1.1). As $-A\ge0$ and $A$ is self-adjoint, Theorem 1.3.9 in \cite{Davies1989} gives 
\[
\|e^{tA}\Phi\|_{X_p}\le\|\Phi\|_{X_p}\,,\quad\forall\ \Phi\in X_p,\ p\in[1,+\infty],\ t\ge0\,.
\]
It follows from \cite{Lamberton1987} Theorem 1 applied to the component $v$ of the strict solution $(u_a,u_b,v)$ that, for all $T>0$ and all $p\in(1,+\infty)$, there exists $C_p^\MR>0$ (not depending on $T$) such that
\beq\label{max reg v}
\|\pa_tv\|_{L^p(\Omega_T)}+\|\Delta v\|_{L^p(\Omega_T)}\le C_p^\MR(\|\Delta v^\init\|_{L^p(\Omega)}+\|f_v(u_a,u_b,v)\|_{L^p(\Omega_T)})\,.
\eeq
Next, by the classical Agmon-Douglis-Nirenberg a priori estimates (see \cite{Lunardi1995}, Theorem 3.1.1), there exists $C_p^\ADN>0$ such that, for all $t\ge0$,
\beq\label{ADN}
\|v(t)\|_{W^{2,p}(\Omega)}\le C_p^\ADN(\|v(t)\|_{L^p(\Omega)}+\|\Delta v(t)\|_{L^p(\Omega)})\,,
\eeq
while, by the Gagliardo-Nirenberg inequality \cite{Nirenberg1966},  there exists $C_p^\GN>0$ such that, for all $t\ge0$ and all $i=1,\dots,N$,
\beq\label{GN}
\|\pa_i v(t)\|_{L^{2p}(\Omega)}\le C_p^\GN\left(\max_{1\le i,j\le N}\|\pa_{i,j} v(t)\|_{L^{p}(\Omega)}\right)^{\f12}\|v(t)\|_{L^\infty(\Omega)}^{\f12}
+C_p^\GN\|v(t)\|_{L^\infty(\Omega)}\,.
\eeq
Finally, by \eqref{bound v}, there exists $C(K_\infty,|\Omega|)>0$, such that
\beq\label{est Lp fv}
\|f_v(u_a,u_b,v)\|_{L^p(\Omega_T)}\le C(K_\infty,|\Omega|)(\|u_a\|_{L^p(\Omega_T)}+\|u_b\|_{L^p(\Omega_T)}+T^{\f1p})\,.
\eeq
Hence, plugging \eqref{est Lp fv} into \eqref{max reg v}, and combining the resulting inequality with \eqref{ADN} integrated over $(0,T)$, and \eqref{bound v}, we obtain that there exists a constant $C_1(C_p^\ADN,C_p^\MR, K_\infty,|\Omega|)>0$ such that
\beq\label{max reg}
\|\partial_{t}v\|_{L^p(\Omega_T)}+\sum_{i,j}\|\partial_{ij} v\|_{L^{{p}}(\Omega_T)}\le C_1(\|\Delta v^\init\|_{L^p(\Omega)}+\| u_a\|_{L^p(\Omega_T)}+\|u_b\|_{L^p(\Omega_T)}+T^{\f1p}).
\eeq
Finally, combining the Gagliardo-Nirenberg inequality above integrated over $(0,T)$ and \eqref{max reg}, give us the existence of $C_2(C_p^\GN,C_1,K_\infty,N)$ such that
\beq\label{est gradv}
\|\nabla v\|_{L^{2p}(\Omega_T)}^{2p}\leq C_2^{2p}(\|\Delta v^\init\|^p_{L^p(\Omega)}+\| u_a\|^p_{L^p(\Omega_T)}+\|u_b\|^p_{L^p(\Omega_T)}+T)\,.
\eeq
\eqref{max reg p=2} follows  taking $p=2$ in \eqref{max reg},\eqref{est gradv} and using \eqref{result p=1}. 
\section{Proof of Lemma \ref{property cross-diff}}\label{appendix B}
By definitions \eqref{def Q},\eqref{def phi psi} and assumption \eqref{(H1)}, $\Lambda(u_a,u_b,v)\ge A^\alpha>0$, for all $(u_a,u_b,v)\in\R^3_+$, 
so that $Q(u_a,u_b,v)=0$ if and only if $q(u_a,u_b,v)=0$. Moreover, if $\tilde u=0$, the unique nonnegative solution of \eqref{nonlinearsys} is $(u_a^*,u_b^*)=(0,0)$, $\forall\ \tilde v\ge0$. 

Let us denote $\Sigma=\{(\tilde u,u_b,\tilde v)\in\R^3:0\le u_b\le \tilde u,\tilde u\ge0,\tilde v\ge0\}$ and
\[
\tilde q(\tilde u,u_b,\tilde v)=q(\tilde u-u_b,u_b,\tilde v)\,,\quad (\tilde u,u_b,\tilde v)\in\Sigma\,.
\]
By \eqref{def Q},\eqref{def phi psi} and \eqref{(H1)} again, it is easily seen that $\tilde q$ is a continuous differentiable function such that, for all $u_b\!\in\!(0,\tilde u),\tilde u>0$,$\tilde v\ge0$,
\[
\pa_2\tilde q(\tilde u,u_b,\tilde v)=-\pa_1q(\tilde u-u_b,u_b,\tilde v)+\pa_2q(\tilde u-u_b,u_b,\tilde v)\ge A^\alpha>0
\]
and $\tilde q(\tilde u,0,\tilde v)=q(\tilde u,0,\tilde v)<0$, $\tilde q(\tilde u,\tilde u,\tilde v)=q(0,\tilde u,\tilde v)>0$. Therefore,  $\forall\ \tilde u>0,\tilde v\ge0$, there exists a unique $U_b\!=\!U_b(\tilde u,\tilde v)$ such that $U_b\!\in\!(0,\tilde u)$, $\tilde q(\tilde u,U_b,\tilde v)=0$, i.e. for all $\tilde u\!>\!0$,$\tilde v\ge0$, there exists a unique solution of \eqref{nonlinearsys} given by $(u_a^*,u_b^*)=(\tilde u-U_b,U_b)$. 

As it holds $\pa_2\tilde q(U_b,\tilde u,\tilde v)>0$, by the implicit function theorem, for all $(\tilde u,\tilde v)\in(0,+\infty)^2$, there exists a neighbourhood ${\cal W}$ of $(\tilde u,\tilde v)$ and a unique continuously differentiable map $u_b^*:{\cal W}\mapsto\R_+$ such that, $\forall(\tilde u,\tilde v)\in {\cal W}$, $u_b^*(\tilde u,\tilde v)\!=U_b$ and $\tilde q(\tilde u,u_b^*(\tilde u,\tilde v),\tilde v)\!=0$. Furthermore, it is easily seen that $u_b^*$ is defined and continuously differentiable over $(0,+\infty)^2$. Hence, defining $u_b^*(0,\tilde v)=0$ and 
\[
u_a^*(\tilde u,\tilde v)=\tilde u-u_b^*(\tilde u,\tilde v),\qquad (\tilde u,\tilde v)\in[0,+\infty)^2\,,
\] 
we have that the pair $(u_a^*(\tilde u,\tilde v),u_b^*(\tilde u,\tilde v))$ is the unique solution of \eqref{nonlinearsys}.

Finally, differentiating the identities below with respect to $\tilde u$ and $\tilde v$, we obtain
\[
\tilde q(\tilde u,u_b^*(\tilde u,\tilde v),\tilde v)=q(\tilde u-u_b^*(\tilde u,\tilde v),u_b^*(\tilde u,\tilde v),\tilde v)=0\quad\text{and}\quad u_a^*(\tilde u,\tilde v)=\tilde u-u_b^*(\tilde u,\tilde v)\,.
\]
\beq\label{grad ub*}
\begin{split}
\pa_{\tilde u}u_b^*(\tilde u,\tilde v)=\f{\pa_1q(u_a^*(\tilde u,\tilde v),u_b^*(\tilde u,\tilde v),\tilde v)}{\pa_1q(u_a^*(\tilde u,\tilde v),u_b^*(\tilde u,\tilde v),\tilde v)-\pa_2q(u_a^*(\tilde u,\tilde v),u_b^*(\tilde u,\tilde v),\tilde v)}\\[1.4ex]
\pa_{\tilde v}u_b^*(\tilde u,\tilde v)=\f{\pa_3 q(u_a^*(\tilde u,\tilde v),u_b^*(\tilde u,\tilde v),\tilde v)}{\pa_1q(u_a^*(\tilde u,\tilde v),u_b^*(\tilde u,\tilde v),\tilde v)-\pa_2q(u_a^*(\tilde u,\tilde v),u_b^*(\tilde u,\tilde v),\tilde v)}\\
\end{split}
\eeq
and
\[
\pa_{\tilde u}u_a^*(\tilde u,\tilde v)=1-\pa_{\tilde u}u_b^*(\tilde u,\tilde v)\,,\qquad \pa_{\tilde v}u_a^*(\tilde u,\tilde v)=-\pa_{\tilde v}u_b^*(\tilde u,\tilde v)\,.
\]
Therefore, \eqref{partial u ua* ub*}, \eqref{partial v ua* ub*} follow taking into account that (see \eqref{def Q})
\beq\label{grad q}
\begin{split}
&\pa_1q(u_a,u_b,v)=-\psi(\bfast u_b+\dfast v)-\afast\,u_a\, \psi'(\afast u_a+\cfast v)\\
&\pa_2q(u_a,u_b,v)=\phi(\bfast u_b+\dfast v)+\bfast\,u_b\, \phi'(\bfast u_b+\dfast v)\\
&\pa_3q(u_a,u_b,v)=\dfast\,u_b\,\phi'(\bfast u_b+\dfast v)-\cfast\,u_a\,\psi'(\afast u_a+\cfast v)
\end{split}
\eeq
and the positivity of $\psi,\phi,\psi',\phi'$. 
\thispagestyle{empty}
\small
\bibliographystyle{siam}\small
\bibliography{BrocchieriCorrias2025}

\begin{thebibliography}{10}

\bibitem{Aubin1987}
{\sc J.-P. Aubin}, {\em Compact sets in the space ${L}^p(0, {T}; {B})$}, Annali
  di Matematica Pura ed Applicata (4), 146 (1987), pp.~65--96.

\bibitem{BreKueSore2021}
{\sc M.~Breden, C.~Kuehn, and C.~Soresina}, {\em On the influence of
  cross-diffusion in pattern formation}, Journal of Computational Dynamics, 8
  (2021), pp.~213--240.

\bibitem{Brocchieri2021}
{\sc E.~Brocchieri, L.~Corrias, H.~Dietert, and Y.-J. Kim}, {\em Evolution of
  dietary diversity and a starvation driven cross-diffusion system as its
  singular limit}, Journal of Mathematical Biology, 83 (2021), p.~58.

\bibitem{Brocchieri2024}
{\sc E.~Brocchieri, L.~Desvillettes, and H.~Dietert}, {\em Study of a class of
  triangular starvation driven cross-diffusion systems}, Ricerche di
  Matematica,  (2024), pp.~1--27.

\bibitem{CazenaveHaraux1998}
{\sc T.~Cazenave and A.~Haraux}, {\em An Introduction to Semilinear Evolution
  Equations}, vol.~13 of Oxford Lecture Series in Mathematics and Its
  Applications, The Clarendon Press, Oxford University Press, New York, 1998.

\bibitem{ChenJuengel2004}
{\sc L.~Chen and A.~J\"ungel}, {\em Analysis of a multidimensional parabolic
  population model with strong cross-diffusion}, SIAM Journal on Mathematical
  Analysis, 36 (2004), pp.~301--322.

\bibitem{ChenJuengel2006}
\leavevmode\vrule height 2pt depth -1.6pt width 23pt, {\em Analysis of a
  parabolic cross-diffusion population model without self-diffusion}, Journal
  of Differential Equations, 224 (2006), pp.~39--59.

\bibitem{Conforto2018}
{\sc F.~Conforto, L.~Desvillettes, and C.~Soresina}, {\em About
  reaction-diffusion systems involving the holling-type ii and the
  beddington--deangelis functional responses for predator-prey models}, NoDEA
  Nonlinear Differential Equations Appl., 25 (2018), pp.~Paper No. 24, 39.

\bibitem{Daus2019}
{\sc E.~S. Daus, L.~Desvillettes, and H.~Dietert}, {\em About the entropic
  structure of detailed balanced multi-species cross-diffusion equations},
  Journal of Differential Equations, 266 (2019), pp.~3861--3882.

\bibitem{Davies1989}
{\sc E.~B. Davies}, {\em Heat Kernels and Spectral Theory}, vol.~92 of
  Cambridge Tracts in Mathematics, Cambridge University Press, Cambridge, 1989.

\bibitem{Desvillettes2025}
{\sc L.~Desvillettes, C.~Kuehn, J.-E. Sulzbach, B.~Q. Tang, and B.-N. Tran},
  {\em Slow manifolds for {PDE} with fast reactions and small cross diffusion},
  arXiv preprint arXiv:2501.16775,  (2025).

\bibitem{Desvillettes2014}
{\sc L.~Desvillettes, T.~Lepoutre, and A.~Moussa}, {\em Entropy, duality, and
  cross diffusion}, SIAM Journal on Mathematical Analysis, 46 (2014),
  pp.~820--853.

\bibitem{Desvillettes2015}
{\sc L.~Desvillettes, T.~Lepoutre, A.~Moussa, and A.~Trescases}, {\em On the
  entropic structure of reaction--cross diffusion systems}, Communications in
  Partial Differential Equations, 40 (2015), pp.~1705--1747.

\bibitem{Desvillettes2015a}
{\sc L.~Desvillettes and A.~Trescases}, {\em New results for triangular
  reaction cross diffusion system}, Journal of Mathematical Analysis and
  Applications, 430 (2015), pp.~32--59.

\bibitem{GAMBINO2012}
{\sc G.~Gambino, M.~C. Lombardo, and M.~Sammartino}, {\em Turing instability
  and traveling fronts for a nonlinear reaction--diffusion system with
  cross-diffusion}, Mathematics and Computers in Simulation, 82 (2012),
  pp.~1112--1132.

\bibitem{Haraux2017}
{\sc A.~Haraux}, {\em A simple characterization of positivity preserving
  semi-linear parabolic systems}, Journal of the Korean Mathematical Society,
  54 (2017), pp.~1817--1828.

\bibitem{Iida2023}
{\sc M.~Iida, H.~Izuhara, and R.~Kon}, {\em Cross-diffusion predator-prey model
  derived from the dichotomy between two behavioral predator states}, Discrete
  Contin. Dyn. Syst. Ser. B, 28 (2023), pp.~6159--6178.

\bibitem{Iida2006}
{\sc M.~Iida, M.~Mimura, and H.~Ninomiya}, {\em Diffusion, cross-diffusion and
  competitive interaction}, Journal of Mathematical Biology, 53 (2006),
  pp.~617--641.

\bibitem{Lamberton1987}
{\sc D.~Lamberton}, {\em \'equations d'\'evolution lin\'eaires associ\'ees \`a
  des semi-groupes de contractions dans les espaces {$L^p$}}, Journal of
  Functional Analysis, 72 (1987), pp.~252--262.

\bibitem{Lep2017}
{\sc T.~Lepoutre and A.~Moussa}, {\em Entropic structure and duality for
  multiple species cross-diffusion systems}, Nonlinear Analysis, 159 (2017),
  pp.~298--315.

\bibitem{LionsMagenes1972}
{\sc J.-L. Lions and E.~Magenes}, {\em Non-Homogeneous Boundary Value Problems
  and Applications, Volume II}, vol.~182 of Die Grundlehren der Mathematischen
  Wissenschaften, Springer-Verlag, New York--Heidelberg, 1972.
\newblock Translated from the French by P. Kenneth.

\bibitem{Lunardi1995}
{\sc A.~Lunardi}, {\em Analytic Semigroups and Optimal Regularity in Parabolic
  Problems}, vol.~16 of Progress in Nonlinear Differential Equations and Their
  Applications, Birkh{\"a}user, Basel, 1995.

\bibitem{Moussa2019}
{\sc A.~Moussa, B.~Perthame, and D.~Salort}, {\em Backward parabolicity,
  cross-diffusion and turing instability}, Journal of Nonlinear Science, 29
  (2019), pp.~139--162.

\bibitem{Nirenberg1966}
{\sc L.~Nirenberg}, {\em An extended interpolation inequality}, Annali della
  Scuola Normale Superiore di Pisa - Classe di Scienze, 20 (1966),
  pp.~733--737.

\bibitem{Shigesada1979}
{\sc N.~Shigesada, K.~Kawasaki, and E.~Teramoto}, {\em Spatial segregation of
  interacting species}, Journal of Theoretical Biology, 79 (1979), pp.~83--99.

\bibitem{Soresina2023}
{\sc C.~Soresina, Q.~B. Tang, and B.~N. Tran}, {\em Fast-reaction limits for
  predator--prey reaction--diffusion systems: improved convergence}, arXiv
  preprint,  (2023).

\bibitem{BaoBao2024}
{\sc B.~Q. Tang and B.-N. Tran}, {\em Rigorous derivation of michaelis--menten
  kinetics in the presence of slow diffusion}, SIAM Journal on Mathematical
  Analysis, 56 (2024), pp.~5995--6024.
\newblock Published by SIAM.

\bibitem{Vazquez1992}
{\sc J.~L. V{\'a}zquez}, {\em An introduction to the mathematical theory of the
  porous medium equation}, in Shape Optimization and Free Boundaries (Montreal,
  PQ, 1990), vol.~380 of NATO Advanced Science Institutes Series C:
  Mathematical and Physical Sciences, Kluwer Academic Publishers, Dordrecht,
  1992, pp.~347--389.

\end{thebibliography}
\noindent
Email addresses: \\
Elisabetta Brocchieri: elisabetta.brocchieri@uni-graz.at\\
Lucilla Corrias: lucilla.corrias@univ-evry.fr
\medskip

\noindent ${^1}$
Department of Mathematics and Scientific Computing, University of Graz, 8010 Graz, Austria.\\
\noindent ${^2}$
Universit\'e Paris-Saclay, CNRS, Univ. Evry, Laboratoire de Math\'ematiques et Mod\'elisation d'Evry, 91037, Evry-Courcouronnes, France. 
%
\end{document}